\documentclass[11pt]{article}
\usepackage[paper=a4paper,
	top=3cm,
	inner=2.54cm,
	outer=2.54cm,
	bottom=3cm,
	headheight=5ex,
	headsep=5ex]{geometry}
\usepackage[pagewise]{lineno}
\usepackage[bookmarks,colorlinks ]{hyperref}
\hypersetup{linkcolor =blue,
	filecolor=red}
\usepackage{enumitem}
\usepackage{mathdots}
\usepackage{tikz-cd}

\usepackage{amsthm}
\newtheorem{thm}{Theorem}[section]
\newtheorem{prop}[thm]{Proposition}
\newtheorem{lem}[thm]{Lemma}
\newtheorem{cor}[thm]{Corollary}

\theoremstyle{definition}

\newtheorem{fact}{Fact}

\theoremstyle{remark}
\newtheorem{rmk}{Remark}[subsection]
\newtheorem{eg}{Example}

\newenvironment{pf}[1][]{{\noindent\it{Proof#1}}. \small}{\hfill \qed
	\vspace{0.5em}}

\usepackage{amsmath, amssymb}
\numberwithin{equation}{section}
\usepackage{mathrsfs}
\newcommand{\N}{\mathbb N}
\newcommand{\Z}{\mathbb Z}
\newcommand{\R}{{\mathbb R}}
\newcommand{\C}{{\mathbb C}}
\renewcommand{\det}{{\mathrm{det}}}
\newcommand{\Tr}{{\mathrm{Tr}}}
\newcommand{\Hom}{\operatorname{Hom}}
\newcommand{\End}{\operatorname{End}}
\newcommand{\Id}{\operatorname{Id}}

\renewcommand{\span}{\operatorname{span}}
\newcommand{\GL}{{\mathrm {GL}}}
\newcommand{\SL}{{\mathrm {SL}}}
\renewcommand{\O}{{\mathrm O}}
\newcommand{\SO}{{\mathrm {SO}}}
\newcommand{\Sp}{{\mathrm {Sp}}}
\newcommand{\la}{\left\langle}
\newcommand{\ra}{\right\rangle}

\newcommand{\lag}[1]{{\mathfrak {#1}}}
\newcommand{\uea}[1]{\mathcal U (#1)}
\newcommand{\Irr}{\operatorname{Irr}}
\newcommand{\Ind}{\operatorname{Ind}}

\title{On the twisted Osborne conjecture}
\author{Chang Huang}
\usepackage{natbib}
\begin{document}
\maketitle
\begin{abstract}
	We aim to prove a twisted version of the Osborne conjecture.
	The untwisted case was proved by Hecht and Schmid in their 1983 Acta Mathematica paper.
	Bergeron and Clozel (2013) have considered a special case, and we generalize their method to our setting.
\end{abstract}

\tableofcontents

\newcommand{\ag}[1]{\mathbf {#1}}
\newcommand{\spl}{\mathbf{spl}}
\newcommand{\Cont}{{\mathrm C}}
\newcommand{\Smt}{\Cont^\infty}
\newcommand{\cpt}{\mathrm c}
\renewcommand{\d}{{\, \mathrm d}}
\newcommand{\Ad}{\operatorname{Ad}}
\newcommand{\Tp}{{\mathrm T}}
\newcommand{\SC}[1]{\mathcal D^{\mathrm{st}}(#1)}
\newcommand{\TC}[1]{\mathcal D^\tau(#1)}
\newcommand{\Tran}{\operatorname{Tran}}
\newcommand{\D}[2][]{\operatorname D_{#1}^{#2}}
\newcommand{\sgn}{\operatorname {sgn}}
\newcommand{\Jac}{\operatorname{Jac}}
\renewcommand{\ss}{{\mathrm{ss}}}

\section{Introduction}
	Let $\ag G$ be a connected real reductive group and $\tau$ be a semisimple $\R$-automorphism of $\ag G$.
	According to \cite[Theorem 7.5]{St68}, $\tau$ preserves a Borel pair $(\ag B, \ag T)$ of $\ag G$.
	We assume that $\tau$ preserves a splitting $\spl = (\ag B, \ag T, \{X_\alpha\})$, and has finite order $d$.

	Suppose $\pi$ is a $\tau$-stable irreducible Casselman-Wallach representation of $G: = \ag G(\R)$.
	Fix an intertwining operator $\tau_\pi: \pi \to \pi^\tau$ such that $\tau_\pi^d = \Id$;  such an operator is unique up to multiplication by a $d$-th root of unity.
	Then we can extend $\pi$ to a representation $\pi^+$ of $G^+:= G \rtimes \la \tau \ra$ by setting $\pi^+(\tau) = \tau_\pi$.
	Consider the distribution $\Smt_\cpt(G^+) \to \C$, $f \mapsto \Tr \int_{G^+} f(g)\pi^+(g) \d g$, which gives rise to a well-defined $G$-invariant generalized function $\Theta(\pi^+)$ with $Z({\lag g})$ acting by scalars.
	By \cite[Theorem 2.1.1]{Bouaziz87}, $\Theta(\pi^+)$ is a locally integrable function and is analytic around regular points.
	In particular, we can define the twisted character of $\pi$ by restriction: $\Theta^\tau(\pi)(g) = \Theta(\pi^+)(g\tau)$.

	Let $K$ be a $\tau$-stable maximal compact subgroup of $G$, and let $P = MN$ be a $\tau$-stable Levi decomposition of a real parabolic subgroup of $G$.
	Denote the underlying $({\lag g}, K)$-module of $\pi$ by $V$.
	Then the ${\lag n}$-homology groups $H_q({\lag n}, V)$ of $V$ are $\tau$-stable Harish-Chandra modules of $M$, where $\lag n$ denotes the complex Lie algebra of $N$.
	It is well-known that the category of Harish-Chandra modules is equivalent to the category of Casselman-Wallach representations (see \cite{Casselman89}, \cite[Section 11]{Wallach92} and \cite{BK14}).
	Thus, those $H_q(\lag n, V)$ can be globalized to Casselman-Wallach representations of $M$, and we denote them by $H_q(\lag n, \pi)$; note that this notation is not standard.
	These $\lag n$-homology groups are also $\tau$-stable, and have twisted characters $\Theta^\tau(H_q(\lag n, \pi))$.

	The aim of this article is to verify the following character identity:
\begin{thm}\label{thm: main}
	There exists a sufficiently large subset $M^- \subset M$ such that
\begin{equation}\label{eq: twisted-Osborne}
	\Theta^\tau (\pi)|_{M^-} = 
	\left. \frac{1}{D^\tau_{\lag n}} \sum_{q \geqslant 0} (-1)^q \Theta^\tau \left(H_q(\lag n, \pi) \right) \right |_{M^-},
	\end{equation}
	where $D^\tau_{\lag n}(m) = \det_{\lag n} (1 - m\tau)$, and the subscript $\lag n$ means that the determinant is taken of the endomorphism $\Id_{\lag n} - \Ad(m) \circ \tau$ on ${\lag n}$.
\end{thm}
	Indeed, the subset $M^-$ consists of elements $m \in M$ such that
\begin{itemize}
\item
	$m$ is $\tau$-regular in $G$, and
\item
	all eigenvalues of $\Ad(m)\circ \tau$ on ${\lag n}$ have modulus less than 1.
\end{itemize}
	The ``sufficiently large'' property is made precise in Corollary \ref{cor: big-enough}: let $A$ be the central split torus of $M$; then, for every $\tau$-regular element $m\in M$, there exists an open subset $A_m^-\subset A$ such that $mA_m^-\subseteq M^-$.

\subsection{Application}
	In a forthcoming paper \cite{DHSX}, we will study theta lifts of unipotent Arthur packets of $\Sp_{2m}(\R)$ and $\O(p, q)$, and give an inductive construction of these packets.
	When constructing unipotent Arthur packets of $\Sp_{2m}(\R)$ from those of $\O(p, q)$ with $p+q<2m$ (stated in \cite[Theorem 1.2(a)]{DHSX}), a key technical step is to prove that a certain isotypic component of $H_0(\lag n, -)$ preserves the unipotent Arthur packets.
	For this result, we primarily use the main result, Theorem \ref{thm: main}, of this article.
	The following is a more detailed explanation.

	A unipotent Arthur parameter of $\Sp_{2m}(\R)$ can be regarded as a semisimple representation $\psi: W_\R \times \SL_2(\C) \to \GL_{2m+1}(\C)$ that decomposes into
\[\bigoplus_{i=1}^r \sgn^{\epsilon_i} \boxtimes S_{m_i},\]
	where $W_\R$ is the Weil group of $\R$, $\epsilon_i \in \{0, 1\}$, $\sgn:W_\R \to \{\pm 1\}$ is its (unique) non-trivial quadratic character, $m_1 \geqslant \cdots \geqslant m_r>0$ are odd integers, and $S_{m_i}$ is the (unique) $m_i$-dimensional irreducible representation of $\SL_2( \C)$.

Using the twisted endoscopic transfer $\Tran_\Sp^\tau: \SC{\Sp_{2m}(\R)} \to \TC {\GL_{2m+1}(\R)}$, Arthur associated a stable character $\eta_\psi \in \SC{\Sp_{2m}(\R)}$ with each Arthur parameter $\psi$ of $\Sp_{2m}(\R)$ by specifying its image under $\Tran_\Sp^\tau$ (\cite[Theorem 2.2.1]{Arthur13}).
	Here $\SC{\Sp_{2m}(\R)}$ consists of the stable characters of ${\Sp_{2m}( \R)}$, and $\TC {\GL_{2m+1}(\R)}$ consists of the twisted characters of ${\GL_{2m+1}( \R)}$.
	The automorphism $\tau$ of $\GL_{2m+1}(\R)$ is given by $g \mapsto J_{2m+1} g^{-\Tp} J_{2m+1}^{-1}$, where $g^{-\Tp}$ is the inverse transpose of $g$, and we write
\[J_N=
	\begin{pmatrix}
		&	&1\\
		&\iddots	&	\\
	(-1)^{N-1}	&	&	
	\end{pmatrix}\]
	for each integer $N\geqslant 1$.
	The Arthur packet $\Pi_\psi(\Sp_{2m}(\R))$ consists of the irreducible representations whose characters occur in $\eta_\psi$.
	When $\psi$ is unipotent, we call $\Pi_\psi(\Sp_{2m}(\R))$ a unipotent Arthur packet of $\Sp_{2m}(\R)$.

	For the technical proofs in \cite{DHSX}, we need the compatibility between twisted endoscopic transfer and $\lag n$-homology.
	Let $P= MN$ be the standard parabolic subgroup of $\Sp_{2m}(\R)$, with Levi component $M = \GL_1(\R) \times \Sp_{2m-2}(\R)$, and let $Q= LU$ be the $\tau$-stable standard parabolic subgroup of $\GL_{2m+1}(\R)$, with Levi subgroup $L = \GL_{1}(\R)\times \GL_{2m-1}(\R) \times \GL_{1}(\R)$.
	Then the Levi subgroup $M$ is an endoscopic group of $(L, \tau)$, and we also have the twisted endoscopic transfer $\Tran_M^\tau: \SC M \to \TC L$.

\begin{thm}\label{thm: compatibility}
	The following diagram commutes:
\begin{equation}\begin{tikzcd}\label{dig: compatibility}
	\SC{\Sp_{2m}( \R)} 
		\arrow[rr, "\Tran_{\Sp}^\tau"] 
		\arrow[d, "{\D[{\lag n}]{}} "'] & 
	&	
	\TC{\GL_{2m+1}( \R)}
		\arrow[d, "{\D[{\lag u}]{}}"] \\
	\SC{M} 
		\arrow[rr, "\Tran_{M}^\tau"]      & 
	&
	\TC{L}
	\end{tikzcd},\end{equation}
	where $\D[{\lag n}]{}$, $\D[{\lag u}]{}$ are induced by $|\det_{\lag n}|^{-\frac 1 2} \otimes \sum_q (-1)^q H_q({\lag n}, -)$, $|\det_{\lag u}|^{-\frac 1 2 } \otimes {\sum_q (-1)^q H_q({\lag u}, -)}$ respectively.
\end{thm}

	The preceding commutative diagram also appears in \cite[Appendix C]{DHSX}.
	Theorem \ref{thm: main} enters the proof of Theorem \ref{thm: compatibility} in two ways:
\begin{itemize}
\item
	the homomorphisms $\D[{\lag n}]{}$, $\D[{\lag u}]{}$ are well-defined since \eqref{eq: twisted-Osborne} holds on a sufficiently large subset, although the twisted character does not necessarily distinguish representations;
\item
	the commutativity of diagram \eqref{dig: compatibility} is deduced from the explicit formulas for $\Tran^\tau$ and $\D{}$.
\end{itemize}

	The following corollary establishes the connections between unipotent Arthur packets via $\lag n$-homology.
	Note that $M = \GL_1(\R) \times \Sp_{2m-2}(\R)$ has a central split torus $A \cong \R_{>0}$.
	For a Harish-Chandra module $V$ of $\Sp_{2m}(\R)$, we define $\D[0] x (V)$ to be the generalized $| \cdot |^x$-isotypic component of $H_0(\lag n, V) \otimes | \det_{\lag n}|^{-\frac 1 2}$ with respect to the action of $A$.
	Then $\D[0] x(V)$ is a Harish-Chandra module of $\Sp_{2m-2}(\R)$.
	Via Casselman-Wallach globalization, this functor induces a functor on Casselman-Wallach representations, denoted by the same symbol $\D[0] x$.

\begin{cor}\label{cor: unipotent}
	If $m_1> m_2$, then for each $\sigma \in \Pi_\psi(\Sp_{2m}(\R))$, all irreducible factors of $\D[0] {-\frac{m_1 -1 }2 }(\sigma)$ appear in $\Pi_{\psi_{-}}(\Sp_{2m-2}(\R))$, and all representations in $\Pi_{\psi_-}(\Sp_{2m-2}(\R))$ arise in this manner.
	Here, $\psi_{-}$ is the unipotent Arthur parameter of $\Sp_{2m-2}(\R)$ which decomposes into
\[\sgn^{\epsilon_1} \boxtimes S_{m_1 -2} \oplus \bigoplus_{i =2}^r \sgn^{\epsilon_i} \boxtimes S_{m_i}\]
	as a representation of $W_\R \times \SL_2(\C)$.
\end{cor}
This result is stated in \cite[Lemma C.1]{DHSX} and proved throughout Appendix C of that paper.
	The main point is to pass to the isotypic component of the diagram \eqref{dig: compatibility} and then derive $\D[0]{-\frac{m_1-1}{2}}\eta_\psi = \eta_{\psi_-}$ from the definition of $\eta_{\psi_-}$.

\subsection{$p$-adic analogy}
	The above story has a $p$-adic analogue, as presented in \cite[Appendix C]{Xu17cusp}.

	In this subsection, let $\ag G$ be a connected quasi-split reductive group over a $p$-adic field $F$ and $\tau$ be a finite-order $F$-automorphism of $\ag G$ preserving an $F$-splitting.
	Denote the order of $\tau$ by $d$.
	Suppose $\pi$ is a $\tau$-stable irreducible smooth representation of $G: = \ag G(F)$.
	Fix an intertwining operator $\tau_\pi: \pi \to \pi^\tau$ with $\tau_\pi^d = \Id$. There is then a twisted character $\Theta^\tau(\pi)$ of $\pi$ defined in \cite{Clozel87}.
	Note that in \cite{Xu17cusp} it is denoted by $\Theta^{G^\tau}_\pi$ and we omit the group $G$ since it is clear from $\pi$.
	The $p$-adic analogue of Theorem \ref{thm: main} is
\begin{thm}[Casselman, Rogawski]
	Let $P=  M N$ be a $\tau$-stable parabolic subgroup of $G$.
	Then for each $m \in M$ that is $\tau$-regular in $G$, there exists $a \in A$ such that
\begin{equation}\label{eq: char-of-Jac}
	\Theta^\tau(\pi)(ma) = \Theta^\tau(\pi_N)(ma).
	\end{equation}
	Here $A$ is the maximal split component of the center of $M$, and $\pi_N$ is the unnormalized Jacquet module of $\pi$ with respect to $P$.
\end{thm}
	This formulation is taken from \cite[Corollary C.6]{Xu17cusp}, which is a corollary of \cite[Theorem 5.2]{Casselman77} when $\tau = \Id$, and of \cite[Proposition 7.4]{Rogawski88} in the twisted case.

\begin{rmk}
	There is flexibility in choosing $a$: there exists $\epsilon_m>0$ such that \eqref{eq: char-of-Jac} holds for $ma$ whenever $|\beta(a)|< \epsilon_m$ for every root $\beta$ in $\ag N$.
	We will see from Corollary \ref{cor: big-enough} that a similar phenomenon also occurs for \eqref{eq: twisted-Osborne}.
\end{rmk}
	There are two differences between \eqref{eq: twisted-Osborne} and \eqref{eq: char-of-Jac}:
\begin{itemize}
\item
	the real analogue of the unnormalized Jacquet module $\pi_N$ seems to be $H_0({\lag n}, -)$, but the latter is not an exact functor, so we have to consider its derived functors, and that is why the Euler alternating sum occurs in \eqref{eq: twisted-Osborne};
\item
	when $|\beta(a)|$ is small enough for every root $\beta$ in $\ag N$, all eigenvalues of $\Id_{\lag n} - \Ad(ma) \circ \tau$ have modulus $1$, which explains why no term like $D^\tau_{\lag n}(ma)$ occurs in \eqref{eq: char-of-Jac}.
\end{itemize}

	As in the Archimedean case, we also have the twisted endoscopic transfer $\Tran_\Sp^\tau: \SC{\Sp_{2m}(F)} \to \TC {\GL_{2m+1}(F)}$, which Arthur uses to characterize the so-called Arthur packet of $\Sp_{2m}(F)$.
	The $p$-adic analogue of our Theorem \ref{thm: compatibility} is
\begin{thm}[Xu]
	Let $P = MN$ be the standard parabolic subgroup of $\Sp_{2m}( F)$ with Levi subgroup $M = \GL_c(F) \times \Sp_{2(m-c)}(F)$, and  let $Q= LU$ be the $\tau$-stable standard parabolic subgroup of $\GL_{2m+1}(F)$ with Levi subgroup $L = \GL_c(F) \times \GL_{2(m-c)+1}(F) \times \GL_c (F)$.
	Then the following diagram commutes:
\[\begin{tikzcd}
	\SC{\Sp_{2m}( F)} 
		\arrow[rr, "\Tran_{\Sp}^\tau"] 
		\arrow[d, "\Jac_P"'] & 
	&	
	\TC{\GL_{2m+1}(F)}
		\arrow[d, "{\Jac_Q}"] \\
	\SC{M} 
		\arrow[rr, "\Tran_{M}^\tau"]      & 
	&
	\TC{L}
	\end{tikzcd},\]
	where $\Jac_P(-)$ and $\Jac_Q(-)$ are induced by the normalized Jacquet modules $(-)_N \otimes \delta_P^{-\frac 1 2}$ and $(-)_U \otimes \delta_Q^{-\frac 1 2}$, respectively.
\end{thm}
	This result can be deduced from \cite[diagram (C.5)]{Xu17cusp}, which asserts the compatibility between twisted endoscopic transfer and normalized Jacquet functor in a slightly more general case.
	The ordinary (non-twisted) endoscopic transfer is also compatible with the normalized Jacquet functor, as proved in \cite[Theorem 5.6]{Hiraga04}.
	The two cases share an analogous argument.

	Note that we also use the normalized version $\D[{\lag n}]{}$ in Theorem \ref{thm: compatibility}.

	In the $p$-adic case, the Arthur parameter of $\Sp_{2m}(F)$ can be regarded as a semisimple representation $\psi: W_F \times \SL_2(\C) \times \SL_2(\C) \to \GL_{2m+1}(\C)$ of orthogonal type (i.e., the image of $\psi$ is contained in $\SO_{2m+1}(\C)$), where $W_F$ is the Weil group of $F$.
	Here we only consider the \textit{anti-tempered and good parity} case,  where $\psi$ decomposes into
\[\bigoplus_{i=1}^r \rho_i \boxtimes S_1 \boxtimes S_{m_i},\]
	and each irreducible factor $\rho_i \boxtimes S_1 \boxtimes S_{m_i}$ is of orthogonal type.

	According to the local Langlands correspondence for general linear groups, each $c$-dimensional irreducible representation $\rho$ of $W_F$ corresponds to a unique irreducible supercuspidal representation of $\GL_c(F)$, which we also denote by $\rho$.
For any $x\in \R$, we define $\Jac_x$ to be the generalized $\rho \otimes | \det |^x$-isotypic component of $\Jac_P$ with respect to the action of $\GL_c(F) \subset M$. This functor sends smooth representations of $\Sp_{2m}(F)$ to smooth representations of $\Sp_{2(m-c)}(F)$.
Now we can state the following $p$-adic analogue of Corollary \ref{cor: unipotent}, where we assume $\psi$ is \textit{multiplicity-free} for convenience:
\begin{cor}\label{thm: anti-temper}
	If $\rho \boxtimes S_1 \boxtimes S_m$ occurs in $\psi$, then for each $\pi \in \Pi_\psi(\Sp_{2m}(F))$, all irreducible factors of $\Jac_{-\frac{m-1}2}(\pi)$ appear in $\Pi_{\psi_-}(\Sp_{2(m-c)}(F))$, and all representations in $\Pi_{\psi_-}(\Sp_{2(m-c)}(F))$ arise in this way.
	Here $\psi_-$ is the Arthur parameter of $\Sp_{2(m-c)}(F)$ obtained from $\psi$ by substituting the factor $\rho \boxtimes S_1 \boxtimes S_m$ with $\rho \boxtimes S_1 \boxtimes S_{m-2}$.
\end{cor}
\begin{pf}
	This result appears to be standard; we include a brief argument.
	Let $\widehat \psi: W_F \times \SL_2( \C) \times \SL_2(\C) \to \SO_{2m+1}(\C)$ be the tempered Arthur parameter with decomposition $\bigoplus_{i=1}^r \rho_i \boxtimes S_{m_i} \boxtimes S_1$.
	It follows from \cite[Lemma 7.2]{Xu17cusp} that $\Jac_{\frac{m-1} 2}$ sends $\Pi_{\widehat \psi}$ to $\Pi_{\widehat \psi_-}$.
	The above result can be deduced using the Aubert involution, as discussed in \cite[Appendix B]{AGIKMS}.
	This involution sends $\Pi_{\widehat \psi}$ to $\Pi_{\psi}$ according to \cite[Theorem 5.9]{LLS24}, and intertwines $\Jac_{\frac{m-1}2}$ with $\Jac_{-\frac{m-1} 2}$ according to \cite[Theorem B.2.3]{AGIKMS}.
	Then the conclusion follows.
\end{pf}

	In Corollary \ref{cor: unipotent}, $\D[0]{-\frac{m_1-1} 2}$ plays a role similar to that of $\Jac_{-\frac{m-1} 2}$ in the above theorem.
	The higher homologies $\D[q]{-\frac{m_1-1}2}$ disappear due to \cite[Proposition 2.32]{HS83}.

	In summary, on both the real and the $p$-adic sides, certain isotypic components of normalized ${\lag n}$-homology and of the Jacquet functor preserve particular types of Arthur packets, and thereby contribute to the explicit construction of Arthur packets; the two stories start from the character formulas \eqref{eq: twisted-Osborne} and \eqref{eq: char-of-Jac}.
	It would be very interesting to consider whether we can construct a comparison functor between representation categories of real and $p$-adic groups, which provides a bridge between the above two parallel stories.

\subsection{Existing results}
	Two special cases of Theorem \ref{thm: main} have been dealt with in earlier work.

	When $\tau$ is trivial, \eqref{eq: twisted-Osborne} reduces to the ordinary Osborne conjecture, which was proposed in \cite{Osborne72} and proved in \cite[Theorem 3.6]{HS83}.
	The proof is based on several fundamental results in the representation theory of real reductive groups: coherent continuation, Langlands classification, Hirai's induced character formula, and Harish-Chandra's estimates of tempered characters.

	When $P$ is a minimal parabolic subgroup of $G$, \cite[Theorem 4.5]{BC13} establishes that \eqref{eq: twisted-Osborne} holds on $T^{\tau, \circ} \cap M^-$, where $T \subseteq M$ is a $\tau$-stable maximally split Cartan subgroup of $G$, and $T^{\tau, \circ}$ is the identity component of the $\tau$-fixed points $T^\tau \subset T$.
	Its proof is adapted from that in \cite{HS83}.

	As mentioned in \cite[footnote (4) on page 128]{BC13}, there may be no $\tau$-stable minimal parabolic subgroups of $G$ in general.
	Nevertheless, the arguments in \cite{HS83} and \cite{BC13} can be extended to prove Theorem \ref{thm: main}, as we do in this article.

\subsection{Organization}
	In Section \ref{sec: preparation} we establish some results about real reductive groups.
	In Section \ref{sec: reduction} we reduce Theorem \ref{thm: main} to Propositions \ref{prop: vanish} and \ref{prop: induce} by the technique of coherent continuation as in \cite[Chapter 5, Section 6.1]{BC13}.
	In Section \ref{sec: ind-character} we calculate the regular value of twisted characters of a parabolically induced module, which leads to a proof of Proposition \ref{prop: induce}.
	In Section \ref{sec: asymtotic} we estimate the twisted character of a representation based on its Langlands classification, which leads to a proof of Proposition \ref{prop: vanish}.

\subsection*{Acknowledgement}
	The author thanks Bin Xu for suggesting that the current formulation of Corollary \ref{thm: anti-temper} is a more appropriate analogue of Corollary \ref{cor: unipotent} than \cite[Lemma 7.2]{Xu17cusp}.
 
\newcommand{\RS}{R}
\newcommand{\Av}{\operatorname{Av}_\tau}
\newcommand{\res}{{\mathrm{res}}}
\newcommand{\der}{{\mathrm{der}}}
\newcommand{\ad}{\operatorname{ad}}
\renewcommand{\min}{{\mathrm{min}}}
\renewcommand{\Re}{\operatorname{Re}}
\renewcommand{\Im}{\operatorname{Im}}
\newcommand{\pr}{\operatorname{pr}}

\section{Preliminaries}\label{sec: preparation}
	The aim of this section is two-fold.
	First, we introduce some notions and fix the notation that will be used throughout this article; for example, admissible Cartan subgroup and twisted Weyl denominator.
	Second, we prove some properties which are well-known to experts, namely the existence of a $\tau$-stable maximal compact subgroup and the $\tau$-fixed representative for a $\tau$-stable infinitesimal character.

\subsection{Notation}\label{subsec: notation}
	We denote algebraic groups, groups of real points, real Lie algebras, and complex Lie algebras by boldface capitals, italic capitals, lowercase fraktur letters with subscript $0$, and lowercase fraktur letters, respectively;
	for example, $\ag G$, $G$, ${\lag g}_0$, ${\lag g}$, and $\ag T_s$, $T_s$, ${\lag t}_{s,0}$, ${\lag t}_s$.
	The algebraic group  $\ag G$ can be identified with its group of $\C$-points $\ag G(\C)$, which is a complex Lie group, while $G$ is viewed as a real Lie group.
	By a superscript $^\circ$ we mean taking the identity component (of an algebraic group or Lie group).
	For example, $\ag G^{\tau, \circ}$ denotes the identity component of $\ag G^\tau$ (either an algebraic group or a complex Lie group), and $T^{\tau, \circ}$ denotes the identity component of $T^\tau$ (a real Lie group).

	Let $H$ be a subgroup of $G$ and act on $G^+$ by conjugation.
  For any subset $S \subset G^+$, we denote the centralizer and normalizer of $S$ in $H$ by $Z(H, S) := \{ g\in H \mid gs = sg, \forall s\in S \}$ and $N(H, S) := \{ g \in H \mid \Ad(g) S = S\}$ respectively.
	The same notation is used at the level of Lie algebras, with one exception: $Z(\lag g)$ denotes the center of the universal enveloping algebra $\uea{\lag g}$ of $\lag g$.

	Let ${\lag h} \subseteq {\lag g}$ be contained in some Cartan subalgebra and $\lag l \subseteq {\lag g}$ be an ${\lag h}$-stable subspace.
	Then $\lag l$ decomposes into weight spaces with respect to ${\lag h}$.
	Denote the set of non-zero weights (with multiplicities) by $\RS(\lag l, {\lag h})$.
	In particular, if ${\lag t} \subseteq {\lag g}$ is a Cartan subalgebra, then $\RS({\lag g}, {\lag t})$ is exactly the root system determined by $({\lag g}, {\lag t})$.
	If ${\lag b} \subseteq {\lag g}$ is a Borel subalgebra containing ${\lag t}$, then it defines a positive root system and simple roots in $\RS({\lag g}, {\lag t})$.
	They are denoted by $\RS^+({\lag g}, {\lag t})$ and $\Pi({\lag g}, {\lag t})$ respectively.

	Denote half the sum of $\RS(\lag l, {\lag h})$ by $\rho(\lag l, {\lag h})$.
	In particular, if $P \subseteq G$ is a real parabolic subgroup containing $T$, then we abbreviate $\rho({\lag p}, {\lag t})$ to $\rho_P$ when ${\lag t}$ is clear from the context.
	It also equals $\rho({\lag n}, {\lag t})$, where ${\lag n}$ is the radical of ${\lag p}$.

	The Borel pair $(\ag B, \ag T)$ also determines a positive root system $\RS^+(\ag G, \ag T) \subset X^*(\ag T)$.
	Roots in $\RS^+({\lag g}, {\lag t})$ are differentials of roots in $\RS^+(\ag G, \ag T)$, and we can identify them.
	For a root $\alpha$, we write $\alpha(t)$ when regarding it as a character of $\ag T$, and write $\la \alpha, X \ra$ when regarding it as a functional on ${\lag t}$.

	If the adjoint action of $g \in G^+$ on ${\lag g}$ induces an endomorphism $\Ad(g)$ on a subquotient space ${\lag f}$, then we denote the determinant of this endomorphism by $\det_{\lag f}(g)$.
	Moreover, let $D_{\lag f}(g) = \det_{\lag f}( 1 - g)$ denote the determinant of the endomorphism $ \Id_{\lag f} - \Ad(g)$ on ${\lag f}$.
	If $g = ts\tau$, and $\Ad(t)$, $\Ad(s\tau)$ both act on ${\lag f}$, then we also write $D^{s\tau}_{\lag f}(t)$ in place of $D_{\lag f}(ts\tau)$.

\subsection{Admissible Cartan subgroup}
	Let $s\tau \in G^+$ be semisimple, so that $\Ad(s) \circ \tau$ acts semisimply on $\lag g_0$.
	By \cite[Theorem 7.5]{St68}, there exists an $s\tau$-admissible Cartan subgroup $\ag T_s \subseteq \ag G$.
	This terminology follows from \cite[Section 1.1]{KS99}, which means $\ag T_s$ belongs to an $s\tau$-stable Borel pair $(\ag B_s, \ag T_s)$ in $\ag G$.

\begin{lem}\label{lem: restricted-roots}
	Since $\ag T_s$ is admissible, any root $\alpha \in \RS({\lag g},{\lag t}_s)$ has a non-zero restriction to ${\lag t}_s^{s\tau}$.
\end{lem}
\begin{pf}
	Let us show that the automorphism $s\tau$ of ${\lag t}_s$ has finite order.
	Take $g_s\in \ag G(\C)$ such that $\Ad(g_s)( \ag B_s, \ag T_s) = (\ag B, \ag T)$. Then $\Ad(g_s)^{-1} \circ \tau \circ \Ad(g_s)$ also preserves $(\ag B_s, \ag T_s)$.
	Let $t_s : = g_s^{-1} \tau(g_s) s^{-1} \in \ag G(\C)$. Then $\Ad(t_s) = \Ad(g_s^{-1}\tau g_s) \circ \Ad(s \tau)^{-1}$ preserves $(\ag B_s, \ag T_s)$, and hence $t_s \in \ag T_s(\C)$.
	Consequently, as an automorphism of ${\lag t}_s$, $s\tau = t_s^{-1} \cdot g_s^{-1} \tau g_s$ coincides with $\tau_s: = g_s^{-1} \tau g_s$, which is of finite order $d$.

	Since $\tau_s$ preserves $(\ag B_s, \ag T_s)$, for each $\alpha \in \RS^+({\lag g}, {\lag t}_s)$, $\frac 1 d \sum_{i =1}^d \tau_s^{i}\alpha$ is a positive combination of positive roots.
	In particular, it never vanishes on ${\lag t}_s$.
	Note that it lies in $({\lag t}_s^*)^{s\tau} \cong ({\lag t}_s^{s\tau})^*$ and coincides with $\alpha|_{{\lag t}_s^{s\tau}}$ on ${\lag t}_s^{s\tau}$, so $\alpha$ has a non-zero restriction to ${\lag t}_s^{s\tau}$.
	The conclusion for negative roots follows directly.
\end{pf}

	Denote the (multi-)set of restrictions of $\RS^+({\lag g}, {\lag t}_s)$ to ${\lag t}_s^{s\tau}$ by $\RS^+_\res({\lag g}, {\lag t}_s)$, and an element thereof by $\alpha_\res$.
	Note that $\RS^+_\res({\lag g}, {\lag t}_s)$ equals $\RS({\lag g}, {\lag t}_s^{s\tau})$ in the notation of Subsection \ref{subsec: notation}.

\subsection{Twisted Weyl denominator}\label{subsec: twisted-Weyl-denominator}
	Let $s \tau \in G^+$ be semisimple, and let $(\ag B_s, \ag T_s)$ be an $s\tau$-stable Borel pair with $\ag T_s$ defined over $\R$.
Denote the identity component of the real Lie group $T_s^{s\tau} = \ag T_s^{s\tau}(\R)$ by $Z_s$, whose complex Lie algebra ${\lag z}_s$ is exactly ${\lag t}_s^{s\tau}$.
	For $t\in Z_s$, define
\[ \Delta_G^{s\tau} (t) = \Delta_{G^+}(t s\tau) = 
	\left| \det_{{\lag u}_s} (ts \tau) \right|^{-\frac 1 2} \cdot 
	\left | \det_{{\lag u}_s}(1 - t s\tau) \right|, \]
	where ${\lag u}_s$ is the nilpotent radical of ${\lag b}_s$.
	Our $\Delta_{G}^\tau$ is the (modulus of) twisted Weyl denominator $\Delta_G = \Delta_{G, \tau}$ in \cite[Subsection 2.3]{BC13}.
	In this subsection, we present two formulas for $\Delta_G^{s\tau}$.

\begin{lem}\label{lem: twisted-Weyl}
\begin{subequations}\label{eq: formulas}
	For any $\alpha_\res \in \RS^+_\res({\lag g}, {\lag t}_s)$, $s\tau$ stabilizes the root space ${\lag g}_{\alpha_\res}$.
	Take the decomposition ${\lag t}_s = {\lag z}_s \oplus \lag c_s$ into $s\tau$-invariant subspaces, and let
\[{\lag q}_s =\lag c_s \oplus \bigoplus_{\alpha \in \RS^+_\res ({\lag g}, {\lag t}_s)}
		{\lag g}_{\alpha_\res} \oplus {\lag g}_{-\alpha_\res}.\]

\begin{enumerate}[label=$(\arabic*)$]
\item
	Denote the eigenvalues (with multiplicity) of $s\tau$ on ${\lag g}_{\alpha_\res}$ by $\lambda_{\alpha_{\res}, j}$, $j = 1,  \cdots, \dim_\C {\lag g}_{\alpha_\res}$. Then
\begin{equation}\label{eq: twisted-Weyl-denominator}
	\Delta_G^{s \tau}(t) = \prod_{\alpha_\res \in \RS^+_\res({\lag g}, {\lag t}_s)}
	\prod_{j=1}^{\dim_\C {\lag g}_{\alpha_\res}} 
		\left|\lambda_{\alpha_\res, j} \alpha_\res( t  ) \right|^{-\frac 1 2}
		\left|1 - \lambda_{\alpha_\res, j} \alpha_\res(t ) \right|.
	\end{equation}

\item
	Recall the notation $D_{{\lag q}_s}^{s\tau}(t) = \det_{{\lag q}_s} \left( 1 - t s \tau \right)$. Then
\begin{equation}\label{eq: det-q}
	\left| D_{{\lag q}_s}^{s\tau}(t)\right| = \left| D_{\lag c_s} (s\tau) \right|
		\Delta_G^{s\tau}(t) ^2.
	\end{equation}
	In particular, $\Delta_G^{s\tau}$ does not depend on the choice of $s\tau$-stable Borel $\ag B_s$ that contains $\ag T_s$.
\end{enumerate}
\end{subequations}
\end{lem}
\begin{pf}
	(1) Since
\[{\lag u}_s = \bigoplus_{\alpha_\res \in \RS_\res^+({\lag g}, {\lag t}_s)} {\lag g}_{\alpha_\res},\]
	and $t$ and $s\tau$ commute on each ${\lag g}_{\alpha_\res}$, with eigenvalues $\alpha_\res(t)$ and $\lambda_{\alpha_\res, j}$, equation \eqref{eq: twisted-Weyl-denominator} follows from the definition.

	(2) Let us calculate the eigenvalues of $ts\tau$ on ${\lag g}_{-\alpha_\res}$.
	We will construct a $G^+$-invariant non-degenerate pairing $\la-, - \ra$ on ${\lag g}$ in Subsection \ref{subsec: tau-stable-subgroup}.
	Then by considering the action of $Z_s$ one deduces that ${\lag g}_{\alpha_\res}$ is paired only with ${\lag g}_{-\alpha_\res}$, and by considering the action of $ts\tau$
one deduces that $ts \tau$ acts on ${\lag g}_{-\alpha_\res}$ by  $\lambda_{\alpha_\res, j}^{-1} \alpha_\res(t)^{-1}$, $j = 1, \cdots, \dim_\C {\lag g}_{\alpha_\res}$.
	Consequently,
\[\begin{aligned}
	\det_{{\lag q}_s}(1 - ts\tau) = &
	\det_{\lag c_s}(1 - ts\tau) \prod_{\alpha_\res \in \RS_\res^+({\lag g}, {\lag t}_s)}
	\prod_{j=1}^{\dim_\C {\lag g}_{\alpha_\res}} \left( 1- \lambda_{\alpha_\res, j}\alpha_\res(t) \right)
	\left( 1 - \lambda_{\alpha_\res, j}^{-1} \alpha_\res(t)^{-1} \right) \\
	= &
	\det_{\lag c_s}(1 - s\tau) \prod_{\alpha_\res \in \RS_\res^+({\lag g}, {\lag t}_s)}
	\prod_{j=1}^{\dim_\C {\lag g}_{\alpha_\res}}
	(-\lambda_{\alpha_\res, j} \alpha_\res(t) )^{-1}
	\left( 1- \lambda_{\alpha_\res, j}\alpha_\res(t) \right)^2,
	\end{aligned}\]
	so $|D_{{\lag q}_s}^{}(ts\tau)| = 
	|D_{\lag c_s}(s\tau)| \Delta_{G^+}(ts\tau)^2.$
\end{pf}

\subsection{$\tau$-regular elements}\label{subsec: tau-regular}
	We call an element $s\in G$ $\tau$-regular if $s \tau\in G^+$ is semisimple and $\dim_\R {\lag g}_0^{s \tau}$ attains the minimum of $\dim_\R {\lag g}_0^{g \tau}$ for $g\tau \in G\tau$.
	We denote the subset of $\tau$-regular elements in $G$ by $G'$.

\begin{lem}\label{lem: adm-Cartan}
	If $s \in G$ is $\tau$-regular, then $\lag g^{s\tau} \subseteq \lag g$ is Abelian, and $\ag T_s = Z(\ag G, {\lag g}^{s\tau})$ is the unique $s\tau$-admissible Cartan subgroup.
\end{lem}
\begin{pf}
	From the proof of Lemma \ref{lem: restricted-roots}, we have $\dim_\R \lag g_0^{s\tau} \geqslant \dim_\C \lag t_s^{s\tau} = \dim_\C \lag t^\tau$.
	On the other hand, from \eqref{eq: twisted-Weyl-denominator} we know that there does exist $s\in \ag T^\tau(\R)^\circ$ such that $\dim_\R \lag g_0^{s\tau}$ attains $\dim_\C \lag t^\tau$.
	Thus, $s \in G$ is $\tau$-regular if and only if $s\tau$ is semisimple and $\dim_\C \lag g^{s\tau}$ attains $\dim_\C \lag t^\tau$.
	In this case, $\lag g^{s\tau}$ has the same dimension as $\lag t_s^{s\tau}$, where $\lag t_s$ is an $s\tau$-admissible Cartan subalgebra in $\lag g$, so $\lag g^{s\tau} = \lag t_s^{s\tau}$ is Abelian.

	Again by Lemma \ref{lem: restricted-roots}, each root space $\lag g_\alpha$ of $\ag T_s$ does not commute with $\lag t_s^{s\tau}$, so $Z(\lag g, \lag g^{s\tau}) = Z(\lag g, \lag t_s^{s\tau}) = \lag t_s$.
	Since $\ag G$ is connected, $Z(\ag G, \lag g^{s\tau})$ is also connected (see, for example, \cite[Theorem 6.7.4]{Springer98}).
	Thus, $Z(\ag G, \lag g^{s\tau})$ coincides with the Cartan subgroup $\ag T_s$, and is the unique $s\tau$-admissible Cartan subgroup.
\end{pf}

\begin{rmk}\label{rmk: criteria-of-regular}
	If $s\in G$ is $\tau$-regular, then the unique $s\tau$-admissible Cartan subgroup $\ag T_s = Z(\ag G, {\lag g}^{s\tau})$ is defined over $\R$.
	In this case, we know that $\lag g^{s\tau} = \lag t_s^{s\tau}$, so $D_{\lag q_s}(s\tau)\not = 0$, and then $\Delta_{G^+}(s\tau) \not = 0$ due to \eqref{eq: det-q}.
\end{rmk}

\subsection{$\tau$-stable maximal compact subgroup}\label{subsec: tau-stable-subgroup}
	In this and the next subsection, we demonstrate the existence of $\tau$-stable maximal compact subgroups and real parabolic subgroups.
	This indicates the generality of our Theorem \ref{thm: main}.
\begin{eg}
	In the case $(\GL_{2m+1}, \tau)$ of Theorem \ref{thm: compatibility},  the standard maximal compact subgroup $\O(2m+1)$ and the standard (block upper-triangular) parabolic subgroup with Levi subgroup $\GL_1(\R) \times \GL_{2m-1}(\R) \times \GL_1(\R)$ are both stable under $\tau$.
\end{eg}

	We follow the definitions in \cite{AT18}:
\begin{itemize}
\item
	a real form of $\ag G(\C)$ is equivalently defined as an algebraic conjugate-linear involution;
\item
	the $\R$-structure of $\ag G$ induces a natural real form $\sigma$;
\item
	a Cartan involution of the real reductive group $\ag G$ is an algebraic involution such that  $\theta$ commutes with $\sigma$, and $\ag G(\C)^{\theta \sigma}$ is compact.
\end{itemize}

\begin{fact}\label{fact: Cartan-involution}
	The real reductive group $\ag G$ has a Cartan involution, and it is unique up to conjugation by an inner automorphism from $G^\circ$.
	A Cartan involution of $\ag G$ induces a Cartan involution of ${\lag g}_{\der, 0} = [{\lag g}_0, {\lag g}_0]$.
\end{fact}
\begin{pf}
	The first claim is stated in \cite[Theorem 3.13(1a)]{AT18}.
	We deduce the second claim from the first claim.
  Let $\theta$ be a Cartan involution of $\ag G$ and $B$ be the Killing form of $\lag g$. 
	We follow the definition as given in \cite[Page 355]{Knapp96}, so it suffices to show that the symmetric bilinear form
$\la X, Y \ra_\der:=- B(\theta(X), Y)$ is positive definite on ${\lag g}_{\der, 0}$.
 By definition, $\sigma_{\mathrm c}: = \theta \sigma$ is a compact form of the complex algebraic group $\ag G$.
 Suppose $\lag g_{\der, 0} = \lag k_0 \oplus \lag p_0$ is the decomposition into eigenspaces of $\theta$.
 Then the compact real form $\lag u_0 := \lag g_\der^{\sigma_{\mathrm c}}$ decomposes as $\lag k_0 \oplus \mathrm i\lag p_0$.
 It is well-known that the Killing form is negative definite on the compact real form of a semisimple complex Lie algebra.
 Then one obtains the positive definiteness of $\la -, - \ra_\der$ easily.
\end{pf}

\begin{lem}\label{lem: tau-maximal-cpt}
	There is a Cartan involution $\theta$ of $\ag G$ which commutes with $\tau$.
\end{lem}

\begin{pf}
	This result is also mentioned in \cite[Subsection 3.5]{BC13}, where it is briefly explained that the finite-order automorphism $\tau$ must have a fixed point on the symmetric space of $G$.
	Here is a more detailed demonstration.

	As mentioned in Fact \ref{fact: Cartan-involution}, $\ag G$ has a Cartan involution $\theta_0$, and  there is a surjection $G^\circ \to \{$Cartan involutions of $\ag G\}$, $g \mapsto \Ad(g) \circ \theta_0 \circ \Ad(g)^{-1}$.
	We claim that it descends to a bijection from $G^\circ / \left(Z(G)K_0  \cap G^\circ \right)$ to the set of Cartan involutions of $\ag G$, where $Z(G)$ is the center of $G$, and $K_0=G^{\theta_0}$ is a maximal compact subgroup of $G$.
	This is equivalent to saying the stabilizer of $\theta_0$ in $G^\circ$ is contained in $Z(G)K_0$.
	Suppose $g \in G$ stabilizes $\theta_0$, i.e., $\Ad(g) \circ \theta_0 \circ \Ad(g)^{-1} = \theta_0$. Then $\Ad(g^{-1} \theta_0(g) ) = \Id_G$.
	Thus $g^{-1} \theta_0(g) = z \in Z(G)$.
	Under the Iwasawa decomposition $G= K_0 A_{\min} N_{\min}$ with respect to $\theta_0$, we write $g = k a n$ and deduce that $n^{-1} a^{-1} \theta(a) \theta(n) = z$.
	Since $\theta(a)=a^{-1}$ and $\theta(N_\min)=N_\min^-$, the opposite of $N_\min$, the Bruhat decomposition of $\ag G$ implies that $n=1$.
	Now $a^2 = z^{-1} \in Z(G)$.
	Take $X \in {\lag a}_{\min, 0}$ such that $a = \exp X$; then, for each $\beta \in \RS({\lag g}, {\lag a}_\min)$, we have $\la \beta, X \ra = \frac 1 2 \la \beta, 2X \ra = 0$.
	Hence, $a= \exp X$ acts trivially on each ${\lag g}_\beta$.
	Thus, $a \in Z(G)$, and $g = k a \in Z(G) K_0$.

	Since $\ag G$ is connected, $\ag G(\C)$ is generated by exponentials of its Lie algebra ${\lag g}$.
	Then $Z(G^\circ) \subseteq Z(\ag G)$, so $Z(G) \cap G^\circ = Z(G^\circ)$.
	On the other hand, $K_0$ intersects every component of $G$, so $(Z(G) K_0) \cap G^\circ = Z(G^\circ) K_0^\circ$.
	Now $G^\circ/ \left( Z(G) K_0 \cap G^\circ\right) =  G^\circ/ Z(G^\circ) K_0^\circ$ is the quotient of the semisimple Lie group $G^\circ/ Z(G^\circ)$ by its maximal compact subgroup, which is known as a symmetric space.
	Recall Casselman's fixed-point theorem for semi-hyperbolic spaces in \cite{Casselman12}:
\begin{itemize}
\item
	the symmetric space of a semisimple group is semi-hyperbolic, since its exponential maps are locally expanding (its Proposition 2.1, Lemma 3.1), and
\item
	any compact group of isometries on a complete semi-hyperbolic space has a fixed point (its Corollary 1.4).
\end{itemize}
	It follows that any finite-order automorphism of $G^\circ/ \left( Z(G) K_0 \cap G^\circ\right)$ has a fixed point.

	The conjugation by $\tau$ on Cartan involutions of $\ag G$ induces a finite-order action on $G^\circ/ \left( Z(G) K_0 \cap G^\circ\right)$ via the previous bijection.
	As a consequence, we obtain a $\tau$-fixed Cartan involution $\theta$ from the above results of Casselman.
\end{pf}

	Let $\theta$ be a Cartan involution that commutes with $\tau$. Then $K: = G^\theta$ is a maximal compact subgroup of $G$ that is $\tau$-stable.
	Define a $G^+$-invariant pairing $\la-, -\ra$ on ${\lag g}_0$ as follows.
	Let $\lag c_0$ be the center of ${\lag g}_0$ and ${\lag g}_{\der, 0} = [{\lag g}_0, {\lag g}_0]$; then ${\lag g}_0 = {\lag g}_{\der, 0} \oplus  \lag c_0$.
	Since $\tau$ is of finite order, there exists a $\la \theta \ra \times \la \tau \ra$-invariant inner product $\la-, -\ra_Z$ on $\lag c_0$.
	The inner product $\la-, -\ra_\der = -\Tr_{{\lag g}_\der}(\ad(\theta(-)) \circ \ad(-) )$ as defined in the proof of Fact \ref{fact: Cartan-involution} is also invariant under $\theta$ and $\tau$.
	Then they induce an inner product $\la-, -\ra_+$ on ${\lag g}_0 = {\lag g}_{\der, 0} \oplus \lag c_0$ naturally.
	One can verify that $\la-, -\ra :=- \la -, \theta(-) \ra_+$ is invariant under $G^+$.

	The 4-tuple $(G, K, \theta, \la-, -\ra)$ defines a reductive group in the sense of \cite[Definition 4.29]{KV95} and admits the action of $\tau$.

\begin{rmk}
	Note that \cite[Lemma 4.19]{Mezo16} also discusses the existence of a maximal compact subgroup that is stable under a semisimple automorphism.
	His approach seems different from ours, and considers only the case where $G$ contains an elliptic torus.
\end{rmk}

\subsection{$\tau$-stable real parabolic subgroup}\label{subsec: stable-real-para}
	Now we discuss the $\tau$-stable parabolic subgroups of $G$.
	Recall that we have fixed a $\tau$-stable Borel pair $(\lag b, \lag t)$ in $\lag g$.
	According to \cite[Proposition 4.57]{KV95}, the parabolic subalgebras ${\lag p} \supseteq {\lag b}$ in ${\lag g}$ are in bijection with the subsets $\Pi$ of simple roots $\Pi({\lag g}, {\lag t})$.
	If $\Pi$ is $\tau$-stable, then the corresponding parabolic subalgebra ${\lag p}$ and its Levi decomposition ${\lag m} \oplus {\lag n}$ are also $\tau$-stable.
	Here ${\lag n}$ is the radical of ${\lag p}$, and ${\lag m} \supseteq {\lag t}$ is determined by specifying its root spaces.

	To obtain real parabolic subgroups, we assume ${\lag t}$ is real and stable under $\theta$.
	According to \cite[Theorem 1.1.A]{KS99}, the $\tau$-admissible Cartan subalgebras of ${\lag g}$ are in bijection with Cartan subalgebras of ${\lag g}^\tau$.
	Since $\theta$ commutes with $\tau$, it also induces a Cartan involution on ${\lag g}_0^\tau$.
	By \cite[Proposition 6.59]{Knapp96}, there are sufficiently many $\theta$-stable Cartan subalgebras in ${\lag g}_0^\tau$.
	After conjugating by an element of $\ag G^\tau$, we can assume that ${\lag t}^\tau \subseteq {\lag g}^\tau$ is $\theta$-stable and is the complexification of ${\lag t}_0^\tau \subseteq {\lag g}_0^\tau$.
	Then ${\lag t} := Z({\lag g}, {\lag t}^\tau)$ is the required $\theta$-stable $\tau$-admissible Cartan subalgebra, and $\ag T = Z( \ag G, {\lag t}^\tau)$ is now defined over $\R$.

	If $\RS({\lag p}, {\lag t}) = \RS^+({\lag g}, {\lag t}) \cup (\RS({\lag g}, {\lag t}) \cap \span \Pi)$ is stable under the complex conjugation induced by ${\lag g}= {\lag g}_0 \otimes \C$, then ${\lag p}$ is real as introduced after \cite[Proposition 4.74]{KV95}: ${\lag p} = {\lag p}_0 \otimes_\R \C$ with ${\lag p}_0 = {\lag p} \cap {\lag g}_0$.
	Following the construction in \cite[pages 656--659]{KV95}, we obtain a Levi decomposition $P=MN$ of the real parabolic subgroup of $G$.
	Explicitly, $M = N(G, {\lag p}) \cap N(G, \theta {\lag p})$, and  $N = \exp {\lag n}_0$.
By the above construction, $P$, $M$, $N$ are real points of algebraic groups $\ag P$, $\ag M$, $\ag N$, where $\ag P \supseteq \ag B$ and $\ag M \supseteq \ag T$.
	Moreover, if $\Pi$ is $\tau$-stable, then the above construction is also $\tau$-stable.
	In particular, the subgroups $\ag P$, $\ag M$, $\ag N$ are all $\tau$-stable.

\begin{lem}\label{lem: tau-Borel}
	Suppose $s\tau \in M^+$ is semisimple; then $(\ag P, \ag M)$ contains an $s\tau$-stable Borel pair $(\ag B_s, \ag T_s)$ of $\ag G$ that is conjugate to the previous $(\ag B, \ag T)$ by an inner automorphism from $m_s \in \ag M(\C)$.
	In particular, it holds on $M\tau$ that
\begin{equation}\label{eq: descent}
	\Delta_{G^+}= |\det_{\lag n}|^{-\frac 1 2} | D_{\lag n} | \Delta_{M^+}.
	\end{equation}
\end{lem}
\begin{pf}
	According to \cite[Theorem 7.5]{St68}, there is an $s\tau$-stable Borel pair $(\ag B_{M, s}, \ag T_s)$ in $\ag M$.
	Since $\ag T_s$ is also a Cartan subgroup in $\ag G$, one can check from the root space decomposition that $\ag B_{M, s}$ and $\ag N$ generate a Borel subgroup $\ag B_s$ of $\ag G$, which is $s\tau$-stable and contains $\ag T_s$ by definition.

	Both $\ag T$ and $\ag T_s$ are Cartan subgroups of $\ag M$, so there exists $m_s' \in \ag M$ such that $\Ad(m_s') \ag T_s = \ag T$.
	Then $\Ad(m_s') \ag B_s$ is also a Borel subgroup of $\ag G$ that contains $\ag T$ and $\ag N$.
	One can find another $m_s'' \in N(\ag M, \ag T)$ such that $\Ad(m_s'')$ sends $\Ad(m_s') {\lag b}_s \cap {\lag m}$ to ${\lag b} \cap {\lag m}$, so $m_s: = m_s'' m_s' \in \ag M$ makes $\Ad(m_s)(\ag B_s, \ag T_s) = (\ag B, \ag T)$.

	For $s\tau \in M'\tau$, we can choose the unique $s\tau$-admissible Cartan subalgebra ${\lag t}_s = Z(\lag m, \lag m^{s\tau})$ as in Lemma \ref{lem: adm-Cartan}, and $\RS^+({\lag g}, {\lag t}_s) = \RS^+({\lag m}, {\lag t}_s) \sqcup \RS({\lag n}, {\lag t}_s)$ as above.
	Then one deduces from the definition that \eqref{eq: descent} holds for such $s\tau$.
	Since both sides of \eqref{eq: descent} are continuous, and $M' \tau$ is dense in $M\tau$, it holds on the whole $M\tau$.
\end{pf}

	The Levi subgroup $M$ has a central split torus $A = \exp {\lag a}_0$, where ${\lag a}_0$ consists of central elements in ${\lag m}_0$ on which $\theta$ acts by $-1$.
	Since $M$ is $\tau$-stable, $A$ is  $s\tau$-stable for any $s\in M$.
	An analogous statement to Lemma \ref{lem: restricted-roots} holds for $\RS^+({\lag g}, {\lag a}) : = \RS({\lag n}, {\lag a})$: each $\beta \in \RS^+({\lag g}, {\lag a})$ has non-trivial restriction to ${\lag a}^\tau$.
	In particular, ${\lag a}^\tau$ and $A^\tau$ are non-trivial, unless $P = G$.
	We will exclude this trivial case in what follows.
\begin{cor}\label{cor: big-enough}
	For any $s \in M'$, there is an open subset $A_s^- \subset A$ such that $sA_s^- \subseteq M^-$.
\end{cor}
\begin{pf}
	For $s \in M'$, which is $\tau$-regular in $M$, let us first  find $a\in A^\tau$ such that $sa \in G'$.
By \eqref{eq: descent}, $\Delta_{G^+}(sa\tau) = |\det_{\lag n}(sa\tau)|^{-\frac 1 2} |D_{\lag n}(sa\tau)| \Delta_{M^+}(s\tau)$ for any $a \in A^\tau$.
	We know from Remark \ref{rmk: criteria-of-regular} that $\Delta_{M^+}(s\tau) \not = 0$.
	Then from \eqref{eq: det-q} we know $sa\tau \in G^+$ is regular if $|D_{\lag n}(sa\tau)| \not = 0$.
Since $s\tau$ commutes with $A^\tau$, it preserves the decomposition
\[{\lag n} = \bigoplus_{\beta_\res \in \RS^+({\lag g}, {\lag a}^\tau)}{\lag g}_{\beta_\res}.\]
	Thus, the eigenvalues of $sa\tau$ on ${\lag g}_{\beta_\res}$ are obtained by scaling those of $s\tau$ by $\beta_\res(a)$.
	When all $\beta_\res(a)$, $\beta_\res \in \RS^+({\lag g}, {\lag a}^\tau)$, are sufficiently small, the eigenvalues of $\Ad(sa) \circ \tau$ on ${\lag n}$ have modulus $<1$.
	Consequently, for such $a\in A^\tau$ we know $\Delta_{G^+}(sa\tau) = |\det_{\lag n}(sa\tau)|^{-\frac 1 2} |D_{\lag n}(sa\tau)| \Delta_{M^+}(s\tau)$ does not vanish, and hence $sa \in M$ is $\tau$-regular in $G$.

	Our choice of $a$ also implies $sa\in M^-$.

	Now let $A^- = \{ a \in A \mid \beta(a) <1, \forall \beta \in \RS({\lag n}, {\lag a})\}$ be an open subset of $A$.
	It suffices to show that for each $s \in M^-$ and $a \in A^-$, all eigenvalues of $\Ad(sa) \circ \tau$ on ${\lag n}$ have modulus $<1$.
	One calculates that $(sa\tau)^d = sa \cdot \tau(sa) \cdots \tau^{d-1}(sa) = 
	\left( s \tau(s) \cdots \tau^{d-1}(s) \right) \cdot
	\left( a \tau(a) \cdots \tau^{d-1}(a) \right) = 
	(s \tau)^d \cdot \left( a \tau(a) \cdots \tau^{d-1}(a) \right)$, where $\left( a \tau(a) \cdots \tau^{d-1}(a) \right) \in A^\tau \cap A^-$.
	It is easy to verify from the definition that, as an automorphism of ${\lag n}$, $(s \tau)^d \cdot \left( a \tau(a) \cdots \tau^{d-1}(a) \right)$ has all eigenvalues of modulus $<1$.
	Then the conclusion also applies to $sa\tau$.
\end{pf}

	In the rest of this article, for every $s\in M$ that is $\tau$-regular in $G$, we fix an $s\tau$-stable Borel pair $(\ag B_s, \ag T_s)$ and $m_s \in \ag M(\C)$ such that $\Ad(m_s) ( \ag B_s, \ag T_s) = (\ag B, \ag T)$ as in Lemma \ref{lem: tau-Borel}.
In particular, we have fixed $t_s \in \ag T_s(\C)$ such that $\tau_s =  t_s s\tau$ preserves the splitting $\spl_s = \Ad(m_s)^{-1} \spl$ of $\ag G$.

\subsection{$\tau$-stable infinitesimal character}\label{subsec: stable-inf-char}
	Consider the infinitesimal character of a $\tau$-stable irreducible representation.
	It is $\tau$-stable when regarded as a multiplicative functional on $Z({\lag g})$.
	By the Harish-Chandra isomorphism $Z({\lag g}) \cong \C[{\lag t}]^{W({\lag g}, {\lag t})}$, this infinitesimal character corresponds to a $W({\lag g}, {\lag t})$-orbit in ${\lag t}^*$.
	Recall that we have assumed ${\lag t}$ to be $\tau$-admissible.
	In this subsection, we describe the $\tau$-stable representatives in the $\tau$-stable orbit.

\begin{lem}\label{lem: stable-infinitesimal-character}
	Since ${\lag t}$ is a $\tau$-admissible Cartan subalgebra of ${\lag g}$, the natural inclusion $({\lag t}^*)^\tau \to {\lag t}^*$ induces a bijection between $W({\lag g}, {\lag t})^\tau$-orbits in $({\lag t}^*)^\tau$ and $\tau$-stable $W({\lag g}, {\lag t})$-orbits in ${\lag t}^*$.
	In other words, a  $\tau$-stable  $W({\lag g}, {\lag t})$-orbit in ${\lag t}^*$ must contain a $\tau$-stable representative $\lambda$, and all its $\tau$-stable representatives belong to $W({\lag g}, {\lag t})^\tau \lambda$.
\end{lem}
\begin{pf}
	This lemma is equivalent to the last sentence of \cite[Subsection 3.2]{BC13}.
	Here we give an explicit proof.
	Suppose the $W({\lag g}, {\lag t})$-orbit of $\lambda \in {\lag t}^*$ is $\tau$-stable.
	This is equivalent to saying that $\tau(\lambda)=\lambda\circ\tau^{-1}$ has the form $w_\tau\lambda$ for some $w_\tau\in W({\lag g}, {\lag t})$. Then $\tau(w\lambda)=\tau(w)\tau(\lambda)=\tau(w)w_\tau\lambda\in W({\lag g}, {\lag t})\lambda$ for any $w\in W({\lag g}, {\lag t})$.
	We have to find some $w \in W({\lag g}, {\lag t})$ such that $\tau(w \lambda) = w\lambda$.
	Let $\RS^+, \RS^- \subset \RS({\lag g}, {\lag t})$ be the $\tau$-stable positive and negative root systems determined by ${\lag b}$, and consider the following Weyl chamber
\[\{ \xi \in {\lag t}^* \mid \forall \alpha \in \RS^+, \textrm{ either } 
	\Re \la \xi, \alpha \ra >0, \textrm{ or }
	\Re \la \xi, \alpha \ra = 0 \textrm{ and } \Im \la \xi, \alpha \ra \leqslant 0\}.\]
	It is easy to check the following properties:
\begin{enumerate}[label=$(\arabic*)$]
\item
	This chamber intersects every $W({\lag g}, {\lag t})$-orbit in ${\lag t}^*$;
\item
	this chamber meets each $W({\lag g}, {\lag t})$-orbit at most once;
\item
	this chamber is preserved by  $\tau$.
\end{enumerate}
	By (1), we may assume the representative $\lambda$ lies in this chamber.
	Then by (3), $\tau(\lambda)$ also lies in this chamber.
	On the other hand, by (2), $w_\tau\lambda$ can lie in this chamber only if $w_\tau\lambda=\lambda$.
	Then this dominant representative $\lambda$ must be $\tau$-stable.

	Let $W({\lag g}, {\lag t})_\lambda$ be the stabilizer of $\lambda$ in $W({\lag g}, {\lag t})$.
	Suppose there is another $\tau$-stable representative $w\lambda$, $w \in W({\lag g}, {\lag t})$; then we have to find some $u \in W({\lag g}, {\lag t})_\lambda$ such that $w^{-1} \tau(w) = u^{-1} \tau(u)$; it follows that $w\lambda = w u^{-1} \lambda$ with $wu^{-1} \in W({\lag g}, {\lag t})^\tau$.
	We find such a $u$ by induction on the length of $w$.
	Let $\RS_0 \subseteq \RS({\lag g}, {\lag t})$ consist of roots orthogonal to $\lambda$ and $\RS_0^+ = \RS_0 \cap \RS^+$.
	Denote the simple roots in $\RS^+$ by $\alpha_1, \cdots, \alpha_n$, and the corresponding simple reflections in $W({\lag g}, {\lag t})$ by $s_1, \cdots, s_n$.
	We may assume that $\la \lambda, \alpha_i \ra =0$ for only $i =1, \cdots, r$.
	Then $\RS_0$ equals the intersection of $\RS({\lag g}, {\lag t})$ and the linear subspace spanned by $\alpha_1, \cdots, \alpha_r$.
	If $\RS_0^+ \cap (w^{-1} \RS^-)$ is non-empty, then there exists $i \in \{1, \cdots, r\}$ such that $w\alpha_i <0$.
	Since $(s_i\RS_0^+) \cap \RS^+ = (s_i\RS_0^+) \cap \RS_0^+ = \RS_0^+ \backslash\{\alpha_i\}$, there is an injection from
\[\RS_0^+ \cap ((ws_i)^{-1} \RS^-) =\{ \alpha \in \RS_0^+ \mid ws_i \alpha <0\}\]
	to
\[\RS_0^+ \cap (w^{-1} \RS^-) = \{\alpha \in \RS_0^+ \mid w \alpha <0\}\]
	via $\alpha \mapsto s_i \alpha$, with the image excluding $\alpha_i$.
	If the result holds for $ws_i$ (note that $ws_i \lambda = w\lambda$ since $\alpha_i$ is orthogonal to $\lambda$), then there exists $u \in W({\lag g}, {\lag t})_\lambda$ such that $s_i w^{-1} \tau(w s_i) = u^{-1} \tau(u)$, which implies that $w^{-1} \tau(w) = s_iu^{-1} \tau(u s_i)$ for $us_i \in W({\lag g}, {\lag t})_\lambda$, and this is exactly the desired result for $w$.
	Consequently, after such a formal induction, we may assume $\RS_0^+ \cap (w^{-1} \RS^-) =\varnothing$.

	Given $\tau(w\lambda) = w\lambda$ and $w\RS_0^+ \subseteq \RS^+$, we claim $w^{-1}\tau(w) = 1$.
	Since $\lambda$ is $\tau$-stable, $\tau(w\lambda) = w \lambda$ implies $w^{-1} \tau(w) \lambda = \lambda$, so $w^{-1} \tau(w) \in W({\lag g}, {\lag t})_\lambda$.
	According to Chevalley's lemma, $W({\lag g}, {\lag t})_\lambda$ is generated by root reflections $s_\alpha$, $\alpha \in \RS_0$.
	In other words, it is a Coxeter group generated by simple reflections $s_1, \cdots, s_r$.
	If $w^{-1}\tau(w) \in W({\lag g}, {\lag t})_\lambda$ is non-trivial, then it must send some simple root $\alpha_i \in \RS_0^+$ to $\RS_0^-$.
	Let us denote $w^{-1}\tau(w)(\alpha_i)$ by $-\beta$, with $\beta \in \RS_0^+$; then $w(\beta) = -\tau(w) \alpha_i$. Since $\tau$ preserves $\lambda$ and ${\lag b}$, it also preserves $\RS_0^+ = \RS_0 \cap \RS^+$.
	Consequently, together with $w\RS_0^+ \subseteq \RS^+$, we have
\[\tau(w) \alpha_i \in \tau(w) \RS_0^+ = \tau(w \RS_0^+) \subseteq \tau(\RS^+) = \RS^+.\]
	Then $w(\beta) = - \tau(w) \alpha_i \in \RS^-$ with $\beta \in \RS_0^+$, which contradicts the assumption $w\RS_0^+ \subseteq \RS^+$.
	This completes the proof.
\end{pf}
\begin{rmk}
	This lemma is applicable to $({\lag t}_s, \tau_s)$, and hence to $({\lag t}_s, s\tau)$.
\end{rmk}

	Fix $\lambda \in ({\lag t}^*)^\tau$ that represents the infinitesimal character of a $\tau$-stable irreducible representation.
	For any $s\in G'$, $\lambda_s: = \lambda \circ \Ad(g_s) \in ({\lag t}_s^*)^{s\tau}$ also represents this infinitesimal character, where $g_s \in \ag G(\C)$ conjugates the admissible Cartan subgroup $\ag T_s$ to $\ag T$ as in the proof of Lemma \ref{lem: restricted-roots}.
	If $s \in M \cap G'$, then we can further require $g_s= m_s\in \ag M(\C)$ as fixed at the end of Subsection \ref{subsec: stable-real-para}.
	The subscript $s$ in $\lambda_s$ will be omitted when the context is clear.

\subsection{Very anti-dominance}\label{subsec: very-anti-dominance}
	The notion of very anti-dominance is crucial to the reduction step in our proof of the twisted Osborne conjecture, as will be presented in the next section.
	A weight $\lambda \in {\lag t}^*$ is called ``very anti-dominant'' with respect to $\RS^+({\lag g}, {\lag t})$ if there exists a sufficiently large $C>0$ such that
\[\Re \la \lambda, \alpha\ra < -C, \quad
	\forall \alpha \in \RS^+({\lag g}, {\lag t}).\]
	As for which $C$ is large enough, it depends on the context.
	In this subsection we discuss several properties of very anti-dominant weights.
\begin{fact}\label{fact: centralizer-of-A}
	Recall that $A \subset M$ is its central split torus.
	Let $\pr_{\lag a} \in \End_\C ({\lag t}^*)$ be the orthogonal projection onto ${\lag a}^*$.
\begin{enumerate}[fullwidth, label=$(\arabic*)$]
\item
	$Z(G, A) = M$.
	In particular, if $\alpha \in \RS({\lag g}, {\lag t})$ vanishes on ${\lag a}$, then $\alpha \in \RS({\lag m}, {\lag t})$.
\item
	If $\alpha \in \RS^+({\lag g}, {\lag t})$ does not vanish on ${\lag a}$, then $\pr_{\lag a} \alpha = \alpha +$ a non-negative linear combination of $\Pi({\lag m}, {\lag t})$.
\end{enumerate}
\end{fact}
\begin{pf}
	(1) In \cite[Proposition 4.76(b)]{KV95}, it is shown that ${\lag m} = Z({\lag g}, {\lag a})$, so $\ag M = Z(\ag G, {\lag a})$, $M = Z(G, A)$.
	If $\alpha \in \RS({\lag g}, {\lag t})$ vanishes on ${\lag a}$, then ${\lag g}_{\alpha}$ commutes with ${\lag a}$, so ${\lag g}_\alpha \subseteq {\lag m}$, and $\alpha \in \RS({\lag m}, {\lag t})$.

	(2) Take the orthogonal decomposition ${} {\lag t}= \lag s \oplus Z({} {\lag m}) = \lag s \oplus (Z({} {\lag m}) \cap \lag k) \oplus {} {\lag a}$.
	Let $\pr_Z$ and $\pr_{\lag s} \in \End_\C({} {\lag t}^*)$ be the orthogonal projections onto $Z({} {\lag m})^*$ and $\lag s^*$, respectively.
	Then $\pr_Z \alpha = \alpha - \pr_{\lag s} \alpha$, and
\[\pr_{\lag s}\alpha = \sum_{\alpha' \in \Pi({} {\lag m}, {} {\lag t})} \la \alpha, \omega' \ra \alpha',\]
	where $\{\omega'\} \subset \lag s^* \subset {} {\lag t}^*$ is dual to the basis $\Pi({} {\lag m}, {} {\lag t}) \subset \lag s^*$.
	Due to the properties of simple roots,
\begin{itemize}
\item
	each $\omega'$ is a non-negative linear combination of $\Pi({\lag m}, {} {\lag t})$, and
\item
	for distinct $\alpha, \alpha' \in \Pi({\lag g}, {}{\lag t})$, $\la \alpha, \alpha'\ra \leqslant 0$.
\end{itemize}
	Then $\pr_Z \alpha$ is a non-negative combination of $\Pi({\lag g}, {} {\lag t})$ for $\alpha \in \Pi({\lag g}, {\lag t})$, and the same is true for $\alpha \in \RS^+({\lag g}, {} {\lag t})$.
	The conclusion follows from $\pr_{{} {\lag a}} \alpha = \frac{1}{2}(\pr_Z\alpha + \pr_Z \bar\alpha)$, since $\bar\alpha|_{{} {\lag a}} = \alpha|_{{} {\lag a}} \in \RS^+({\lag g}, {} {\lag a})$ implies $\bar\alpha \in \RS^+({\lag g}, {} {\lag t})$ due to the compatibility between $({} {\lag b}, {} {\lag t})$ and ${} {\lag p} = {} {\lag m} \oplus {} {\lag n}$.
\end{pf}

\begin{lem}\label{lem: very-anti-dominant}\label{lem: mu-component}
	Suppose $\lambda \in {\lag t}^*$ is very anti-dominant. Then for each $w \in W(\ag G, \ag T) = W(\lag g, \lag t)$,
\begin{enumerate}[fullwidth, label=$(\arabic*)$]
\item
	if $w\lambda - \lambda \in {\lag t}^*$ is written as a linear combination of $\Pi({\lag g}, {\lag t})$,  then each non-zero coefficient has positive real part, and
\item
	if $w\lambda |_{{\lag a}} = \lambda|_{{\lag a}}$, then $w \in W(\ag M, \ag T)$.
\end{enumerate}
\end{lem}
\begin{pf}
	(1) Let $w= s_l \cdots s_1$ be a reduced expression, where each $s_i$ is the simple reflection corresponding to $\alpha_i \in \Pi({\lag g}, {\lag t})$.
	We prove by induction that
\[w\lambda - \lambda = \sum_{i=1}^l c_i \alpha_i\]
	with $\Re c_i >0$.
	If $l=0$, there is nothing to do.
	If $l>0$, then since
\[w\lambda - \lambda = 
	w\lambda - s_l w\lambda + s_l w \lambda -\lambda = 
	2 \frac{\la w\lambda, \alpha_l \ra}{|\alpha_l|^2} \alpha_l + s_l w\lambda - \lambda,\]
	it suffices to show $\Re \la w\lambda, \alpha_l\ra >0$.
	It is well-known that $w^{-1}\alpha_l <0$ iff $w^{-1} s_l$ is shorter than $w^{-1}$.
	Since we have chosen $w = s_l\cdots s_1$ to be reduced, it follows that $w^{-1}\alpha_l < 0$ and $\la w\lambda, \alpha_l\ra = \la\lambda, w^{-1}\alpha_l\ra$ has positive real part due to the anti-dominance of $\lambda$.
	Now the conclusion follows.

	(2) By (1), we may write $w= s_l \cdots s_1$ as a reduced expression. Then $w\lambda - \lambda = \sum_{i=1}^l c_i \alpha_i$, where $\Re c_i >0$ and $\alpha_i \in \Pi({\lag g}, {\lag t})$ is the simple root corresponding to $s_i$.
	It follows from Fact \ref{fact: centralizer-of-A}(2) that there exist $d_{i, \alpha} \geqslant 0$ such that
\[\pr_{\lag a} (w\lambda - \lambda) = \sum_{\alpha_i|_{\lag a} \not =0} c_i \left( \alpha_i 
		+ \sum_{\alpha' \in \Pi({\lag m}, {\lag t})} d_{i, \alpha'} \alpha' \right)
	=\sum_{\alpha_i|_{\lag a} \not =0} c_i \alpha_i   + 
	\sum_{\alpha' \in \Pi({\lag m}, {\lag t})} \left( \sum_{\alpha_i|_{\lag a} \not =0} 
		c_i d_{i, \alpha'}\right) \alpha'.\]
	Suppose $w\lambda - \lambda$ vanishes on ${\lag a}$. Then $\pr_{\lag a}( w \lambda - \lambda) = 0 \in {\lag t}^*$.
	Since $\Pi({\lag g}, {\lag t}) \subset {\lag t}^*$ is linearly independent, we know from the above equation that there is no $i$ such that $\alpha_i|_{{\lag a}} \not =0$.
	Thus, each $\alpha_i$ lies in $\RS({\lag m}, {\lag t})$ according to Fact \ref{fact: centralizer-of-A}(1).
	It follows that $w \in W({\lag m}, {\lag t})$.
\end{pf}

	Let $\Lambda = X^*(\ag T)^\tau \subset (\lag t^*)^\tau$ be the lattice of $\tau$-stable integral weights.
  Then any $\lambda \in ({\lag t}^*)^\tau$ can be $\Lambda$-translated into a very anti-dominant weight.
	Indeed, if $\lambda$ is $\tau$-stable, the above condition is also equivalent to
\begin{equation}\label{eq: very-anti}
	\Re \la \lambda, \alpha_\res\ra < -C, \quad
	\forall \alpha_\res \in \RS_\res^+({\lag g}, {\lag t}).
	\end{equation}

	Note that if $\lambda$ is very  anti-dominant with respect to $\RS^+({\lag g}, {\lag t})$, then  so is $\lambda_s = \lambda \circ \Ad(m_s)$ with respect to $\RS^+({\lag g}, {\lag t}_s)$, since the latter positive root system is determined by $({\lag b}_s, {\lag t}_s) = \Ad(m_s) ({\lag b}, {\lag t})$.
	Thus, omitting the subscript $s$ from $\lambda_s$ does not affect the validity of the very anti-dominance condition.
 
\newcommand{\IC}[1]{\mathcal D(#1)}
\newcommand{\CW}{\mathcal{CW}}
\newcommand{\GG}[2][]{\mathrm K_{#1}(#2)}
\newcommand{\FG}[1]{\mathrm F(#1)}
\newcommand{\FC}[1]{\mathcal R(#1)}
\newcommand{\wt}[1]{\mathrm{wt}(#1)}
\newcommand{\IF}[2][]{\mathcal C_{#1}(#2 \tau)}
\newcommand{\HC}{\mathcal{HC}}

\section{Reduction}\label{sec: reduction}
	Our strategy to prove Theorem \ref{thm: main} involves modifying the argument of \cite{BC13} in the minimal parabolic case.
	The main purpose of this section is to reformulate the twisted Osborne conjecture and reduce it via coherent continuation.

	In Subsection \ref{subsec: Grothendieck} we set up notation for the Grothendieck group and the twisted character space.
	In Subsection \ref{subsec: character-of-fd-rep} we introduce a nice class of finite-dimensional representations of $G^+$.
	In Subsection \ref{subsec: local-expression} we establish the local expression \eqref{eq: local-expression-Theta} for the character at regular elements.
	After these preparations, we define the notion of coherent family in Subsection \ref{subsec: coherent-family}, and reduce the twisted Osborne conjecture to the case where the infinitesimal character is represented by a very anti-dominant weight in Subsection \ref{subsec: coherent-continuation}.
	In Subsection \ref{subsec: mu-component} we further reduce the equation to its isotypic component \eqref{eq: reduction-components}.
	For this special case, we outline a proof in the final Subsection \ref{subsec: outline}.

\subsection{Reformulation of the twisted Osborne conjecture}\label{subsec: Grothendieck}
	In this subsection we reformulate the twisted Osborne conjecture in terms of Grothendieck groups.

	Let $\CW(G^+)$ denote the category of Casselman-Wallach representations of $G^+ = G \rtimes \la \tau \ra$.
	Note that the objects in $\CW(G^+)$ are automatically of finite length.
	Denote the Grothendieck group of $\CW(G^+)$ by $\GG{G^+}$.
	For $\lambda \in \lag t^*$, define $\CW_\lambda(G^+)$ to be the full subcategory of $\CW(G^+)$ consisting of representations $\pi^+$ such that each irreducible $G$-subquotient of $\pi^+$ has an infinitesimal character represented by $\tau^i(\lambda)$ for some $i$.
	Denote the Grothendieck group of $\CW_\lambda(G^+)$ by $\GG[\lambda]{G^+}$.
	We have the following natural decomposition
\begin{equation}\label{eq: decompose-of-Grothendieck}
\GG{G^+} = \bigoplus_\lambda \GG [\lambda]{G^+},
	\end{equation}
	where $\lambda$ ranges over the $W(\lag g, \lag t) \rtimes \la \tau \ra$-orbits in $\lag t^*$.

	Suppose $\pi^+$ is an irreducible Casselman-Wallach representation of $G^+$.
	It has a well-defined character $\Theta(\pi^+)$, which is a locally integrable analytic function defined on regular points of $G^+$.
	Recall from Subsection \ref{subsec: tau-regular} that $G'$ consists of the $\tau$-regular elements in $G$.
	Let $\IC{G'\tau}$ be the $\C$-linear space spanned by $\Theta(\pi^+)|_{G'\tau}$, where $\pi^+$ runs over the irreducible Casselman-Wallach representations of $G^+$; then we obtain a natural homomorphism $\GG{G^+} \to \IC{G'\tau}$.
	Note that if $\pi^+|_G$ is reducible, then it is necessarily a direct sum of irreducible $G$-representations, and the action of $\tau \in G^+$ permutes the direct summands.
	In this case, $\Theta(\pi^+)|_{G'\tau}$ vanishes.
	Thus, the function space $\IC{G'\tau}$ can be naturally identified with the space $\TC G$ of twisted characters.

	For a Casselman-Wallach representation $\pi^+$ of $G^+$ with underlying $(\lag g, K^+)$-module denoted by $V$, the $\lag n$-homology groups $H_q(\lag n, V)$ are Harish-Chandra modules of $M^+$.
	We denote the Casselman-Wallach globalization of each $H_q(\lag n, V)$ by $H_q(\lag n, \pi^+)$, as in the Introduction, and define the function
\[\Theta_{\lag n}(\pi^+) = \frac 1 {|D_{\lag n}|} \sum_{q\geqslant 0} (-1)^q \Theta(H_q(\lag n, \pi^+) )\]
	following the notation of \cite{HS83} and \cite{BC13}.
It is easy to see that $\pi^+ \mapsto \Theta_{\lag n}(\pi^+)$ induces a homomorphism from $\GG{G^+}$ to some function space on $(M\cap G')\tau$, and we denote this homomorphism also by $\Theta_{\lag n}$.

	For convenience, we introduce a function space $\IF[\lambda]{\Omega}$ as follows, where $\lambda \in (\lag t^*)^\tau$ and $\Omega \subseteq G'$ has the property that for each $s \in \Omega$, there is a neighborhood $U_s \subset Z_s$ of the identity such that $U_s s \subset \Omega$.
	The space $\IF[\lambda]{\Omega}$ consists of the functions $\Theta: \Omega\tau \to \C$ such that for each $s \in \Omega$, $\Theta$ admits the local expression
\[\left( \Delta_{G^+}\Theta \right)(s \tau \exp X) = \sum_{w \in W(\lag g, \lag t)^\tau} c(s\tau, w) e^{\la w\lambda, X \ra}, \]
	when $X$ lies in a neighborhood of $0 \in \lag z_{s, 0}$; note that we should use $(w\lambda)_s = (w\lambda) \circ \Ad(g_s)$, where $\Ad(g_s)$ sends $(\lag b_s, \lag t_s)$ to $(\lag b, \lag t)$, and we omit the subscript $s$ for convenience.
	We will demonstrate in Subsection \ref{subsec: local-expression} that each twisted character $\Theta(\pi^+)|_{G'\tau}$ whose infinitesimal character is represented by $\lambda$ has a local expression \eqref{eq: local-expression-Theta} as given above.
	It follows that the image of $\GG[\lambda]{G^+}$ under $\GG{G^+} \to \IC{G'\tau}$, $\pi^+ \mapsto \Theta(\pi^+)|_{G'\tau}$, is a subspace of $\IF[\lambda]{G'}$.
	Moreover, the assignment $\pi^+ \mapsto \Theta_{\lag n}(\pi^+)|_{M^-\tau}$ also induces a homomorphism  $\Theta_{\lag n}(-)|_{M^-\tau}: \GG[\lambda]{G^+} \to \IF[\lambda]{M^-}$ according to \eqref{eq: local-expression-n}.
	Thus, we reformulate Theorem \ref{thm: main} as follows:
\begin{thm}\label{thm: main-new}
	For any $\lambda \in (\lag t^*)^\tau$, the homomorphism $\GG[\lambda]{G^+} \to \IF[\lambda]{M^-}$, $\pi^+ \mapsto \Theta_{\lag n}(\pi^+)|_{M^-\tau}$, coincides with the one given by $\pi^+ \mapsto \Theta(\pi^+)|_{M^- \tau}$.
\end{thm}
\begin{pf}[ of the equivalence between Theorem \ref{thm: main} and Theorem \ref{thm: main-new}]
	Suppose $\pi^+$ is an irreducible representation of $G^+$.
	If the restriction of $\pi^+$ to $G$ is also irreducible, then the assertion for $\pi^+$ is exactly the twisted Osborne conjecture in Theorem \ref{thm: main};
  note that $D^\tau_{\lag n}(m) >0$ for any $m \in M^-$, so $|D_{\lag n}(m \tau)| = D^\tau_{\lag n}(m)$.
	If the restriction of $\pi^+$ to $G$ is reducible, then the restriction is a direct sum of irreducible $G$-subrepresentations, and $\tau$ permutes the direct summands.
	In this case, both $\Theta(\pi^+)|_{G'\tau}$ and $\Theta_{\lag n}(\pi^+)|_{M^-\tau}$ vanish.
	Thus, Theorem \ref{thm: main-new} is equivalent to Theorem \ref{thm: main} together with the trivial assertion in the case where $\pi^+|_G$ is reducible.
\end{pf}

\subsection{Finite-dimensional representations}\label{subsec: character-of-fd-rep}
	In this subsection we consider $\tau$-stable finite-dimensional representations of the complex reductive group $\ag G$.

	Let $\FG{G^+}$ be the Grothendieck group of the category of finite-dimensional algebraic representations of $\ag G \rtimes \la \tau \ra$, naturally embedded as a subgroup of $\GG{G^+}$.
	The tensor product induces a multiplication on $\FG{G^+}$ and an action of $\FG{G^+}$ on $\GG{G^+}$.
	Denote the image of $\FG{G^+}$ under $\GG{G^+} \to \IC{G'\tau}$ by $\FC G$.
	Then the multiplication on $\FG{G^+}$ is the multiplication of functions on $\FC G$, and $\IC{G'\tau}$ is also a module over $\FC G$.
	One can verify easily from the definition that for each $\sigma^+ \in \FG{G^+}$ and $\pi^+ \in \GG{G^+}$, we have
\[\Theta(\sigma^+ \otimes \pi^+) = \Theta(\sigma^+) \Theta(\pi^+).\]
	It follows that the homomorphism $\GG{G^+} \to \IC{G'\tau}$ is linear over $\FG{G^+}$.
	The same holds for $\pi^+ \mapsto \Theta_{\lag n}(\pi^+)$:
\begin{lem}\label{lem: F-linear-of-Theta-n}
	For each $\sigma^+ \in \FG{G^+}$ and $\pi^+ \in \GG{G^+}$, we have
\[\Theta_{\lag n}(\sigma^+ \otimes \pi^+) = \Theta(\sigma^+) \Theta_{\lag n}(\pi^+).\]
\end{lem}
\begin{pf}
	This is essentially \cite[Lemma 3.59]{HS83}.
	For a proof, it suffices to assume that both $\sigma^+$ and $\pi^+$ are irreducible representations of $G^+$.
	Denote the representation spaces of $\sigma^+$ and $\pi^+$ by $F$ and $V$, respectively, and introduce a filtration of $P^+$-representations $0 = F_0 \subset F_1 \subset \cdots \subset F_n = F$ by requiring
\[F_i/F_{i-1} = \{ v \in F/ F_{i-1} \mid \lag n \cdot v =0 \} .\]
	It follows that $H_q( \lag n, (F_i/ F_{i-1}) \otimes V) = (F_i/ F_{i-1}) \otimes H_q(\lag n, V)$, so we have the following equality in $\GG{M^+}$:
\[\begin{aligned}
	\sum_{q\geqslant 0} (-1)^q \left[H_q(\lag n, F \otimes V ) \right]=	&
	\sum_{q\geqslant 0} (-1)^q \sum_{i =1}^n [H_q(\lag n,  (F_i/ F_{i-1}) \otimes V) ]	\\
	=	&
	\sum_{q\geqslant 0} (-1)^q \sum_{i=1}^n [(F_i/ F_{i-1})] \otimes [H_q(\lag n, V)]=		
	[F] \otimes \sum_{q\geqslant 0} (-1)^q [H_q(\lag n, V)].
	\end{aligned}\]
	Then the conclusion follows.
\end{pf}

\begin{subequations}
	Suppose $\varphi \in \FC G$ is the restriction of $\Theta(\sigma^+)$ to $G'\tau$, denote the representation space of $\sigma^+$ by $F$, and
	decompose $F$ into weight spaces with respect to $\ag T$ as
\[F = \bigoplus_{\xi \in X^*(\ag T)} F_\xi.\]
	Then we have the following expression for $\varphi(t \tau \exp X)$ with $t \in T$, $X\in \lag t_0^\tau$:
\begin{equation}\label{eq: local-expression-of-fd-T}
	\varphi(t \tau \exp X) = \sum_{\xi \in X^*(\ag T)^\tau} m_\xi \xi(t) e^{\la \xi, X \ra},
	\end{equation}
	where $m_\xi$ is the trace of $\tau$ on the weight space $F_\xi$.
	Only $\tau$-stable weights contribute to $\varphi = \Theta(\sigma^+)|_{G'\tau}$, since $F_\xi$ is permuted to $F_{\tau(\xi)}$ under $\tau$.

	For convenience, set $\wt F = \{\xi \in X^*(\ag T) \mid F_\xi \not=0\}$.
	Introduce a convenient order on $\lag t^* \supset X^*(\ag T)$ such that $\lambda \leqslant \lambda'$ if
\[\lambda' - \lambda= \sum_{\alpha \in \RS^+(\lag g, \lag t)} n_\alpha \alpha\]
	with $n_\alpha$ non-negative integers.
	The following properties of $\wt F$ are well-known; see for example \cite[Theorem 5.5]{Knapp96}:
\begin{itemize}
\item
	there exists a unique highest weight $\xi \in \wt F$ such that for each $\eta \in \wt F$, $\eta \leqslant \xi$;
\item
	$\dim_\C F_{w\eta} = \dim_\C F_\eta$ for any $\eta \in X^*(\ag T)$ and $w \in W(\lag g, \lag t)$, and in particular, $\wt F$ is stable under $W(\lag g, \lag t)$;
\item
	$\dim_\C F_\eta = 1$ for the highest weight $\eta \in \wt F$.
\end{itemize}

\begin{lem}\label{lem: fd-alg-rep}
	Let $\Lambda = X^*(\ag T)^\tau \subset (\lag t^*)^\tau$ be the lattice of $\tau$-stable integral weights.
\begin{enumerate}[label=$(\arabic*)$]
\item
	The group $\FC G$ is a free $\Z[\zeta_d]$-module, with a basis indexed by the dominant weights in $\Lambda$.

\item
	For each $\xi \in \Lambda$, there exists $\psi_\xi \in \FC G$ such that for $X \in \lag t_0^\tau$
\[\psi_\xi(\tau \exp X) = \sum_{w \in W(\lag g, \lag t)^\tau} e^{\la w \xi, X \ra}.\]

\item
	For each $\varphi \in \FC G$, let $\varphi(\tau \exp X)$, $X \in \lag t_0^\tau$, be expressed as in \eqref{eq: local-expression-of-fd-T}; then for each $s \in G'$ and $X \in \lag t_{s, 0}^{s\tau}$
\begin{equation}\label{eq: local-expression-of-fd}
	\varphi(s \tau \exp X) = \sum_{\xi \in \Lambda} m_\xi \xi_s(t_s)^{-1} e^{\la \xi_s, X \ra},
	\end{equation}
	where $\xi_s = \xi \circ \Ad(g_s)$ and $t_s = g_s^{-1}\tau(g_s)s^{-1}$ as defined at the end of Subsection \ref{subsec: stable-inf-char}.
\end{enumerate}
\end{lem}
\begin{pf}
	(1) We have observed in Subsection \ref{subsec: Grothendieck} that an irreducible representation $\sigma^+ \in \FG{G^+}$ has $\Theta(\sigma^+)|_{G'\tau}\not=0$ only when it extends a $\tau$-stable irreducible representation of $\ag G$.
	The action of $d$th roots of unity $\zeta_d$ amounts to a different choice of extension to $G^+$.
	Thus, $\FC G$ is a module over $\Z[\zeta_d]$, and is generated by the set $\mathfrak F$ of finite-dimensional $\tau$-stable irreducible algebraic representations of $\ag G$.

	It is well-known that the isomorphism classes of finite-dimensional irreducible algebraic representations $F$ of $\ag G$ are in bijection with the dominant integral weights $\xi$ in $X^*(\ag T)$, where $\xi$ is the highest weight in the corresponding representation $F$; see for example \cite[Corollary 2.7 in Part II]{Jantzen87}.
	Since the Borel pair $(\ag B, \ag T)$ is stable under $\tau$, the subset of dominant weights in $X^*(\ag T)$ is also stable under $\tau$.
	Thus, the irreducible representation $F$ is $\tau$-stable if and only if its highest weight $\xi$ is $\tau$-stable.
	Then the generating set $\mathfrak F$ of $\FC G$ over $\Z[\zeta_d]$ is also indexed by dominant weights in $\Lambda$.

	One can deduce the linear independence of $\mathfrak F$ over $\Z[\zeta_d]$ from the expression \eqref{eq: local-expression-of-fd-T} together with Dedekind's theorem on the linear independence of characters.

	(2) Without loss of generality, we can assume $\xi\in \Lambda$ is dominant, and $\ag G$ is semisimple.
	Let $F$ be the $\tau$-stable irreducible representation of $\ag G$ with highest weight $\xi$, and  extend it to an irreducible representation $\sigma^+$ of $\ag G \rtimes \la \tau \ra$.
According to Lemma \ref{lem: stable-infinitesimal-character}, for each $\eta \in  \wt F^\tau$, the $\tau$-stable weights in $W(\lag g, \lag t) \eta$ are $W(\lag g, \lag t)^\tau \eta$, so $\psi_\xi$ occurs in the expression \eqref{eq: local-expression-of-fd-T} of $\Theta(\sigma^+)$, and any other $\tau$-stable dominant weight $\eta$ of $F$ satisfies $\eta< \xi$.
	We can eliminate $\eta$ by the irreducible representation with the highest weight $\eta$, and finally obtain $\psi_\xi$ by induction, as desired.

	(3) This is because $\varphi \in \Smt (G\tau)$ is the restriction of a $\ag G$-invariant function on $\ag G \rtimes \la \tau \ra$.
	Since $t_s s\tau = g_s^{-1} \tau g_s$, we have $\varphi(s\tau \exp X) = \varphi(\Ad(g_s) (s\tau \exp X)) = \varphi(g_s t_s g_s^{-1} \cdot \tau \exp(\Ad(g_s)(X)))$.
	Then, \eqref{eq: local-expression-of-fd} follows from \eqref{eq: local-expression-of-fd-T}.
\end{pf}

\end{subequations}

\subsection{Local expression}\label{subsec: local-expression}
	In the previous subsection we explained the local expression of $\varphi \in \FC G$.
	In this subsection we aim to find a local expression for $\Theta \in \IC{G'\tau}$.

	Let $\pi^+$ be an irreducible representation of $G^+$ whose restriction to $G$ is irreducible.
	Then the infinitesimal character of $\pi^+$ can be represented by $\lambda \in (\lag t^*)^\tau$.
	For $s \in G'$, we have observed in Remark \ref{rmk: criteria-of-regular} that $\ag T_s = Z(\ag G, \lag g^{s\tau})$ is defined over $\R$.
	In Subsection \ref{subsec: twisted-Weyl-denominator} we introduced $Z_s = T_s^{s\tau, \circ}$, whose Lie algebra is $\lag z_{s, 0} = \lag t_{s, 0}^{s\tau}$.

\begin{fact}
	For any $s\in G'$, there exists a neighborhood $\mathcal U$ of $0 \in {\lag z}_{s,0}$ such that for each $X \in \mathcal U$,
\begin{equation}\label{eq: local-expression-Theta}
	\left( \Delta_{G^+}\Theta (\pi^+)\right) ( s\tau \exp X) 
	= \sum_{w \in W(\lag g, \lag t)^{\tau}} c_\lambda(s\tau, w) e^{\la w \lambda, X\ra}.
	\end{equation}
	where we omit $\Ad(g_s)$ as in the last paragraph of Subsection \ref{subsec: stable-inf-char}.
\end{fact}
\begin{pf}
\begin{subequations}
	This result essentially coincides with \cite[Theorem 3.2]{BC13}.
	Each $z \in Z(\lag g)^\tau$ is realized as an invariant differential operator on $G \tau$, and the radial part of $z$ on $Z_s s\tau \subset G \tau$ is calculated in \cite[Theorem 2.4.1]{Bouaziz87}, which takes the form $\Delta_{G^+} \circ \gamma(z) \circ \Delta_{G^+}^{-1}.$ Here  $\gamma: Z(\lag g) \to \C [ \lag t]^{W(\lag g, \lag t)}$ is the Harish-Chandra isomorphism, which  sends $z \in Z(\lag g)^\tau$ to $\gamma(z) \in \C[\lag z]^{W(\lag g, \lag t)^\tau} \subset \C[\lag t]$, and $\gamma(z) \in \C[\lag z]$ is regarded as an invariant differential operator on $Z_s$.
	Since $\pi^+$ has infinitesimal character represented by $\lambda \in (\lag t^*)^\tau$, we define $F(t) = \left( \Delta_{G^+}\Theta(\pi^+)\right) ( ts\tau )$ for $t \in Z_s \cap G's^{-1}$, and $F$ is then an eigenfunction for the algebra $\C[\lag z]^{W(\lag g, \lag t)^\tau}$ of differential operators on $Z_s$.
	In other words,
\[\forall z \in \C[\lag z]^{W(\lag g, \lag t)^\tau}, \quad 
	z F = \lambda(z) F.\]
	This differential equation has been studied by Harish-Chandra, and its solution has the following local expression
\[F( \exp X) = \sum_{w \in W(\lag g, \lag t)^{\tau} / W(\lag g, \lag t)_\lambda^\tau} P_w(X) e^{\la w \lambda, X \ra},\]
	when $X$ lies in a neighborhood $\mathcal U$ of $0 \in \lag z_{s, 0}$, where $W(\lag g, \lag t)^{\tau}_{\lambda}$ is the stabilizer of $\lambda$ in $W(\lag g, \lag t)^\tau$, and $P_w(X)$ are polynomials in $X$ of degree less than the cardinality of $W(\lag g, \lag t)^{\tau}_{\lambda}$; see \cite[Corollary of Lemma 27]{HC56} and \cite[Subsection 3.3]{BC13}.

	If $\lambda$ is regular, then the conclusion follows.
	As for singular $\lambda \in (\lag t^*)^\tau$, Bergeron-Clozel used a tensoring argument, due to Fomin-Shapovalov \cite{FS74}, to show that the polynomials $P_w(X)$ have degree zero.
	Although they only consider the case where $s \in T^{\tau, \circ} \cap G'$, their argument applies to any $s \in G'$.
\end{subequations}
\end{pf}

\begin{lem}
	For any $s\in M \cap G'$, there exists a neighborhood $\mathcal U$ of $0 \in {\lag z}_{s,0}$ such that for each $X \in \mathcal U$, we have
\begin{equation}\label{eq: local-expression-n}
	\left( \Delta_{G^+}\Theta_{\lag n}(\pi^+)\right) ( s\tau \exp X) 
	= \sum_{w \in W(\lag g, \lag t)^{\tau}} c_{\lag n, \lambda}(s\tau, w) e^{\la w \lambda, X\ra}.
	\end{equation}
\end{lem}
\begin{pf}
	Let $V$ denote the underlying Harish-Chandra module of $\pi^+$.
	According to \cite[Corollary 3.32]{HS83}, the (generalized) infinitesimal characters of $H_q({\lag n}, V)$ are given by $w \lambda + \rho_P$, where $w\in W({\lag g},{\lag t})$.
	When one restricts these characters to $M\tau$, only the factors with $\tau$-stable infinitesimal characters contribute.
	Due to Lemma \ref{lem: stable-infinitesimal-character}, the $\tau$-stable generalized infinitesimal characters of $H_q({\lag n}, V)$ are exhausted by $w \lambda + \rho_P$, $w\in W({\lag g}, \lag t)^{\tau}$.
	Then \eqref{eq: local-expression-n} follows from \eqref{eq: descent} and \eqref{eq: local-expression-Theta}.
\end{pf}

\subsection{Coherent family}\label{subsec: coherent-family}
	In this subsection we introduce the notion of coherent family.
Our definition is a natural generalization of \cite[Section 5]{BC13}.

	Recall that $\Lambda = X^*(\ag T)^\tau \subset ({\lag t}^*)^\tau$ in Lemma \ref{lem: fd-alg-rep}.
	Let $\Omega \subset G'$ be such that for each $s \in \Omega$, there is a neighborhood $U_s \subset Z_s$ of the identity such that $U_s s \subset \Omega$.
	We have defined the function space $\IF[\lambda] \Omega$ for each $\lambda \in (\lag t^*)^\tau$ in Subsection \ref{subsec: Grothendieck}.
	By a coherent family over $\Omega$ based at $\lambda_0 \in (\lag t^*)^\tau$, we mean a family of functions $\{ \Theta_\lambda \in \IF[\lambda]{\Omega} \mid \lambda \in \lambda_0 + \Lambda\}$ such that
\begin{itemize}
\item
	for each $\varphi \in \FC G$, if $\varphi$ has an expression as in \eqref{eq: local-expression-of-fd-T}, then for each $\lambda \in \Lambda + \lambda_0$,
\begin{equation}\label{eq: coherent}
	\varphi \Theta_\lambda = \sum_{\xi \in \Lambda} m_\xi \Theta_{\lambda + \xi}.
	\end{equation}
\item
	$\Theta_{w\lambda} = \Theta_\lambda$ for all $w \in W(\lag g, \lag t)^\tau$ that stabilizes $\lambda_0$.
\end{itemize}

\begin{fact}\label{fact: coherency}
	Suppose a family of functions $\{ \Theta_\lambda \in \IF[\lambda]{\Omega} \mid \lambda \in \lambda_0 + \Lambda\}$ has the local expression
\begin{equation}\label{eq: local-expression}
	\left( \Delta_{G^+}\Theta_\lambda \right) (s \tau \exp X) = \sum_{w \in W(\lag g, \lag t)^\tau} 
		c_\lambda(s\tau, w) e^{\la w\lambda, X \ra},
	\end{equation}
	and require $c_\lambda(s\tau, w) = c_\lambda(s\tau, v)$ whenever $w \lambda = v \lambda$.
	Then this family is coherent if and only if, for each $\xi \in \Lambda$
\begin{equation}\label{eq: coefficient-of-coherent-family}
	c_{\lambda_0+\xi}(s\tau, w) = (w\xi)(t_s)^{-1} c_{\lambda_0}(s\tau, w),
	\end{equation}
	where $t_s = g_s^{-1} \tau(g_s) s^{-1} \in \ag T_s$, and $g_s\in \ag G$ conjugates $(\ag B_s, \ag T_s)$ to $(\ag B, \ag T)$.
\end{fact}
\begin{pf}
	The proof is similar to that of \cite[Lemma 3.44]{HS83}.
	Given \eqref{eq: coefficient-of-coherent-family}, it is straightforward to verify the coherence from \eqref{eq: local-expression-of-fd}.
	Now suppose the family $\{\Theta_\lambda\}$ is coherent.
	Let $\xi \in \Lambda \subset (\lag t^*)^\tau$ be very anti-dominant with respect to $\lambda_0$.
	It is known from the proof of Lemma \ref{lem: very-anti-dominant}(1) that for each $w \in W(\lag g, \lag t)^\tau$, $\lambda_0+ w\xi$ is not in the $W(\lag g, \lag t)^\tau$-orbit of $\lambda_0+ \xi$ unless  $w\lambda_0 = \lambda_0$.
	Taking $\psi_\xi$ as in Lemma \ref{lem: fd-alg-rep}(2), one can easily calculate from \eqref{eq: local-expression} that $\psi_\xi \Theta_{\lambda_0}$ decomposes uniquely into a sum of $\Theta_{\lambda_0+ w \xi} \in \IF[\lambda_0+ w\xi]{G'}$ over $w \in W(\lag g, \lag t)^\tau_{\lambda_0} \backslash W(\lag g, \lag t)^\tau$.
	In particular, we can obtain $\Theta_{\lambda_0 + \xi}$ from $\psi_\xi \Theta_{\lambda_0}$.
	It follows from the definition of coherent family and \eqref{eq: local-expression-of-fd} that \eqref{eq: coefficient-of-coherent-family} holds for such a very anti-dominant weight $\xi$.
	To obtain \eqref{eq: coefficient-of-coherent-family} for general $\xi$,  we can take $\lambda \in \lambda_0 + \Lambda$ such that both $\lambda - \lambda_0$ and $\lambda - \lambda_0 - \xi$ are very anti-dominant.
	It follows from the previous argument that
\[c_\lambda(s\tau, w) = w(\lambda- \lambda_0)(t_s)^{-1}c_{\lambda_0}(s\tau, w) , \quad
	c_\lambda(s\tau, w) = w(\lambda- \lambda_0 - \xi)(t_s)^{-1}c_{\lambda_0 + \xi}(s\tau, w).\]
	Then we obtain \eqref{eq: coefficient-of-coherent-family} from the coherence.
\end{pf}

	As a corollary, each $\Theta \in \IF[\lambda_0]{\Omega}$ can be inserted into a unique coherent family over $\Omega$ based at $\lambda_0$.

\subsection{Coherent continuation}\label{subsec: coherent-continuation}
	For each $\pi^+ \in \GG[\lambda_0]{G^+}$, we know from Fact \ref{fact: coherency} that $\Theta(\pi^+)|_{M^-\tau}$ and $\Theta_{\lag n}(\pi^+)|_{M^-\tau}$ can be inserted into unique coherent families $\{\Theta_\lambda\}$ and $\{\Theta_{\lag n, \lambda}\}$, respectively, and the twisted Osborne conjecture $\Theta_{\lambda_0} = \Theta_{\lag n, \lambda_0}$ holds iff $\Theta_{\lambda} = \Theta_{\lag n, \lambda}$ for some $\lambda \in \lambda_0 + \Lambda$.
	In this subsection, we show that for a very anti-dominant $\lambda \in \lambda_0 + \Lambda$, $\Theta_{\lambda} = \Theta_{\lag n, \lambda}$ also follows from a twisted Osborne conjecture.

\begin{lem}\label{lem: reduction-by-coherent-continuation}
	Suppose $\lambda \in \lambda_0 + \Lambda$ is very anti-dominant.
	There exists $\pi_\lambda \in \GG[\lambda]{G^+}$ such that
\[\Theta(\pi_\lambda)|_{M^-\tau} = N_{\lambda_0} \Theta_\lambda, \quad
	\Theta_{\lag n}(\pi_\lambda)|_{M^- \tau} = N_{\lambda_0} \Theta_{\lag n, \lambda},\]
	where $N_{\lambda_0}$ is the cardinality of $W(\lag g, \lag t)_{\lambda_0}^\tau$.
\end{lem}
\begin{pf}
	Recall from Lemma \ref{lem: fd-alg-rep}(2) that for each $\xi \in \Lambda$, there exists a $\sigma^+ \in \FG{G^+}$ whose character $\psi_\xi = \Theta(\sigma^+)|_{G\tau}$ has local expression
\[\psi_\xi (\tau \exp X) = \sum_{w\in W(\lag g, \lag t)^\tau} e^{\la w \xi, X \ra} \]
	for $X \in \lag t_0^\tau$.
	Take the decomposition
\[ \sigma^+ \otimes \pi^+ = \sum_{\lambda'} \pi_{\lambda'} \in \bigoplus_{\lambda'} \GG[\lambda']{G^+}\]
	following \eqref{eq: decompose-of-Grothendieck}, where $\lambda'$ ranges over a set of representatives for $W(\lag g, \lag t) \rtimes \la \tau \ra$-orbits in $\lag t^*$.
	Then
\begin{equation}\label{eq: decomposition-of-tensor}
	\Theta(\sigma^+ \otimes \pi^+) = \sum_{\lambda'} \Theta(\pi_{\lambda'}),	\quad
	\Theta_{\lag n}(\sigma^+ \otimes \pi^+) = \sum_{\lambda'} \Theta_{\lag n}(\pi_{\lambda'}).
	\end{equation}
	On the other hand, we know from Lemma \ref{lem: F-linear-of-Theta-n} and the definition of coherent family that
\begin{equation}\label{eq: coherent-expansion-of-tensor}
	\begin{aligned}
	\Theta(\sigma^+ \otimes \pi^+) |_{M^- \tau} =	& 
	\psi_\xi \Theta(\pi^+)|_{M^-\tau} = 
	\sum_{w\in W(\lag g, \lag t)^\tau} \Theta_{\lambda_0 + w \xi},	\\
	\Theta_{\lag n}(\sigma^+ \otimes \pi^+) |_{M^- \tau} = 	&
	\psi_\xi \Theta_{\lag n}(\pi^+)|_{M^-\tau} = 
	\sum_{w\in W(\lag g, \lag t)^\tau} \Theta_{\lag n, \lambda_0 + w \xi},
	\end{aligned}
	\end{equation}
	Let $\xi$ be very anti-dominant with respect to $\lambda_0$, so that $\lambda: = \lambda_0 + \xi$ is very anti-dominant and $\lambda_0 + w \xi$, $\lambda_0 + v \xi$ belong to the same $W(\lag g, \lag t)^\tau$-orbit iff $v^{-1} w \in W(\lag g, \lag t)_{\lambda_0}^\tau$.
	Consider the term $e^{\la w\lambda, X \ra}$ in the local expressions of $\Delta_{G^+} \Theta(\sigma^+ \otimes \pi^+)$ and $\Delta_{G^+} \Theta_{\lag n}(\sigma^+ \otimes \pi^+)$ respectively.
	In \eqref{eq: decomposition-of-tensor}, this term comes from the terms $\Theta(\pi_{\lambda})$ and $\Theta_{\lag n}(\pi_{\lambda})$, whereas in \eqref{eq: coherent-expansion-of-tensor}, it comes from
\[\sum_{w\in W(\lag g, \lag t)^\tau_{\lambda_0}} \Theta_{\lambda_0 + w \xi} = N_{\lambda_0} \Theta_{\lambda}, \quad
	\sum_{w\in W(\lag g, \lag t)^\tau_{\lambda_0}} \Theta_{\lag n, \lambda_0 + w \xi} = N_{\lambda_0} \Theta_{\lag n, \lambda}.\]
	Then the conclusion follows.
\end{pf}

	As a corollary, to prove Theorem \ref{thm: main-new}, it suffices to consider the case where $\lambda$ is very anti-dominant.

\subsection{$\mu$-component}\label{subsec: mu-component}
	Recall that $A \subset M$ is the central split torus, so ${\lag a}$ lies in ${\lag t}_s$ for any $s\in M \cap G'$.
	Define $\mu : = \lambda_s |_{\lag a} \in ({\lag a}^*)^{s\tau} = ({\lag a}^*)^\tau$, which is independent of $s$.

	For $\Theta \in \IF[\lambda]{\Omega}$, we can define a function $\Theta_\mu$ by requiring
\[\left( \Delta_{M^+} \Theta_\mu \right)( s\tau \exp X) = 
	| \det_{\lag n}(s\tau)|^{\frac 1 2} \sum_{w \in W(\lag g, \lag t)^\tau, w\lambda|_{\lag a} = \mu} c(s\tau, w) e^{\la w \lambda + \rho_P, X \ra },\]
	when
\[\left( \Delta_{G^+}\Theta \right)(s \tau \exp X) = \sum_{w \in W(\lag g, \lag t)^\tau} c(s\tau, w) e^{\la w\lambda, X \ra}.\]
	It is uniquely determined by $\Theta$ due to the local independence of characters on $A^\tau$.
	Indeed, $|\det_{\lag n}|^{-\frac 1 2} \Delta_{M^+} \Theta_\mu$ is exactly the $\mu$-isotypic component of $\Delta_{G^+} \Theta$ with respect to the action of $A^\tau \subset Z_s$.
	Moreover, if $\Omega \subset G$ is a submanifold, then $\Theta$ is smooth on $\Omega$, and $\Theta_\mu$ is likewise smooth.

\begin{prop}\label{prop: reduction}
	For each $\lambda \in (\lag t^*)^\tau$ and $\pi^+ \in \GG[\lambda]{G^+}$, fix the following expression of characters:
\[\begin{aligned}
	\left( \Delta_{G^+} \Theta(\pi^+) \right) (s\tau \exp X) = &
	\sum_{w \in W(\lag g, \lag t)^\tau} c_\lambda(\pi^+; s \tau, w) e^{\la w \lambda, X \ra},	\\
	\left( \Delta_{G^+} \Theta_{\lag n}(\pi^+) \right) (s\tau \exp X) = &
	\sum_{w \in W(\lag g, \lag t)^\tau} c_{\lag n, \lambda}(\pi^+; s \tau, w) e^{\la w \lambda, X \ra},
	\end{aligned}\]
	where we require $c_\lambda (-, w) = c_\lambda ( -, v)$ and $c_{\lag n, \lambda} (-, w) = c_{\lag n, \lambda} (-, v)$ whenever $w \lambda = v\lambda$.
	Then the following are equivalent:

\begin{enumerate}[label=$(\arabic*)$]
\item
	for each $\lambda \in (\lag t^*)^\tau$ and $\pi^+ \in \GG[\lambda]{G^+}$, we have $\Theta(\pi^+)|_{M^-\tau} = \Theta_{\lag n}(\pi^+)|_{M^-\tau}$;

\item
	for each $\lambda \in (\lag t^*)^\tau$, $\pi^+ \in \GG[\lambda]{G^+}$ and $w \in W(\lag g, \lag t)^\tau$, we have $c_\lambda(\pi^+; -, w) |_{M^-\tau} = c_{\lag n, \lambda}(\pi^+; -, w) |_{M^-\tau}$;

\item
	for each  $\lambda \in (\lag t^*)^\tau$, $w \in W(\lag g, \lag t)^\tau$, and $\pi^+ \in \GG[w\lambda]{G^+}$, we have $c_{w\lambda}(\pi^+; -, 1) |_{M^-\tau} = c_{\lag n, w\lambda}(\pi^+; -, 1)|_{M^-\tau}$;

\item
	for each $\lambda \in \mathcal F$ and $\pi^+ \in \GG[\lambda]{G^+}$, we have
\begin{equation}\label{eq: reduction-components}
	\Theta(\pi^+)_\mu|_{M^- \tau} = \Theta_{\lag n}(\pi^+)_\mu|_{M^- \tau},
	\end{equation}
	where $\mu = \lambda|_{\lag a}$, and $\mathcal F$ is a subset of $(\lag t^*)^\tau$ that consists of very anti-dominant weights and satisfies $\mathcal F + \Lambda = (\lag t^*)^\tau$.
\end{enumerate}
\end{prop}
\begin{pf}
	Trivially, (1) is equivalent to (2), and (2) implies (3) and (4).
	We can also deduce (2) from (3) by the straightforward observation that $\GG[\lambda]{G^+} = \GG[w\lambda]{G^+}$, $c_\lambda(\pi^+; -, w) = c_{w\lambda}(\pi^+; -, 1)$, and $c_{\lag n, \lambda}(\pi^+; -, w) = c_{\lag n, w\lambda}(\pi^+; -, 1)$.

	We derive (3) from (4) by coherent continuation as follows.
	Let $\pi^+ \in \GG[\lambda]{G^+}$, where the representative $\lambda \in (\lag t^*)^\tau$ is not necessarily anti-dominant.
	Following Fact \ref{fact: coherency}, we can insert $\Theta(\pi^+)|_{M^-\tau}$ and $\Theta_{\lag n}(\pi^+)|_{M^-\tau}$ into unique coherent families $\{\Theta_{\lambda'} \}$ and $\{\Theta_{\lag n, \lambda'}\}$, respectively.
By Lemma \ref{lem: reduction-by-coherent-continuation}, when $\lambda' \in \lambda + \Lambda$ is very anti-dominant, there exists $\pi_{\lambda'} \in \GG[\lambda']{G^+}$ such that $\Theta_{\lambda'} = \frac 1 {N_{\lambda}} \Theta(\pi_{\lambda'})|_{M^-\tau}$ and $\Theta_{\lag n, \lambda'} = \frac 1{N_{\lambda}} \Theta_{\lag n}(\pi_{\lambda'})|_{M^-\tau}$.
	The condition (4) implies that $c_{\lambda'}(\pi_{\lambda'}; -, 1)|_{M^-\tau} = c_{\lag n, \lambda'}(\pi_{\lambda'}; -, 1)|_{M^-\tau}$ for $\lambda' \in \mathcal F \cap (\lambda + \Lambda)$, which is non-empty and consists of very anti-dominant weights.
	Then we obtain $c_\lambda(\pi^+ ; -, 1) = c_{\lag n, \lambda}(\pi^+ ; -, 1)$ from the coherence according to \eqref{eq: coefficient-of-coherent-family}.
	The original $\lambda$ need not be anti-dominant, so (3) is implied by (4).
\end{pf}

	In the above proposition, (1) is exactly the twisted Osborne conjecture stated in Theorem \ref{thm: main-new}, and (4) is what we will verify in the rest of this article; we will take $\mathcal F \subset (\lag t^*)^\tau$ to be the subset of very anti-dominant weights satisfying \eqref{eq: extra-on-lambda}.
	Moreover, it suffices to consider the case where $\pi^+$ is a representation of $G^+$ whose restriction to $G$ is irreducible.

\subsection{Proof outline}\label{subsec: outline}
	Let $\pi$ be a $\tau$-stable irreducible representation of $G$ with infinitesimal character represented by a very anti-dominant weight $\lambda \in (\lag t^*)^\tau$.
	Fix an extension $\pi^+$ of $\pi$ to $G^+$. We present an outline for our proof of \eqref{eq: reduction-components}, $\Theta(\pi^+)_\mu|_{M^-\tau} = \Theta_{\lag n}(\pi^+)_\mu|_{M^-\tau}$ where $\mu = \lambda|_{\lag a}$, in this subsection.

	Let $V$ be the underlying Harish-Chandra module of $\pi^+$.
	Recall the primary decomposition of $H_q({\lag n}, V)$ in \cite[(2.29)-(2.30)]{HS83}:
\[H_q ({\lag n}, V) = \bigoplus_{\nu \in {\lag a}^*} H_q({\lag n}, V)_\nu,\]
	where $H_q({\lag n}, V)_\nu$ is the $(\nu + \rho_P)$-primary subspace with respect to ${\lag a}$.
	We denote the Casselman-Wallach globalization of $H_q(\lag n, V)_\nu$ by $H_q(\lag n, \pi^+)_\nu$, and it is easy to check that the latter representation is the $(\nu| \det_{\lag n}|^{\frac 1 2})$-primary subspace of $H_q(\lag n, \pi^+)$ with respect to $A$.
The following lemma generalizes \cite[Lemma 6.1]{BC13} and \cite[Lemma 5.24]{HS83}.

\begin{lem}\label{lem: n-character}
	Since $\lambda$ is very anti-dominant, $\Theta_{\lag n}(\pi^+)_\mu = \Theta (H_0(\lag n, \pi^+)_\mu)$ on $M^-\tau$.
\end{lem}

\begin{pf}
	Recall from Subsection \ref{subsec: mu-component} that $|\det_{\lag n}|^{-\frac 1 2} \Delta_{M^+} \Theta_{\lag n}(\pi^+)_\mu$ is the $\mu$-isotypic component of $\Delta_{G^+} \Theta_{\lag n}(\pi^+)$, and that, by \eqref{eq: descent}, $\Delta_{G^+} = | \det_{\lag n}|^{-\frac 1 2} | D_{\lag n} | \Delta_{M^+}$. 
  Thus, $\Theta_{\lag n}(\pi^+)_\mu$ is the $(\mu| \det_{\lag n}|^{\frac 1 2})$-isotypic component of
\[|D_{\lag n}|\Theta_{\lag n}(\pi^+) = \sum_{q\geqslant 0} (-1)^q \Theta(H_q(\lag n, \pi^+) )\]
	with respect to $A^\tau$.
	We will prove that the following statements hold on $M^- \tau$:
\begin{enumerate}[label=$(\arabic*)$]
\item
	the $\mu$-isotypic component of $\Theta(H_q(\lag n, \pi^+))$ with respect to $A^\tau$ is $\Theta(H_q(\lag n, \pi^+)_\mu)$, and
\item
	$\Theta(H_q(\lag n, \pi^+)_\mu)$ vanishes when $q > 0$.
\end{enumerate}
	Then the conclusion follows.

	Only $\tau$-stable components of $H_q(\lag n, \pi^+)$ contribute to its character on $M\tau$, so we have the following equality
\[\Theta(H_q(\lag n, \pi^+)) |_{M'\tau}= \sum_{\nu \in ({\lag a}^*)^\tau} \Theta(H_q(\lag n,  \pi^+)_\nu) |_{M'\tau}.\]
	The restriction $({\lag a}^*)^\tau \to ({\lag a}^\tau)^*$ is an isomorphism, and  $\Theta(H_q(\lag n, \pi^+)_\nu)|_{M'\tau}$ is exactly the $\nu$-isotypic component of  $\Theta(H_q(\lag n, \pi^+))|_{M'\tau}$ with respect to $A^\tau$.
	Then (1) follows.

	We prove (2) by contradiction.
	If $\Theta(H_q(\lag n, \pi^+)_\mu)$ does not vanish for some $q> 0$, then $H_q({\lag n}, V)_\mu \not =0$.
	According to \cite[Proposition 2.32]{HS83}, there exists another $\nu \in {\lag a}^*$ such that $\nu < \mu$ and $H_0({\lag n}, V)_\nu \not =0$.
	Here ``$<$'' is the partial order on ${\lag a}^*$ defined in \cite[(2.31)]{HS83}; that is, $\mu - \nu$ is a non-zero linear combination of $\RS^+({\lag g}, {\lag a})$ with non-negative integer coefficients.
	However, due to Casselman-Osborne's lemma (see \cite[Corollary 3.32]{HS83}), $H_0({\lag n}, V)_\nu \not =0$ implies $\nu = w \lambda|_{\lag a}$ for some $w \in W({\lag g}, {\lag t})$.
	Since $\lambda$ is very anti-dominant, it follows from Lemma \ref{lem: very-anti-dominant} that $\nu - \mu = ( w\lambda - \lambda)|_{\lag a}$ is a linear combination of $\RS^+({\lag g}, {\lag a})$ with non-zero coefficients having positive real parts.
	Now we obtain a contradiction.
\end{pf}

	It would be more difficult to deal with $\Theta(\pi^+)_\mu$.

\begin{lem}\label{lem: n-homology}
	Suppose $H_0({\lag n}, V)_\mu \not =0$.
	Let $I = \Ind_P^G(H_0(\lag n, \pi^+)_\mu \otimes |\det_{\lag n}|^{-\frac 1 2 })$ and let $E$ be its underlying $({\lag g}, K)$-module.
	Since $\lambda$ is very anti-dominant, $V$ is the unique irreducible factor of $E$ such that  $H_0({\lag n}, -)_\mu \not =0$; in particular, $H_0({\lag n}, V)_\mu = H_0({\lag n}, E)_\mu$.
\end{lem}
\begin{pf}
	Let $W$ be an irreducible quotient of $H_0({\lag n}, V)_\mu$. Then its Casselman-Wallach globalization $\xi$ is an irreducible quotient of $H_0(\lag n, \pi^+)_\mu$.
	Let $I' = \Ind_P^G( \xi \otimes | \det_{\lag n} |^{-\frac 12})$, and let $E'$ be its underlying $({\lag g}, K)$-module.
	Due to the Frobenius reciprocity \cite[(4.11)]{HS83}
\[\Hom_G\left(V, E'\right) 
	\cong \Hom_{M}\left(H_0({\lag n}, V), W \right),\]
	the natural projection $H_0(\lag n, V) \twoheadrightarrow H_0(\lag n, V)_\mu \twoheadrightarrow W$ induces a non-zero map $V \to E'$.
	Since $V$ is irreducible, this map must be injective.

	Denote the Abelian category of Harish-Chandra modules of $G$ with generalized infinitesimal character $\lambda$ by $\HC_\lambda(G)$.
	By the anti-dominance of $\lambda$, the functors $H_q({\lag n}, -)_\mu$, $q>0$, vanish on $\HC_\lambda(G)$ as explained in the proof of the previous lemma; see also the argument in \cite[pages 120--121]{HS83}.
	Thus, $H_0({\lag n}, -)_\mu$ is exact on $\HC_\lambda(G)$.
	Then the injection $V \to E'$ induces an injection $H_0({\lag n}, V)_\mu \to H_0({\lag n}, E')_\mu$.
	According to \cite[Lemma 7.10]{HS83}, $H_0({\lag n}, E')_\mu = W$ is irreducible, so the non-zero module $H_0({\lag n}, V)_\mu$ is also irreducible.
	Then by definition, $W = H_0(\lag n, E')_\mu = H_0({\lag n}, V)_\mu$, and $\xi = H_0(\lag n, \pi^+)_\mu$ as an $M$-representation; in particular, $I'$ must be $\Ind_P^G (H_0(\lag n, \pi^+)_\mu \otimes | \det_{\lag n}|^{-\frac 1 2}) =: I$.
	We have verified that $H_0({\lag n}, V)_\mu$ coincides with $H_0({\lag n}, E)_\mu$, so $E/ V$ and all its subquotients have $H_0({\lag n}, -)_\mu =0$.
\end{pf}

	Note that the character of an admissible representation is essentially determined by its underlying Harish-Chandra module, and only those $\tau$-stable factors contribute to the twisted character.
	Then we can deduce $\Theta(\pi^+)_\mu = \Theta(I^+)_\mu$ or $0$ from the above lemma, once we generalize \cite[Proposition 6.2]{BC13} and \cite[Proposition 6.15]{HS83} as follows.
	Note that an extra condition on the infinitesimal character is needed.

\begin{prop}\label{prop: vanish}
	Suppose $H_0({\lag n}, V)_\mu =0$ and the very anti-dominant weight $\lambda$ satisfies an additional condition \eqref{eq: extra-on-lambda}. Then $\Theta(\pi^+)_\mu$ vanishes on $M^-\tau$.
\end{prop}

	We will prove this proposition in Section \ref{sec: asymtotic}, and show in Subsection \ref{subsec: extra-condition} that the condition \eqref{eq: extra-on-lambda} on $\lambda$ can also be satisfied by $\Lambda$-translation.
	Consequently, it suffices to calculate the twisted character of parabolically induced modules.
	Now we come to the generalization of \cite[Lemma 6.4]{BC13} and \cite[Lemma 5.17]{HS83}, which is established in Section \ref{sec: ind-character}.

\begin{prop}\label{prop: induce}
	Suppose $\xi$ is a $\tau$-stable irreducible Casselman-Wallach representation of $M$ with infinitesimal character represented by the very anti-dominant weight $\lambda \in (\lag t^*)^\tau$.
	Fix an intertwining operator $\tau_\xi : \xi \to \xi^\tau$ with $\tau_\xi^d = \Id$ such that $\xi$ and $I := \Ind_P^G(\xi)$ extend to $M^+$ and $G^+$, respectively.
	Then $\Theta(I^+)_\mu = \Theta(\xi^+ \otimes |\det_{\lag n}|^{\frac 1 2})$ on $M^- \tau$.
\end{prop}

\begin{pf}[ of Theorem \ref{thm: main-new}]
	Let $\mathcal F\subset (\lag t^*)^\tau$ consist of the very anti-dominant weights satisfying \eqref{eq: extra-on-lambda}.
	We will explain at the end of Subsection \ref{subsec: extra-condition} that $\mathcal F + \Lambda = (\lag t^*)^\tau$.
	Then, by Proposition \ref{prop: reduction}, it suffices to verify $\Theta(\pi^+)|_{M^- \tau} = \Theta_{\lag n }(\pi^+)|_{M^- \tau}$ for each extension $\pi^+$ of a $\tau$-stable irreducible $G$-representation $\pi$ with infinitesimal character $\lambda \in \mathcal F$.

	Let $\lambda \in \mathcal F$, $\pi^+$, and $\pi$ be defined as above.
	If $H_0(\lag n, \pi)_\mu= 0$, then we know from Proposition \ref{prop: vanish} that $\Theta(\pi^+)_\mu|_{M^- \tau} = 0$, and from Lemma \ref{lem: n-character} that $\Theta_{\lag n}(\pi^+)_\mu$, which equals $\Theta(H_0(\lag n, \pi^+)_\mu)$, also vanishes on $M^-\tau$.
	Suppose now that $\xi: = H_0(\lag n, \pi)_\mu \otimes |\det_{\lag n}|^{-\frac 1 2}$ does not vanish.
	According to the proof of Lemma \ref{lem: n-homology}, $\xi$ is irreducible and $\pi$ is the unique irreducible subrepresentation of $I: = \Ind_P^G(\xi)$;
  moreover, $H_0(\lag n, I/\pi)_\mu = 0$.
  The extension $\pi^+$ of $\pi$ naturally induces an extension of $\xi$ to $M^+$ and an extension $I^+$ of $I$ to $G^+$, and $\pi^+$ can be viewed as a subrepresentation of $I^+$.
  Only $\tau$-stable factors of $I/\pi$ contribute to its twisted character,
  but Proposition \ref{prop: vanish} shows that their $\mu$-components vanish on $M^-\tau$,
  so we have $\Theta(I^+/\pi^+)_\mu|_{M^- \tau} = 0$.
  It follows that $\Theta(\pi^+)_\mu|_{M^- \tau} = \Theta(I^+)_\mu|_{M^- \tau}$.
	Furthermore, as established in Proposition \ref{prop: induce}, we have $\Theta(I^+)_\mu|_{M^- \tau} = \Theta(\xi^+\otimes |\det_{\lag n}|^{\frac 1 2})_\mu|_{M^- \tau}$.
	Finally, from Lemma \ref{lem: n-character}, we find that $\Theta(\xi^+\otimes |\det_{\lag n}|^{\frac 1 2})_\mu = \Theta_{\lag n}(\pi^+)_\mu$.
	Now we obtain $\Theta(\pi^+)_\mu |_{M^-\tau} = \Theta_{\lag n}(\pi^+)_\mu |_{M^- \tau}$ in both cases.
\end{pf}
 
\newcommand{\df}{{\mathrm d}}
\newcommand{\supp}{\operatorname{supp}}
\newcommand{\Tg}{{\mathrm T}}

\section{Twisted character of induced module}\label{sec: ind-character}
	In this section we prove Proposition \ref{prop: induce}.
	Our argument will follow the approach of \cite[Lemma 5.17]{HS83}, with its Theorem 5.7 replaced by the twisted version \eqref{eq: regular-value}.

\subsection{Character as distribution}\label{subsec: parabolic-descent}
	Here we review the definition of parabolic induction, and present an integration formula due to Hirai \cite{Hirai68}.

	Let $P$ be a real parabolic subgroup of $G$, with unipotent radical $N$ and Levi component $M$.
	Let $\xi$ be an irreducible Casselman-Wallach representation of $M$, with underlying $({\lag m}, M \cap K)$-module $W$.
	Following the terminology of \cite[Introduction]{KV95}, we suppose that $W$ is also the underlying $({\lag m}, M \cap K)$-module of a continuous representation of $M$ on a Hilbert space.
	After averaging, $W$ is equipped with an inner product $\la-, - \ra_W$ such that the compact group $M \cap K$ acts on $W$ by unitary operators.

	Define the (normalized) parabolic induction $I = \Ind_P^G(\xi)$ as follows.
	For convenience, we use the Hilbert-space version of the representation.
The underlying $({\lag g}, K)$-module $E$ of $I$ consists of the functions $f: K \to W$ such that $f$ is finite under left translation by $K$, and $f(k m) = \xi(m)^{-1} f(k)$ for each $m \in M \cap K$, $k \in K$.
	Let $L^2(K, W)$ denote the completion of $E$ with respect to the inner product
\[\la f, f' \ra = \int_K \la f(k), f'(k) \ra_W \d k,\]
	and take this space as the representation space for $I$.
	To define the action of $g \in G$ on $L^2(K, W)$, we fix $(u(g), m(g), n(g)) \in K \times M \times N$ for each $g \in G$ such that $u(g) m(g) n(g) = g$, as in \cite[(11.21d)]{KV95}.
	Then $g \in G$ sends $f \in E$ to
\begin{equation}\label{eq: para-ind}
	(I(g) f)(k) := \left(\xi \otimes  |\det_{\lag n}|^{\frac 1 2} \right) (m(g^{-1} k))^{-1} f( u(g^{-1}k)).
	\end{equation}
	The uniqueness of the Iwasawa decomposition makes this action well defined.
	By \cite[Proposition 11.41]{KV95}, $I(g)$ extends to a bounded operator on $L^2(K, W)$.
	Thus, $I = \Ind_P^G(\xi)$ is a continuous representation of $G$ on a Hilbert space.
	It follows from \cite[(11.42)]{KV95} that any $K$-type occurs in $I$ with finite multiplicity, so the representation $I$ is admissible.

	The well-known Langlands classification (see \cite[Theorem 5.4.4]{Wallach88}, \cite[Theorem 11.198]{KV95}) asserts that any irreducible Harish-Chandra module can be realized as a submodule of the underlying Harish-Chandra module of some continuous representation on a Hilbert space.
	Thus we can define $\Ind_P^G(\xi)$ for an irreducible Casselman-Wallach representation $\xi$ of $M$.

	Since the representation $\xi$ in Proposition \ref{prop: induce} is $\tau$-stable, and we have fixed an intertwining operator $\tau_\xi: \xi \to \xi^\tau$, it naturally induces an intertwining operator $\tau_I: I \to I^\tau$ with $(\tau_If)(k) = \tau_\xi f(\tau^{-1}(k))$.
	Consequently, $\xi$ and $I$ both extend to $M^+ = M \rtimes \la \tau \ra$ and $G^+ = G \rtimes \la \tau \ra$, respectively, as $\xi^+$ and $I^+$.

	The following character formula relating $\Theta(I^+)$ and $\Theta(\xi^+)$ is a generalization of \cite[Theorem 1 in Section 4]{Hirai68}.
	For convenience, fix the Haar measures $\d k$, $\d m$, $\d n$ on $K$, $M$, $N$, respectively, and require the volume of $K$ to be $1$.
	By \cite[(11.21g), (11.21h)]{KV95}, there exists a Haar measure $\d g$ on $G$ such that for each $f \in \Smt_\cpt(G)$,
\begin{equation}\label{eq: Haar-measure-decompose}
	\int_G f(g) \d g =
	\int_K \d k \int_N \d n \int_M f(knm) \d m = 
	\int_K \d k \int_N \d n \int_M f(k m n) | \det_{\lag n}(m)|  \d m.
	\end{equation}
	The above Haar measures $\d g$ and $\d m$ can be extended to $G^+$ and $M^+$, respectively.
	Define the parabolic descent $\Smt_\cpt(G^+) \to \Smt_\cpt(M^+)$, $f \mapsto f_P$ by
\[f_P(m) =| \det_{\lag n}(m)|^{\frac 1 2} \int_K \d k \int_N f(k^{-1} mn k) \d n=
	| \det_{\lag n}(m)|^{-\frac 1 2} \int_K \d k \int_N f(k nm k^{-1}) \d n.\]

\begin{fact}\label{fact: Hirai}
	For each $f\in \Smt_\cpt(G^+)$,
\begin{equation}\label{eq: distribution}
	\int_{G^+}f(g) \Theta(I^+)(g) \d g=
	\int_{M^+} f_P(m) \Theta(\xi^+)(m) \d m.
	\end{equation}
\end{fact}
\begin{pf}
	See \cite[(7.1.5)]{Bouaziz87}.
	Hirai's original argument is valid for the disconnected group $G^+$.
\end{pf}

\subsection{Evaluation at regular points}
	Here we review Bouaziz's formula for evaluating the induced twisted character at regular elements.

	For $s \in G'$, $Z_s$ is defined as the identity component of $T_s^{s \tau}$ in Subsection \ref{subsec: twisted-Weyl-denominator}.
	Let $Z_s' = Z_s \cap G' s^{-1}$, and denote the union of $gZ_s' s\tau g^{-1}$ by $G[ Z_s'  s\tau]$.
	We will prove in Lemma \ref{lem: finite-cover} that $G[ Z_s'  s\tau]$ is open in $G\tau$.

	Following the notation in \cite[Subsection 7.1]{Bouaziz87}, take
\[\mathscr A_G^{M^+}\left( Z_ss \tau \right) =
	\{g \in G \mid g Z_s s\tau g^{-1} \subseteq M \tau\}.\]
	It admits a left translation by $M$, and a right translation by $N(G, Z_s s \tau)$.
	Take a set of representatives for $M \times N(G, Z_s s\tau)$-orbits in $\mathscr A_G^{M^+}\left( Z_s s \tau \right)$, and denote it by $\mathcal X$.
	We require the representative for the orbit $M N(G, Z_s s\tau)$ to be the identity $1 \in G$.
\begin{lem}\label{lem: G-conjugate-Cartan}
	Under the above notation,
\[G[Z_s' s\tau] \cap M^+ = \bigsqcup_{x \in \mathcal X} 
		M[ x Z_s's \tau x^{-1} ].\]
\end{lem}
\begin{pf}
	This is \cite[Lemma 7.1.2]{Bouaziz87} in the present context.
	The right-hand side is contained in the left-hand side by definition.
	Now let $t\in Z_s'$, so $ts \tau \in G^+$ is regular.
	Then for each $g \in G$, $g ts \tau g^{-1}$ is also regular.
	According to Lemma \ref{lem: adm-Cartan}, the $gts\tau g^{-1}$-admissible maximal torus in $\ag G$ is uniquely determined as $Z(\ag G, {\lag g}^{\Ad(g) (ts \tau)})$.
	If $g ts\tau g^{-1} \in M^+$, then from Lemma \ref{lem: tau-Borel}, we know that the $gts\tau g^{-1}$-admissible maximal torus $Z(\ag M, {\lag m}^{\Ad(g)(ts\tau)})$ in $\ag M$ coincides with the $gts\tau g^{-1}$-admissible torus in $\ag G$.
	In particular, $ {\lag m} \supseteq {\lag g}^{gts\tau g^{-1}}= \Ad(g) {\lag g}^{ts\tau} = \Ad(g) \lag z_s$, and then $\Ad(g)(Z_s ) \subseteq M$.
	Consequently, $\Ad(g)(Z_s s\tau) = \Ad(g)(Z_s) \cdot \Ad(g)(ts\tau) \subseteq M^+$, so $g$ lies in $\mathscr A_G^{M^+}\left( Z_ss \tau \right) $.
\end{pf}

	Denote $W_{G, s\tau} = N(G, Z_s s\tau)/ Z(G, Z_s s\tau)$.
	By \cite[Lemma 1.3.2]{Bouaziz87}, this is a finite group.
	It is also easy to deduce $Z(G, Z_s s\tau) = Z(G, Z_s)^{s\tau} = T_s^{s\tau}$ from Lemma \ref{lem: adm-Cartan}.
	For $x \in \mathcal X$, $x s\tau x^{-1} \in x Z_s s\tau x^{-1} \subseteq M\tau$ by definition, so we also have
\[W_{M, xs\tau x^{-1}} = \frac{N(M, xZ_s s\tau x^{-1})}
		{Z(M, xZ_s s\tau x^{-1})}.\]

\begin{fact}\label{fact: regular-value}
	Let $s\tau \in G' \tau$.
\begin{enumerate}[label=$(\arabic*)$]
\item
	If $\mathscr A_G^{M^+}(Z_s s\tau) = \varnothing$, then $\Theta(I^+)(s\tau) =0.$
\item
	If $s\tau \in M\tau$, then
\begin{equation}\label{eq: regular-value}
	\left( \Delta_{G^+} \Theta (I^+)\right) (s\tau) = \sum_{w \in W_{G, s\tau}} \sum_{x \in \mathcal X} 
		\frac 1 {|W_{M, xs\tau x^{-1}}|} 
		\left( \Delta_{M^+} \Theta(\xi^+) \right)( \left( \Ad(x) \circ w \right)(s\tau) )
	\end{equation}
\end{enumerate}
\end{fact}
\begin{subequations}
\begin{pf}
	This is the formula \cite[(7.1.13)]{Bouaziz87} adapted to the present context.
	Its proof needs three ingredients:
	Hirai's character formula \eqref{eq: distribution}, Lemma \ref{lem: G-conjugate-Cartan}, and the twisted Weyl integration formula \eqref{eq: twisted-Weyl-integration-formula} stated later.
	Consider the map $N \to N$ given by $n \mapsto n m n^{-1}m^{-1}$, where $m \in M^+$.
	It has differential $\Id_{\lag n}- \Ad(m)$, so for each $\gamma \in \Smt_\cpt(N)$,
\[\int_N \gamma(n) \d n = \left| D_{\lag n}(m) \right| \int_N \gamma( n m n^{-1} m^{-1}) \d n\]
	whenever $|D_{\lag n}(m)| \not=0$.
	Then for each $f\in \Smt_\cpt(G[Z_s's\tau])$, applying this change of variables to the parabolic descent, we obtain
\begin{equation}\label{eq: parabolic-descent}
	f_P(m) = | \det_{\lag n}(m)|^{-\frac 1 2}  \left| D_{\lag n}(m) \right|\int_K \d k \int_N f(kn mn^{-1} k^{-1}) \d n.
	\end{equation}
	In particular, $\supp(f_P)$ is included in $G[Z_s' s\tau] \cap M^+$.

	(1) If $s \tau \in G\tau$ satisfies $\mathscr A_G^{M^+}(Z_s s\tau) = \varnothing$, then $G[Z_s' s\tau] \cap M^+ = \varnothing$ by Lemma \ref{lem: G-conjugate-Cartan}.
	Consequently, for each $f \in \Smt_\cpt( G[Z_s' s\tau])$, we know from \eqref{eq: parabolic-descent} that $f_P$ vanishes on $M^+$, and from \eqref{eq: distribution} that $\int_{G^+} f \Theta(I^+ ) \d g= \int_{M^+}f_P \Theta(\xi^+) \d m = 0$.
	Thus, $\Theta(I^+)$ must vanish at $s\tau$.

	(2) Now let $s\tau \in M \tau$, and $f \in \Smt_\cpt( G[Z_s' s\tau])$.
	By Lemma \ref{lem: G-conjugate-Cartan}, the character formula \eqref{eq: distribution} becomes
\begin{equation}\label{eq: Hirai-formula}
	\int_{G\tau} f \Theta(I^+) \d g = \sum_{x \in \mathcal X} \int_{M[ x Z_s' s \tau x^{-1}]} f_P \Theta(\xi^+) \d m.
	\end{equation}
	For convenience, denote $\tau_x: = x s \tau x^{-1} = x s \tau(x)^{-1} \cdot \tau$; one verifies from the definition that $x Z_s s \tau x^{-1} = Z_{xs \tau(x)^{-1}} \tau_x$.
	According to the twisted Weyl integration formula \eqref{eq: twisted-Weyl-integration-formula} on $M\tau$, we have
\[\begin{aligned}
	\int_{M[x Z_s s\tau  x^{-1} ]} f_P \Theta(\xi^+) \d m = &
	\frac 1 {|W_{M, \tau_x}|} \int_{Z_{xs \tau(x)^{-1}}}\left( |D_{{\lag m}/ \lag z_{xs\tau(x)^{-1}}}| \Theta(\xi^+) \right)  (t \tau_x) \d t
	\int_{M/ xT_s^{s\tau} x^{-1} }f_P \left(\Ad(\bar m) (t \tau_x)  \right)  \d \bar m.
	\end{aligned}\]
	Note that
\[\begin{aligned}
	\int_{M/ x T_s^{s\tau} x^{-1}}f_P \left(\Ad(\bar m) (t \tau_x)  \right)  \d \bar m	=	& 
	 | \det_{\lag n}(t\tau_x )|^{-\frac 1 2}  \left| D_{\lag n}(t\tau_x) \right|
	\int_{M/ x T_s^{s\tau} x^{-1}}  \d \bar m \int_K \d k \int_N
	f \left(\Ad(k n \bar m) (t \tau_x)  \right)  \d n 	\\
	=	&
	 | \det_{\lag n}(t\tau_x )|^{-\frac 1 2}  \left| D_{\lag n}(t\tau_x) \right|
	\int_{G/ x T_s^{s\tau} x^{-1}} f \left(\Ad(\bar g) (t \tau_x)  \right) \d \bar g,	
	\end{aligned}\]
	where we have used \eqref{eq: parabolic-descent} and \eqref{eq: Haar-measure-decompose}, so
\[\begin{aligned}
	\int_{M[x Z_s s\tau  x^{-1} ]} f_P \Theta(\xi^+) \d m = 	&
	\frac 1 {|W_{M, x s \tau x^{-1}}|} \int_{xZ_s x^{-1}}\left( |D_{{\lag m}/ \lag z_{xs\tau(x)^{-1}}}| |\det_{\lag n}|^{-\frac 1 2} | D_{\lag n}| \Theta(\xi^+) \right)  (t \tau_x) \d t
	\int_{G/ x T_s^{s\tau} x^{-1}} f \left(\Ad(\bar g) (t \tau_x)  \right)  \d \bar g
	\\
	=	&
	\frac 1 {|W_{M, x s \tau x^{-1}}|} \int_{Z_s}\left( 
		|D_{{\lag m}/ \lag z_{xs\tau(x)^{-1}}}| |\det_{\lag n}|^{-\frac 1 2} | D_{\lag n}| \Theta(\xi^+) 
	\right)  (xt s\tau x^{-1}) \d t
	\int_{G/ T_s^{s\tau}} f \left(\Ad(\bar g) (t s\tau)  \right)  \d \bar g
	\end{aligned}.\]
	Since  $\Ad(x)$ sends ${\lag t}_s = \lag z_s \oplus \lag c_s$ to  ${\lag t}_{xs\tau(x)^{-1}}=  \lag z_{xs\tau(x)^{-1}} \oplus \lag c_{xs \tau(x)^{-1}}$, we know from \eqref{eq: det-q} and \eqref{eq: descent} that
\[\begin{aligned}
	\left( |D_{{\lag m}/ \lag z_{xs\tau(x)^{-1}}}| |\det_{\lag n}|^{-\frac 1 2} | D_{\lag n}| \right) (x t s \tau x^{-1}) =	&
	\left( |D_{\lag c_{xs\tau(x)^{-1}}}|\Delta_{M^+}^2 |\det_{\lag n}|^{-\frac 1 2} | D_{\lag n}| \right) (x ts \tau x^{-1}) =
	\left(|D_{\lag c_{xs\tau(x)^{-1}}}| \Delta_{M^+} \Delta_{G^+} \right) (x ts \tau x^{-1}) 	\\
	=	&
	\Delta_{M^+}(x ts \tau x^{-1}) \cdot \left( |D_{\lag c_s}| \Delta_{G^+} \right)(ts \tau) =
	\Delta_{M^+}(x ts \tau x^{-1}) \cdot \left(\Delta_{G^+}^{-1} |D_{\lag q_s}|  \right)(ts \tau) 
	\end{aligned}\]
	Thus, the right-hand side of \eqref{eq: Hirai-formula} equals
\begin{equation}\label{eq: through-m}
	\sum_{x \in \mathcal X} \frac 1 {|W_{M, x s \tau x^{-1}}|} 
	\int_{Z_s} \left( \Delta_{M^+} \Theta(\xi^+) \right)  (xt s\tau x^{-1}) \cdot (\Delta_{G^+}^{-1}|D_{\lag q_s} | )(ts \tau)\d t
	\int_{G/ T_s^{s\tau}} f \left(\Ad(\bar g) (t s\tau)  \right)  \d \bar g.
	\end{equation}
	On the other hand, applying the twisted Weyl integration formula \eqref{eq: twisted-Weyl-integration-formula} directly on $G\tau$, we obtain that the left-hand side of \eqref{eq: Hirai-formula} equals
\begin{equation}\label{eq: directly}
	\frac 1 {|W_{G, s\tau}|} \int_{Z_s} \Theta(I^+) (ts\tau) \cdot |D_{{\lag q}_s} (ts\tau)| \d t  \int_{G/ T_s^{s\tau}}
	f\left( \Ad(\bar g)(ts\tau)\right)  \d \bar g .
	\end{equation}
	From \cite[Lemma 1.5.2]{Bouaziz87}, we know that $f \mapsto \int_{G/T_s^{s\tau}} f(\Ad(\bar g)(ts\tau)) \d \bar g$ maps $\Smt_\cpt(G[Z_s's \tau])$ onto the $W_{G, s\tau}$-fixed subspace of $\Smt_\cpt(Z_s')$.
	Then we can deduce
\[\left( \Delta_{G^+}\Theta(I^+) \right)(t s\tau) = 
	\sum_{w \in W_{G, s\tau}} \sum_{x\in \mathcal X} \frac 1 {|W_{M, x s\tau x^{-1}}|} 
	\left({\Delta_{M^+}}\Theta(\xi^+)\right)\left (\Ad(x)(w( ts\tau))\right)\]
	by directly comparing \eqref{eq: through-m} and \eqref{eq: directly}.
\end{pf}
\end{subequations}

\subsection{$\mu$-isotypic component}
	For $s \in M^-$, take a neighborhood $\mathcal U$ of $0 \in {\lag z}_{s, 0}$ such that $\exp: {\lag z}_{s, 0} \to Z_s$ maps $\mathcal U$ diffeomorphically onto a neighborhood of $1 \in Z_s'$.
	We can expand  $\left( \Delta_{M^+} \Theta(\xi^+) \right)( \left( \Ad(x) \circ w \right) 
	\left( s\tau \exp X \right))$  in \eqref{eq: regular-value} for $X \in \mathcal U$ as
\[\sum_{v \in W({\lag m}, {\lag t}_{xs\tau(x)^{-1}})^{xs\tau x^{-1}}}
	d_{\lambda_x}(  (\Ad(x)w) ( s\tau) , v) 
	e^{\la v\lambda_x,  \Ad(x) w X \ra}.\]
	Here we need the following observations:
\begin{itemize}
\item
	by the definitions of $w\in W_{G, s\tau}$, $w(s\tau) \in Z_s' s\tau$, and $w$ preserves $Z_s$, ${\lag z}_s = {\lag t}_s^{s\tau}$;
\item
	by the definition of $x \in \mathcal X\subseteq \mathscr A_G^{M^+}(Z_s s\tau)$, we have that $\Ad(x)(s\tau) \in x Z_s' s\tau x^{-1} \subseteq M\tau$, $\Ad(x) {\lag t}_s^{s\tau} = \Ad(x) {\lag g}^{s\tau} = {\lag g}^{xs\tau x^{-1}} = {\lag t}_{x s\tau(x)^{-1}}^{xs \tau x^{-1}}$, $\Ad(x) Z_s = Z_{x s\tau (x)^{-1}}$, and $\Ad(x) Z_s' = Z_{x s\tau (x)^{-1}}'$;
\item
	$\lambda_x \in {\lag t}_{x s\tau(x)^{-1}}^*$ is chosen to be $\lambda_{x s\tau(x)^{-1}} = \lambda_s \circ \Ad(m_s^{-1} m_{x s\tau(x)^{-1}})$ under the notation at the end of Subsection \ref{subsec: stable-inf-char}.
\end{itemize}
	Recall that Proposition \ref{prop: induce} requires $\lambda \in (\lag t^*)^\tau$ to be very anti-dominant.
	To pick out the $\mu$-isotypic components, we have to determine those $(w, x, v)$ such that
\[\left(w^{-1} \Ad(x)^{-1} v \right) \lambda _ x | _{\lag a} = \mu.\]
\begin{lem}
	The above equation holds if and only if $w \in W_{M, s\tau}$ and $x = 1$.
\end{lem}
\begin{pf}
	Under the above notation, we have to find all $(w, x, v)$ such that
\begin{equation}\label{eq: test}
	\left( w^{-1} \Ad(x)^{-1} v \Ad(m_{x s\tau(x)^{-1}}^{-1}m_s ) \right) \lambda_s  |_{\lag a} = \lambda_s |_{\lag a}.
	\end{equation}
	Recall the definitions of $m_s$ at the end of Subsection \ref{subsec: stable-real-para}; we deduce that $\Ad(m_s) ({\lag t}_s^*)^{s\tau} = \Ad(m_s) ({\lag t}_s^*)^{\tau_s} = ({\lag t}^*) ^\tau$, and similarly $\Ad(m_{xs \tau(x)^{-1}}) ({\lag t}_{xs\tau(x)^{-1}}^*)^{xs\tau x^{-1}}= ({\lag t}^*)^\tau$.
	Since $v$ preserves $({\lag t}_{x s\tau(x)^{-1}}^*)^{x s\tau x^{-1}}$ by definition, and it is already observed that
\begin{itemize}
\item
	$\Ad(x)^{-1}$ sends $({\lag t}_{x s\tau(x)^{-1}}^*)^{x s\tau x^{-1}}$ to $ ({\lag t}_s^*)^{s\tau}$,
\item
	the actions of $w$ and $s\tau$ on ${\lag t}_s$ commute,
\end{itemize}
	we obtain that $w^{-1} \Ad(x)^{-1} v \Ad(m_{x s\tau(x)^{-1}}^{-1}m_s ) \in W(\ag G, \ag T_s)^{s\tau}$.

	Since $\lambda$ is very anti-dominant, \eqref{eq: test} implies that $w^{-1} \Ad(x)^{-1} v \Ad(m_{x s\tau(x)^{-1}}^{-1}m_s ) \in W(\ag M, \ag T_s)^{s\tau}$, according to Lemma \ref{lem: mu-component}(2).
	The rest of the argument is the same as the proof of \cite[Lemma 5.17]{HS83}.

	The actions of $v$ and $\Ad(m_{xs\tau(x)^{-1}}^{-1} m_s)$ on ${\lag a}$ are trivial by definition.
	Then the previous result implies that $\Ad(x) w$ acts trivially on ${\lag a}$.
	Take $\dot w \in N(G, Z_s s\tau)$ to be a lift of $w \in W_{G, s\tau}$; then $x \dot w \in G$ commutes with ${\lag a}$.
	According to Fact \ref{fact: centralizer-of-A}(1), $Z(G, A) = M$, so $x = (x \dot w) \cdot \dot w^{-1} \in M  N(G, Z_s s\tau)$.
	We have required that this orbit be represented by the identity element; thus $x = 1$.
	Consequently, $w$ acts trivially on ${\lag a}$, hence $w \in W_{M, s\tau}$.

	It is easy to see that $w \in W_{M, s\tau}$ and $x=1$ are sufficient for \eqref{eq: test}.
\end{pf}

	Now we know from \eqref{eq: regular-value} that the $\mu$-isotypic component of $\Delta_{G^+}\Theta(I^+)( s\tau \exp X)$ is
\[\begin{aligned}
	&\sum_{w \in W_{M, s\tau}} \frac 1 {|W_{M, s\tau}|} 
	\sum_{v \in W({\lag m}, {\lag t}_s)^{s\tau}} d_{\lambda_s}( w(s\tau), v) e^{\la v \lambda_s, w X \ra}	\\
	=	&
	\sum_{w \in W_{M, s\tau}} \frac 1 {|W_{M, s\tau}|} \left( \Delta_{M^+} \Theta(\xi^+) \right) ( w(s\tau)) =
	\left( \Delta_{M^+}\Theta(\xi^+) \right)(s\tau) .
	\end{aligned}\]
	Recall the definition in Subsection \ref{subsec: mu-component}; we deduce
\[\Theta(I^+)_\mu(s \tau) = 
	\left( \Delta_{M^+}^{-1} | \det_{\lag n} |^{\frac 1 2} \right)( s\tau) 
	\left( \Delta_{M^+} \Theta(\xi^+) \right) (s\tau)
	=\Theta( \xi^+ \otimes |\det_{\lag n}|^{\frac 1 2} )( s\tau).\]
	This is exactly the conclusion of Proposition \ref{prop: induce}.

\subsection{Twisted Weyl integration formula}
	The twisted Weyl integration formula is necessary for Fact \ref{fact: regular-value} and for the results in the next section.
	We state it in this subsection.

\begin{lem}\label{lem: finite-cover}
	Let $s\in G'$.
	Define $p: G/ T_s^{s\tau} \times Z_s' \to G\tau$, $(\bar g, t) \mapsto \Ad(\bar g)(ts\tau)$.
\begin{enumerate}[label=$(\arabic*)$]
\item
	$p$ is a local diffeomorphism.
	In particular, its image $G[Z_s' s\tau]$ is open in $G\tau$.
\item
	$p:  G/ T_s^{s\tau} \times Z_s' \to G[Z_s' s\tau]$ is a normal covering with finite deck transformation group $W_{G, s\tau}$.
\end{enumerate}
\end{lem}
\begin{pf}
	For (1), we calculate the tangent map of $p$ at $(\bar g, t) \in G/ T_s^{s\tau} \times Z_s'$ and show that it is an isomorphism.
	Identify $\Tg_t Z_s' \cong {\lag z}_{s, 0}$, $\Tg_{\bar g} (G/ T_s^{s\tau}) \cong {\lag q}_{s, 0}$, and $\Tg_{\Ad(\bar g) (ts\tau)} G\tau \cong {\lag g}_{0}$ as follows:
\begin{itemize}
\item
	let $Y\in {\lag z}_{s, 0}$ send $\psi \in \Smt\left( Z_s \right)$ to
\[Y \psi= \left. \frac \d {\d \varepsilon} \right|_{\varepsilon =0} \psi( t \exp \varepsilon Y);\]
\item
	fix a lift $g \in G$ of $\bar g$, regard $\gamma \in \Smt(G/T_s^{s\tau})$ as a right $T_s^{s\tau}$-invariant function on $G$, and let $Z \in {\lag q}_{s, 0}$ send $\gamma$ to
\[Z\gamma = \left. \frac{\d}{\d \varepsilon} \right|_{\varepsilon = 0}\gamma ( g \exp \varepsilon Z);\]
\item
	let $X \in {\lag g}_{0}$ send $f\in \Smt(G \tau)$ to
\[Xf = \left. \frac{\d}{\d \varepsilon} \right|_{\varepsilon = 0}
	f(g ts\tau g^{-1} \exp \varepsilon X).\]
\end{itemize}
	Then for each $Y \in {\lag z}_{s, 0}$, $p_* Y \in \Tg_{\Ad(\bar g)ts\tau}G\tau$ sends $f\in \Smt(G\tau)$ to
\[\begin{aligned}
	Y (f \circ p)=	&
	\left. \frac \d{\d \varepsilon} \right|_{\varepsilon =0} (f\circ p) ( \bar g, t \exp \varepsilon Y)	 
	=	
	\left. \frac \d{\d \varepsilon} \right|_{\varepsilon =0} f(g( t \exp \varepsilon Y) s\tau g^{-1})	\\
	=	&
	\left. \frac \d{\d \varepsilon} \right|_{\varepsilon =0} f(g ts \tau g^{-1} 
		\cdot \Ad(g)( \exp \varepsilon Y))	
	=	
	\left(\Ad(g) Y\right) f,
	\end{aligned}\]
	and for each $Z \in {\lag q}_{s, 0} \cong \Tg_{\bar g}(G/ T_s^{s\tau})$, $p_*Z \in \Tg_{\Ad(\bar g)ts\tau}G\tau$ sends $f\in \Smt(G\tau)$ to
\[\begin{aligned}
Z( f \circ p) =	&
	\left. \frac \d {\d \varepsilon} \right|_{\varepsilon =0} (f\circ p) ( g  \exp \varepsilon Z, t)=	
	\left. \frac \d {\d \varepsilon} \right|_{\varepsilon =0} f\left (
		g  \cdot \exp \varepsilon Z \cdot 
		ts\tau \cdot \exp (-\varepsilon Z) \cdot g^{-1}
		\right)	\\
	=	&
	\left. \frac \d {\d \varepsilon} \right|_{\varepsilon =0} f\left (
		g ts\tau g^{-1}  \cdot  g (ts\tau)^{-1}\exp \varepsilon Z \cdot 
		ts\tau \cdot \exp (-\varepsilon Z) \cdot g^{-1}
		\right)	\\
	=	&
	\left. \frac \d {\d \varepsilon} \right|_{\varepsilon =0} f\left (
		\Ad(\bar g)(ts\tau)  \cdot  \Ad(g) \left( 
			\Ad(ts\tau)^{-1}\left(\exp \varepsilon Z\right) \cdot 
			\exp (-\varepsilon Z)\right)
		\right)	=
	\left( \Ad(g) \left( \Ad(ts\tau)^{-1} Z - Z  \right)\right)f.
	\end{aligned}\]
	In conclusion, the tangent map of $p$ at $(\bar g, t)$ is the composition
\[\begin{tikzcd}
	{{\lag z}_{s, 0}\oplus {\lag q}_{s, 0}} \arrow[r] & 
	{\lag g}_0 \arrow[r, "\Ad(g)"] & 
	{\lag g}_0 \\
	{(Y, Z)} \arrow[r, maps to]           & 
	Y+(\Ad(ts\tau)^{-1} -1)Z &     
	\end{tikzcd}\]
	The above formula coincides with the formulas in \cite[(50)]{Mezo13}, \cite[(4.45)]{Knapp96}, and \cite[Proposition 2.4]{Renard97}.
	Since $t \in Z_s' = Z_s \cap G's^{-1}$, $ts\tau$ is regular in $G\tau$, so the above map is an isomorphism.

	To prove (2), it suffices to show that for each $t \in Z_s'$, $p^{-1}(ts\tau) \subset G/T_s^{s\tau}$ is in bijection with $W_{G, s\tau} = N(G, Z_s s\tau)/ T_s^{s\tau}$.
	Suppose $g t's\tau g^{-1} = t s \tau$. Then ${\lag g}_0^{ts\tau} = {\lag g}_0^{\Ad(g)(t's\tau)} = \Ad(g) {\lag g}_0^{t's\tau}$.
	Due to the regularity, ${\lag g}_0^{ts\tau}$ and ${\lag g}_0^{t's\tau}$, both of which contain $\lag z_{s, 0}= {\lag g}_0^{s\tau}$, coincide.
	Thus, $\Ad(g) Z_s = Z_s$, and $g Z_s s\tau g^{-1} = g Z_s g^{-1} \cdot g t' s\tau g^{-1}= Z_s s\tau$, so $g \in N(G, Z_s s\tau)$.
	Conversely, for each $g \in N(G, Z_s s\tau)$, $(g, g^{-1} t \Ad(s\tau)(g))$ lies in $p^{-1}(ts\tau)$.
	Consequently, we obtain $p^{-1}(ts\tau) \cong W_{G, s\tau}$.
\end{pf}

	The following equality \eqref{eq: twisted-Weyl-integration-formula} is stated in \cite[Proposition 1 in Subsubsection 5.4.1]{Mezo13}.
\begin{lem}
Choose the invariant measures on $G$, $T_s^{s \tau}$ and $G/ T_s^{s \tau}$ such that for each $f \in \Smt_\cpt(G^+)$,
\[\int_G f(g) \d g =  \int_{T_s^{s \tau}}\d t \int_{G/ T_s^{s \tau}} f(g t)  \d \bar g.\]
	Then for each $f \in \Smt_\cpt(G^+)$ supported in $G[Z_s' s\tau]$,
\begin{equation}\label{eq: twisted-Weyl-integration-formula}
	\int_{G\tau} f(g) \d g = 
	\frac 1 {|W_{G, s\tau}|} \int_{Z_s} | D_{{\lag q}_s}(ts \tau)| \d t
	\int_{G/ T_s^{s\tau}} f( \Ad(\bar g)(ts \tau)) \d \bar g  .
	\end{equation}
	Here ${\lag q}_s$ is defined in Lemma \ref{lem: twisted-Weyl}(2), and $D_{{\lag q}_s}(ts\tau) = \det_{{\lag q}_s}(1 - ts\tau)$ in the notation of Subsection \ref{subsec: notation}.
\end{lem}
\begin{pf}
	We have shown in Lemma \ref{lem: finite-cover} that $p$ is a $|W_{G, s\tau}|$-sheet covering.
	By \cite[(8.59) and Proposition 8.19]{Knapp96},
\[\int_{G\tau} f \d g = \frac 1 {|W_{G, s\tau}|} \int_{G/ T_s^{s\tau} \times Z_s'} \left( p^*f  \right) \cdot \left( p^*\df g \right),\]
	where $(p^*f)(\bar g, t) = f( \Ad(\bar g)(ts\tau))$, so it suffices to verify $p^*\df g = |D_{{\lag q}_s}^{s\tau}(t)| \d \bar g \d t$.

	We have identified the tangent spaces at $\bar g \in G/ T_s^{s\tau}$, $t\in Z_s'$, $p(\bar g, t) \in G\tau$ with ${\lag q}_{s, 0}, {\lag z}_{s, 0}, {\lag g}_0$ respectively.
	The canonical decomposition ${\lag g}_0 \cong {\lag z}_{s, 0} \oplus {\lag q}_{s, 0}$ makes the invariant measures on $G/ T_s^{s\tau}$, $T_s^{s\tau}$, $G\tau$ compatible as required.
	Then the determinant of the tangent map $p_*: {\lag z}_{s, 0} \oplus {\lag q}_{s, 0} \to {\lag g}_0$ at $(\bar g, t)$ is $\det_{\lag g}(g) \cdot \det_{{\lag q}_s}( (ts\tau)^{-1} -1 )$.
	Since $\Ad(g) \in \GL({\lag g})$ lies in the semisimple subgroup $\Ad({G})$, its determinant must have finite order.
	Thus, $| \det_{{\lag g}}(g)| =1$.
	Similarly, $|\det_{{\lag q}_s}( (ts\tau)^{-1} -1 )| = |\det_{{\lag q}_s}(1 - ts\tau)|$.
	Now the conclusion follows.
\end{pf}

\begin{rmk}
	In the next section we will use a slightly different version: choose the invariant measure as in the previous lemma; then for each $F \in \Smt_\cpt(G/T_s^{s\tau} \times Z_s')$ and $f \in \Smt_\cpt(G[Z_s' s\tau])$, we have
\begin{equation}\label{eq: pull-back-int}
	\int_G f\cdot  p_*F \d g=
	\int_{Z_s'} |D_{{\lag q}_s}^{s\tau}| \d t
	\int_{G/T_s^{s\tau}} p^*f  \cdot F \d \bar g .
	\end{equation}
	It is deduced as follows.
	The above formula \eqref{eq: twisted-Weyl-integration-formula} implies that
\[\int_G f\cdot  p_*F \d g=
	\frac 1 {|W_{G, s\tau}|} \int_{Z_s'} |D_{{\lag q}_s}^{s\tau}| \d t
	\int_{G/T_s^{s\tau}} p^*\left( f  \cdot p_*F \right) \d \bar g .\]
	Note that $p^*(f \cdot p_* F) = p^*f \cdot p^*p_*F$, and
\[(p^*p_*F)(\bar g, t) = \sum_{w \in W_{G, s\tau}} F(\bar g w^{-1}, w(t))\]
	is the sum of $F$ over deck transformations in $W_{G, s\tau}$.
	Then \eqref{eq: pull-back-int} follows, since the integration over $G/T_s^{s\tau} \times Z_s'$ is invariant under the deck transformations.
\end{rmk}
 
\newcommand{\dual}[1]{{#1}^\vee}

\section{Asymptotics of twisted character}\label{sec: asymtotic}
	In this section we will prove Proposition \ref{prop: vanish}, namely that $\Theta(\pi^+)_\mu$ vanishes on $M^-\tau$ under certain conditions.
	Our argument follows the approach in \cite{HS83} and \cite{BC13}.

	Take
\[ {\lag a}^{\tau,-}_0 = 
	\{ X \in {\lag a}^\tau_0 \mid \la \beta_\res, X \ra < 0, \forall \beta_\res \in \RS^+({\lag g}, {\lag a}^\tau) \},\]
	and set $A^{\tau, -} = \exp {\lag a}^{\tau, -}_0$. Then it is easy to check that $M^- A^{\tau, -} \subseteq M^-$.
	For $s \in M^-$, it is easy to deduce from \eqref{eq: local-expression-Theta} that
\[\left( \Delta_{G^+}\Theta(\pi^+)\right) ( m s\tau \exp X) 
	= \sum_{i=1}^r \theta_i(ms\tau) e^{\la \nu_i ,  X\ra},\]
	where $X \in {\lag a}_0^{\tau, -}$, $ms\tau$ lies in a neighborhood of $s\tau \in M^-\tau$, each $\theta_i$ is smooth in $m$, and $\mu = \nu_1, \cdots, \nu_r$ are representatives of $W({\lag g}, {\lag t})^{\tau} \lambda|_{{\lag a}^\tau}$.

	If $\Theta(\pi^+)_\mu$ does not vanish at some $s\tau \in M^-\tau$, then from the definition $\theta_1 = |\det_{\lag n} |^{-\frac 1 2} \Delta_{M^+} \Theta(\pi^+)_\mu|_{M^- \tau}$, we know that $\theta_1(s\tau) \not =0$.
	Hence, the main part of  $\left( \Delta_{G^+}\Theta(\pi^+)\right) (s \tau \exp X)$ would be $e^{\la \mu, X \ra}$ as $X \to -\infty$.
	Let us make it more precise.
	Take a neighborhood $U \subseteq Z_s \cap M^- s^{-1}$ of $1 \in Z_s$, and a function $\psi \in \Smt_\cpt(U)$ such that $\int_U \theta_1(ts\tau) \psi(t) \d t =1$.
	Define
\[\Psi(\exp X):= e^{-\la \mu, X \ra} \int_U \left(\Delta_{G^+} \Theta(\pi^+) \right)(t s \tau \exp X) \psi(t) \d t = 1+ \sum_{i=2}^r c_i e^{\la \nu_i - \mu, X \ra}.\]
	According to Lemma \ref{lem: very-anti-dominant}(1), each non-zero $\nu_i - \mu \in ({\lag a}^{\tau})^*$ is a linear combination of $\beta_\res \in \RS^+({\lag g}, {\lag a}^\tau)$ with coefficients having positive real parts.
	Choose a sequence $\{a_k = \exp X_k \}\subset A^{\tau, -} = \exp {\lag a}_0^{\tau, -}$ such that for every $\beta_\res \in \RS^+({\lag g}, {\lag a}^\tau)$, $\la \beta_\res, X_k \ra \to - \infty $ as $k \to \infty$.
	Then for every $i =2, \cdots, r$, we have $e^{\la \nu_i - \mu, X_k \ra} \to 0$ as $k \to \infty$.
	Hence, $\Psi(a_k) \to 1$.

	On the other hand, by analyzing the matrix coefficient functions of $\pi$, we can obtain another asymptotic estimate for $\left( \Delta_{G^+}\Theta(\pi^+)\right) (s \tau \exp X)$; this is Lemma \ref{lem: asymptotic-restriction}.
	Its proof occupies the bulk of this section.
	Through this approach, the condition of Proposition \ref{prop: vanish} would imply that $\left( \Delta_{G^+}\Theta(\pi^+)\right) ( t s \tau \exp X)$ is strictly dominated by $e^{\la \mu, X\ra}$, as a corollary of Lemma \ref{lem: exponents-of-V}.
	Then we obtain a contradiction, as presented in Subsection \ref{subsubsec: extract}.

\subsection{$\tau$-stability of Langlands embedding}
	The goal of this subsection is to establish the ``$\tau$-stability'' of the Langlands data of $\pi$, which is necessary for estimating matrix coefficients of $\pi$.
	Fix an Iwasawa decomposition $G = K A_\min N_\min$. Then there exists a minimal parabolic subgroup $P_\min = L_\min A_\min N_\min$ with $L_\min = Z(K, A_\min)$.
	This decomposition need not be stable under $\tau$, but $G = K \tau(A_\min) \tau(N_\min)$ is still an Iwasawa decomposition.
	By the uniqueness of Iwasawa decomposition (see, for example, \cite[Theorem 6.51 and Corollary 6.55]{Knapp96}), there exists $k_\tau \in K$ such that $\Ad(k_\tau)(A_\min N_\min) = \tau^{-1}(A_\min N_\min)$.

	According to the Langlands classification \cite[Theorem 5.4.4]{Wallach88}, $V$ is the unique irreducible submodule of (the underlying $({\lag g}, K)$-module of) a standard continuous-series representation $\bar I : =  \Ind_{\bar P}^G(\eta \boxtimes \nu)$, where
\begin{itemize}
\item
	$\bar P = \bar L \bar A \bar N$ is a real parabolic subgroup of $G$ satisfying $\bar P \supseteq P_\min$, $\bar L \bar A\supseteq L_\min A_\min$, and $\bar N \subseteq N_\min$.
\item
	$\nu\in \bar {\lag a}^*$ is regarded as a character of $\bar A$, satisfying $\la \Re \nu, \beta \ra < 0$ for each $\beta \in \RS^+({\lag g}, \bar {\lag a})$, and
\item
	$\eta$ is an irreducible tempered representation of $\bar L$.
\end{itemize}
	For convenience, we denote the Levi subgroup $\bar L \bar A$ of $\bar P$ also by $\bar M$, and its representation $\eta \boxtimes \nu$ simply by $\xi$ if $\nu$ is clear.
	Since $\pi$ is $\tau$-stable, the standard continuous-series representation $\bar I$ is also $\tau$-stable by the uniqueness of Langlands classification.
	Let us make it clear.

\begin{lem}\label{lem: extend-standard}
	The standard continuous-series representation $\bar I$ can be extended to $G^+ = G \rtimes \la \tau \ra$.
\end{lem}
\begin{pf}
	We will provide a series of isomorphisms
\[\begin{tikzcd}
	\bar I^\tau & 
	\Ind_{\tau^{-1}(\bar P)}^G(\xi^\tau) \arrow[l, "\tau"'] & 
	\Ind_{\bar P_\tau}^G(\xi_\tau) \arrow[l, "k_\tau"'] & 
	\Ind_{\bar P}^G(\xi)  = \bar I \arrow[l, "\tau_\xi"'] ,
	\end{tikzcd}\]
	where
\begin{itemize}
\item
	$\bar I^\tau$ has the same representation space as $\bar I$, while $g \in G$ acts by $\bar I^\tau(g) : = \bar I (\tau(g))$
\item
	$\Ind_{\tau^{-1}(\bar P)}^G(\xi^\tau)$ is the parabolic induction of the $\tau^{-1}(\bar M)$-representation $\xi^\tau: = \xi \circ \tau$ with respect to $\tau^{-1}(\bar P)$, and
\item
	$\bar P_\tau = \Ad(k_\tau)^{-1} \tau^{-1} \bar P \supseteq P_\min$, $\xi_\tau = \xi^\tau \circ \Ad(k_\tau)$.
\end{itemize}
	We have explained the representation space of $\bar I$ in Subsection \ref{subsec: parabolic-descent}.
	The isomorphism $\tau:  \Ind_{\tau^{-1}(\bar P)}^G(\xi^\tau) \to \bar I^\tau$ sends $f$ to $\tau f: k \mapsto f (\tau^{-1}(k))$.
	The isomorphism $k_\tau: \Ind_{\bar P_\tau}^G (\xi_\tau) \to \Ind_{\tau^{-1}(\bar P)}^G(\xi^\tau) $ sends $f$ to $k \mapsto f(k k_\tau)$.
	The isomorphism $\tau_\xi: \Ind_{\bar P}^G(\xi) \to \Ind_{\bar P_\tau}^G (\xi_\tau)$ comes from \cite[Theorem 5.4.1(3)]{Wallach88}: since $V \stackrel{\tau_\pi} \to V \hookrightarrow 
	\bar I^\tau \stackrel{\sim} \to 
	\Ind_{\tau^{-1}(\bar P)}^G(\xi^\tau) \stackrel{\sim}\to
	\Ind_{\bar P_\tau}^G(\xi_\tau)$ is also a Langlands embedding with respect to $P_\min = L_\min A_\min N_\min$, $\bar P_\tau$ must coincide with $\bar P$, and $\xi_\tau$ must be unitarily equivalent to $\xi$.
	Then it is easy to induce $\tau_\xi:  \Ind_{\bar P}^G(\xi) \to \Ind_{\bar P_\tau}^G(\xi_\tau)$ from $\tau_\xi: \xi \cong \xi_\tau$.

	Now let $\tau_{\bar I}: \bar I \to \bar I^\tau$ send $f$ to $\tau_{\bar I} f : k \mapsto \tau_\xi f(\tau^{-1}(k) k_\tau)$. Then it fits into the following commutative diagram:
\[\begin{tikzcd}
	\bar I \arrow[d, "\tau_{\bar I}"']  \arrow[r, "\tau_\xi"] & 
	\Ind_{\bar P_\tau}^G(\xi_\tau)   \arrow[d, "k_\tau"] \\
	\bar I^\tau         &
	\Ind_{\tau^{-1}(\bar P)}^G(\xi^\tau)  \arrow[l, "\tau"']    
	\end{tikzcd}.\]
	Moreover, $\left( \tau^{d}_{\bar I}f \right)(k) = \tau_\xi^{d} f(k \cdot \tau^{-(d-1)}(k_\tau) \cdots k_\tau)$.
	Since $\Ad(\tau^{-(d-1)}(k_\tau) \cdots k_\tau) = \Ad(\tau(k_\tau) \cdots \tau^d (k_\tau) ) = \left( \tau \Ad(k_\tau) \right)^d$ preserves the Iwasawa decomposition $G = K A_\min N_\min$, it follows from \cite[Corollary  6.55 and Theorem 6.57]{Knapp96} that $\tau^{-(d-1)}(k_\tau) \cdots k_\tau =: l_\tau \in Z(K, A_\min)  = L_\min$.
	Thus, we have $\left( \tau^{d}_{\bar I}f \right)(k) = \tau_\xi^{d} \cdot \xi(l_\tau)^{-1} f(k)$.
	Since $\tau_\xi^d$ intertwines the irreducible representation $\xi$ with $\xi\circ \left( \tau \circ \Ad(k_\tau) \right)^d = \xi\circ \Ad(l_\tau )$, it differs from $\xi(l_\tau)$ by a scalar.
	Then $\tau_{\bar I}^d$ is also a scalar.
	Now it is straightforward to scale $\tau_\xi$ such that $\tau_{\bar I}^d = \Id_{\bar I}$ and $\tau_{\bar I}|_{\pi} = \tau_\pi$.
	We obtain the conclusion.
\end{pf}

	It follows from the above argument that $\tau\Ad(k_\tau)$ also preserves $\bar M$, $\bar A$, $\bar N$, and the central character $\nu \in \bar {\lag a}^*$ of $\xi$.
	Moreover, according to \cite[(11.8) and (11.20)]{KV95}, $\bar L$ is uniquely determined by $\bar P = \bar L \bar A \bar N$, so it is stable under $\tau \Ad(k_\tau)$.

	We give an alternative interpretation of the ``$\tau$-stability'' of $\nu$.
	Since $\tau \Ad(k_\tau)$ is a semisimple automorphism of $\ag G$ that commutes with $\theta$, after an argument similar to that in Subsection \ref{subsec: stable-real-para}, we obtain a Borel pair $(\bar {\lag b}, \bar {\lag t})$ in ${\lag g}$ such that
\begin{itemize}
\item
	$(\bar {\lag b}, \bar {\lag t})$ is stable under $\tau \Ad(k_\tau)$, $\bar {\lag t}$ is real and stable under $\theta$, and
\item
	$(\bar {\lag b},\bar {\lag t})$ is compatible with $\bar P = \bar M \bar N$, in the sense that $\bar {\lag b} \subseteq \bar {\lag p}$, $\bar {\lag t} \subseteq \bar {\lag m}$, and hence $\bar {\lag b} \supseteq \bar {\lag n}$, $\bar {\lag t} \supseteq \bar {\lag a}$.
\end{itemize}
	Let us fix $g \in \ag G(\C)$ such that $\Ad(g) (\bar {\lag b}, \bar {\lag t}) = ({\lag b}, {\lag t})$.
	Recall that we have fixed a very anti-dominant $\lambda \in {\lag t}^*$ that represents the infinitesimal character of $\pi$ in Subsection \ref{subsec: very-anti-dominance}.
\begin{lem}\label{lem: nu-stable-tau}
	There is a $w \in W({\lag g}, {\lag t})^\tau$ such that $(w \lambda) \circ \Ad(g)$ coincides with $\nu$ on $\bar {\lag a}$.
\end{lem}
\begin{pf}
	Let $\bar \lambda = \lambda \circ \Ad(g)$. We will find $w \in W({\lag g}, \bar {\lag t})^{g^{-1} \tau g}$ such that $w \bar \lambda|_{\bar {\lag a}} = \nu$ as follows.
	Since $g^{-1}\tau^{-1}(g) k_\tau = \left( g^{-1} \tau g \right)^{-1} \cdot \tau k_\tau$ is an element of $\ag G(\C)$ stabilizing $(\bar {\lag b}, \bar {\lag t})$, it must act trivially on $\bar {\lag t}$.
	In particular, $W({\lag g}, \bar {\lag t})^{g^{-1} \tau g} = W({\lag g}, \bar {\lag t})^{\tau k_\tau}$.
	Since $\xi \cong \xi \circ \tau\Ad(k_\tau)$ as representations of $\bar M$, their infinitesimal character can be represented by $\bar \lambda' \in (\bar {\lag t}^*)^{\tau k_\tau}$, with $\bar \lambda' |_{\bar {\lag a}} = \nu$; here we have used Lemma \ref{lem: stable-infinitesimal-character}.
	According to \cite[Proposition 11.43]{KV95}, $\bar \lambda'$ also represents the infinitesimal character of $\bar I$.
	Since $\bar \lambda'$ and $\bar \lambda \in \bar {\lag t}^*$ are both stable under $g^{-1} \tau g$, it follows from Lemma \ref{lem: stable-infinitesimal-character} again that $\bar \lambda' = w \bar \lambda$ for some $w \in W({\lag g}, \bar {\lag t})^{g^{-1} \tau g}$.
	Now the conclusion follows.
\end{pf}

\subsection{Decomposition into types}
	Decompose $V$ into its $K^+= K \rtimes \la \tau \ra$-isotypic components
\[V = \bigoplus_{\sigma{} \in \Irr(K^+)} V_{\sigma{}};\]
	then $\Theta(\pi^+)$ is a sum of matrix coefficient functions for each $V_{\sigma{}}$.
	To be more precise, take the contragredient $\dual V$ of $V$.
	For a \(K^+\)-type \(\sigma \in \Irr(K^+)\), it is easy to verify that \(V_\sigma\) is paired with \((\dual V)_{\dual \sigma}\).
	Choose dual bases $\{f_1, \cdots, f_m\} \subset V_{\sigma{}}$ and $\{\dual f_1, \cdots, \dual f_m\} \subset (\dual V)_{\dual{\sigma{}}}$, and define the partial trace
\[\Tr_{\sigma}(g) = \sum_{r =1}^m \la \dual f_r, g \cdot f_r \ra\]
	for $g \in G^+$.
	One can check by definition that $\Theta(\pi^+) = \sum_\sigma \Tr_\sigma$, at least in the sense of distributions.
	In this subsection we estimate these $\Tr_\sigma$.

	Since $G = K P_\min = K \bar P$, and the compact group $K \cap \bar P$ has trivial intersection with $\bar A \bar N$, any element $g \in G$ can be expressed as $u(g) l(g) \left( \exp H(g) \right) n(g)$, where $H(g) \in \bar {\lag a}$ and $n(g) \in \bar N$ are uniquely determined by $g$.
	Hence, for \(\bar \nu \in \bar {\lag a}^*\), Harish-Chandra's spherical function is defined by
\[\varphi_{\bar \nu}(g) = \int_K e^{-\la \bar \nu + \rho_{\bar P}, H(g^{-1} k)\ra} \d k.\]
	The domain of $\varphi_{\bar \nu}$ extends naturally to $G^+= G\rtimes \la \tau \ra$ by $\varphi_{\bar\nu} (g\tau) = \varphi_{\bar\nu}(g)$.
	Now we present the following generalization of \cite[Lemma 6.23]{HS83}.

\begin{lem}\label{asymptotic-matrix-coefficient}
	For any compact subset $S \subset G$, there is a positive constant $C(S, \nu)$ such that for each $g\in S, x \in G^+$, and $\sigma{} \in \Irr(K^+)$
\[\left |\Tr_{\sigma{}} (g xg^{-1}) \right | \leqslant C(S, \nu) d_{\sigma{}}^3 \varphi_{\Re \nu}(x).\]
	Here $d_{\sigma{}}$ denotes the dimension of $\sigma{}$.
\end{lem}
\begin{pf}
	We reproduce the argument of \cite[Lemma 6.23]{HS83}.
	Denote the underlying $(\bar {\lag m}, \bar M \cap K)$-module of $\xi$ by $W$, and the underlying $({\lag g}, K)$-module of $\bar I$ by $E$.
	We have mentioned in Subsection \ref{subsec: parabolic-descent} that $E$ is admissible.
	Thus, for each $\sigma \in \Irr(K^+)$, $E_\sigma \supseteq V_\sigma$ is finite-dimensional.

	The contragredient of $\bar I$ is $\bar I^\vee = \Ind_{\bar P}^G(\xi^\vee)$, where the pairing is given by
\[\la f, \dual f \ra = \int_K \la f(k), \dual f (k) \ra \d k,\]
	and $(E^\vee)_{\sigma^\vee}$ is paired with $E_\sigma$.
	We have chosen a dual basis $\{ f_1, \cdots, f_m \} \subset V_\sigma$ and $\{ f_1^\vee, \cdots, f_m^\vee\} \subset (V^\vee)_{\sigma^\vee}$, where $V_\sigma \subseteq E_\sigma$, and $\{f_r^\vee\}$ can be lifted to $(E^\vee)_{\sigma^\vee}$ via the canonical splitting of $(E^\vee)_{\sigma^\vee} \twoheadrightarrow (V^\vee)_{\sigma^\vee}$ as $K$-representations.
	Then
\[\begin{aligned}
	\Tr_{\sigma{}}( gx g^{-1}) =	&
	\sum_{r=1}^m \la \dual f_r,  \bar I (g x g^{-1}) f_r \ra = 
	\sum_{r=1}^m \la \dual{\bar I}(g^{-1}) \dual f_r, \bar I (x g^{-1}) f_r \ra 	\\
	=	&
	\sum_{r=1}^m \int_K \la \left( \dual{\bar I} (g^{-1}) \dual f_r \right)(k), \left({\bar I} (x g^{-1}) f_r \right) (k) \ra \d k . 
	\end{aligned}\]
	Recall the definition of $\bar I(g)$ in \eqref{eq: para-ind} and $\tau_{\bar I}$ in Lemma \ref{lem: extend-standard}: for $ x \in G$, $f \in E$, $k \in K$,
\[\left( \bar I(x \tau) f \right)(k) = e^{-\la \nu + \rho_{\bar P}, H(x^{-1} k) \ra } \eta(l(x^{-1} k)) (\tau_{\bar I}f) (u(x^{-1}k)).\]
	Since the $\bar L$-representation $\eta$ is unitary, there exists a $\bar L$-invariant Hermitian inner product $\la -, - \ra_W$ on $W$.
	Define the $L^\infty$-norm on $E$ by
\[\| f \|_\infty = \sup_{u \in K} |f(u) |_W.\]
	Then for each $x \in G$, $f \in E$, and $k \in K$
\[\left | \left ( \bar I (x \tau ) f \right )(k) \right |_W \leqslant
	e^{-\la \Re \nu + \rho_{\bar P}, H(x^{-1} k) \ra} \|\tau_{\bar I} f\|_\infty =
	e^{-\la \Re \nu + \rho_{\bar P}, H(x^{-1} k) \ra} \| f\|_\infty.\]
	Hence, for any $x \in G^+ = G \rtimes \la \tau \ra$,
\[\begin{aligned}
	|\Tr_{\sigma{}}( gx g^{-1}) | 	\leqslant 	& 
	\sum_{r=1}^m \int_K   \left |\left( \dual {\bar I} (g^{-1}) \dual f_r \right) (k)\right |_{\dual W}
		\cdot  \left |\left( \bar I (x) \bar I (g^{-1}) f_r \right)(k)\right |_W \d k \\
	\leqslant 	& 
	\sum_{r=1}^m \left \| \dual {\bar I}(g^{-1}) \dual f_r \right \|_\infty \varphi_{\Re \nu}(x) \left\| {\bar I} (g^{-1})f_r \right \|_\infty
	\leqslant	
	C(S, \nu) \varphi_{\Re \nu}(x) \sum_{r=1}^m \|f_r \|_\infty \|\dual f_r \|_\infty.
	\end{aligned}\]

	The Hermitian inner product $\la \cdot, \cdot \ra_W$ induces a $K$-invariant inner product on $E$:
\[\la f, f' \ra_2= \int_K \la f(k), f'(k) \ra_W \d k.\]
	We can require $\{f_1, \cdots, f_m\} \subset E_\sigma$ to be orthonormal.
	Then, using the same argument as in \cite[pages 104--105]{HS83} for the proof of Lemma 6.23, one can deduce from Schur's orthogonality relations that $\|f_i\|_\infty \leqslant d_\sigma^{\frac 1 2}$.
	There is also a $K$-invariant inner product on $(E^\vee)_{\sigma^\vee}$, which induces a conjugate-linear isomorphism $E_\sigma \cong (E^\vee)_{\sigma^\vee}$ of Hilbert representations of $K$.
	In particular, the dual basis $\{f_1^\vee, \cdots f_m^\vee\} \subset (E^\vee)_{\sigma^\vee}$ of $\{f_1, \cdots, f_m\}$ is also orthonormal.
	Then we also have $\| f_i^\vee \|_\infty \leqslant d_\sigma^{\frac 1 2}$.
	Now we obtain $|\Tr_{\sigma{}}( gxg^{-1}) | \leqslant
	C(S, \nu) \varphi_{\Re \nu}(x)  m d_{\sigma{}}.$

	According to \cite[(6.31)]{HS83}, the dimension $m$ of $V_{\sigma{}}$ is at most $d_{\sigma{}}^2$, which can be deduced as follows.
	Due to Casselman's subrepresentation theorem (see \cite[Subsubsection 3.8.3]{Wallach88}), we can realize $V$ as a submodule of (the underlying Harish-Chandra module of) some $\Ind_{P_\min}^G(\xi_\min)$, where $\xi_\min$ is an irreducible representation of $M_\min = L_\min \times A_\min$.
	We know from \cite[Formula (11.42)]{KV95} that the multiplicity of any $K$-type in $\Ind_{P_\min}^G(\xi_\min)$ is no larger than the dimension of this $K$-type.
	Then the same bounded multiplicity holds for $V$, so $m \leqslant d_\sigma^2$.
	Now we conclude that  $|\Tr_{\sigma{}}( gxg^{-1}) | \leqslant
	C(S, \nu)  d_{\sigma{}}^3\varphi_{\Re \nu}(x)$ as desired.
\end{pf}

\begin{rmk}
	We use the unitarity of $\xi$ rather than temperedness.
	In particular, if $\pi$ is unitary, we can apply the above argument to $\pi$ in place of $\xi$, and obtain that $|\Tr_\sigma(g x g^{-1})| \leqslant C(S) d_\sigma^3.$
\end{rmk}

\subsection{Restriction to closed subsets}
	Recall that the task of this section is to estimate $\Delta_{G^+}\Theta(\pi^+)(s\tau \exp X)$ for $X \in \lag a_0^{\tau, -}$.

	It is not trivial to restrict the equality $\Theta(\pi^+) = \sum_{\sigma{}} \Tr_{\sigma{}}$ of distributions to regular points or closed subsets.
	We have to deal with the pointwise convergence of the summation over $\sigma\in \Irr(K^+)$.
	Consider the Casimir element $\Omega_K \in \uea{\lag k}$.
	It is $\tau$-invariant and acts on each $\sigma{}$ by a scalar $\omega_{\sigma{}}\geqslant 0$.
	This scalar is explicitly determined from the highest weight of $\sigma$, as in \cite[Proposition 5.28(b)]{Knapp96}.
	Moreover, the dimension $d_\sigma$ of $\sigma$ can also be determined from the highest weight by the Weyl dimension formula, see \cite[Theorem 5.84]{Knapp96}.
	From these two formulas, we can see that there exists an $n \in \N$ such that $\sum_\sigma \frac{d_{\sigma{}}^3}{(1+\omega_{\sigma{}})^n}$ converges absolutely for $\sigma$ running over $
\Irr(K^+)$.
	Consequently, the following equality also holds in the sense of distributions:
\[\Theta (\pi^+) = (1+ \Omega_K)^n \sum_{\sigma{}} \left( 1+ \omega_{\sigma{}} \right)^{-n} \Tr_{\sigma{}},\]
	where $(1+ \Omega_K)^n \in \uea{\lag k} \subseteq \uea{\lag g}$ is realized as a left-invariant differential operator $\Omega$ on $G^+$, and $\sum_{\sigma{}} \left( 1+ \omega_{\sigma{}} \right)^{-n} \Tr_{\sigma{}}$ now converges to a well-defined continuous function $h$ on $G^+$.
	Moreover, $h$ satisfies an estimate similar to that for $\Tr_{\sigma{}}$ in Lemma \ref{asymptotic-matrix-coefficient}: for each compact subset $S \subset G$, there is a positive constant $C(S, \nu)$ such that for each $g\in S$, $x \in G^+$,
\begin{equation}\label{estimate-on-h}
	|h(g x g^{-1})| \leqslant C(S, \nu) \varphi_{\Re \nu}(x).
	\end{equation}

	Let $s \in M^-$.
	For convenience, we introduce the notation $\Delta_G^{s\tau}$ and $\Theta^{s\tau}(\pi)$ by $\Delta_G^{s\tau}(t): = \Delta_{G^+}(ts\tau)$ and $\Theta^{s\tau}(\pi)(t): = \Theta(\pi^+)(ts\tau)$ for $t \in Z_s' = Z_s \cap G' s^{-1}$.
	The following lemma generalizes \cite[Lemma 6.39]{HS83}.

\begin{lem}\label{lem: asymptotic-restriction}
	As a distribution on $Z_s^-: = Z_s \cap M^-s^{-1}$,
\[\Delta_G^{s\tau} \Theta^{s\tau}(\pi) = \sum_j W_j \chi_j,\]
	where
\begin{itemize}
\item
	each $W_j$ is a translation-invariant differential operator on $Z_s$, and
\item
	each $\chi_j$ is a continuous function on $Z_s^-$, with the estimate that for each compact subset $U \subset Z_s^-$, there is $C(U, \nu)>0$ such that for each $t \in U, X\in {\lag a}^{\tau, -}_0$,
\begin{equation*}\label{estimate-on-chi_j}
	|\chi_j(t \exp X)| \leqslant C(U, \nu)  \left(\Delta_{G^{+}}\varphi_{\Re \nu} \right) (t s \tau \exp X).
\end{equation*}
\end{itemize}
\end{lem}
\begin{pf}
\begin{subequations}
	We repeat the argument from the proof of \cite[Lemma 6.39]{HS83}: pull back $\Theta(\pi^+) = \Omega h$ via the covering map $p: G/ T_s^{s\tau} \times Z_s' \to G[Z_s' s\tau]$.
	The main issue is to pull back the operator $\Omega$ and estimate the coefficients of the resulting operator; this will
be done in the next subsection.

	Fix $\gamma \in \Smt_\cpt(G/T_s^{s\tau})$ such that $\int_{G/T_s^{s\tau}} \gamma \d \bar g=1$.
	For any $\psi \in \Smt_\cpt\left( Z_s' \right)$, define
\[ F_\psi(\bar g, t) :=  |D_{\lag c_s}(s\tau)|^{-1} \gamma(\bar g) \left( |\Delta_G^{s\tau}|^{-1} \psi\right)(t), \]
and we have
\begin{equation}\begin{aligned}\label{eq: restrict-Theta-to-Z-1}
	\int_{Z_s}  \Delta_G^{s\tau} \Theta^{s\tau}(\pi)  \psi \d t=&
	\int_{Z_s} \Delta_G^{s\tau} \psi\d t
	\int_{G/T_s^{s\tau}}\Theta(\pi^+)(\Ad(\bar g)( t s\tau))\gamma(\bar g)\d \bar g\\
	=&
	\int_{Z_s}  |D_{{\lag q}_s}^{s\tau}| \d t
	\int_{G/T_s^{s\tau}}\left( p^*\Theta(\pi^+)\cdot  F_\psi\right)(\bar g, t) \d \bar g 
	\stackrel{\eqref{eq: pull-back-int}}=
	\int_{G\tau} \Theta(\pi^+) \cdot p_* F_\psi \d g.
	\end{aligned}\end{equation}
	Since $\Theta(\pi^+) = \Omega h$, where $\Omega = (1+\Omega_K)^n$ is a self-adjoint differential operator on $G\tau$ and $h$ is continuous on $G\tau$,
\begin{equation}\begin{aligned}\label{eq: restrict-Theta-to-Z-2}
	\int_{G\tau} \Theta(\pi^+) \cdot p_* F_\psi \d g =&
	\int_{G\tau} \Omega h \cdot p_* F_\psi \d g =
	\int_{G\tau} h \cdot \Omega (p_*F_\psi) \d g \\
	\stackrel{\textrm{Lemma \ref{lem: differential-geometry}(2)}}=& 
	\int_{G\tau} h \cdot p_*\left( (p^*\Omega) F_\psi \right) \d g 
	\stackrel{\eqref{eq: pull-back-int}}=
	\int_{Z_s} |D_{{\lag q}_s}^{s\tau}| \d t
	\int_{G/T_s^{s\tau}} p^* h \cdot  (p^*\Omega) F_\psi \cdot  \d \bar g.
	\end{aligned}\end{equation}
	Decompose the differential operator $p^*\Omega$ on $G/ T_s^{s\tau} \times Z_s'$ as
\[p^*\Omega = \sum_r H_r Z_r \boxtimes Y_r,\]
	where each $H_r $ is a smooth function, each $Z_r$ is a differential operator on $G/T_s^{s\tau}$, and each $Y_r$ is a translation-invariant differential operator on $Z_s$.
	Then
\begin{equation}\begin{aligned}\label{eq: restrict-Theta-to-Z-3}
	\int_{G/T_s^{s\tau}} p^* h \cdot  (p^*\Omega) F_\psi\d \bar g 	
	=&
	\int_{G/T_s^{s\tau}} p^* h \cdot \sum_r  H_r \cdot  (Z_r \boxtimes Y_r) 
	\left(|D_{\lag c_s}(s\tau)|^{-1} \gamma \boxtimes (|\Delta_G^{s\tau}|^{-1} \psi)\right) \d \bar g\\
	 =&
	|D_{\lag c_s}(s\tau)|^{-1} \sum_r Y_r(|\Delta_G^{s\tau}|^{-1} \psi) \cdot
	\int_{G/T_s^{s\tau}} p^*h \cdot H_r \cdot Z_r \gamma \d \bar g.
	\end{aligned}\end{equation}
	Write
\[|\Delta_G^{s\tau}| Y_r |\Delta_G^{s\tau}|^{-1} = \sum_{j} \zeta_{r, j} X_j,\]
	where each $\zeta_{r, j}$ is a smooth function on $Z_s^-$ and each $X_j$ is a translation-invariant differential operator on $Z_s$.
  The above calculation shows that
\begin{equation}\begin{aligned}\label{eq: restrict-Theta-to-Z-4}
	\int_{Z_s} \Delta_G^{s\tau}\Theta^{s\tau}(\pi) \psi \d t \stackrel{\eqref{eq: restrict-Theta-to-Z-1}} =	&
	\int_{G\tau} \Theta(\pi^+) p_*F_\psi \d g \stackrel{\eqref{eq: restrict-Theta-to-Z-2}} = 
	\int_{Z_s} |D_{{\lag q}_s}^{s\tau}| \d t
	\int_{G/T_s^{s\tau}} p^* h \cdot  (p^*\Omega) F_\psi \cdot  \d \bar g	\\
	\stackrel{\eqref{eq: restrict-Theta-to-Z-3}} =	&
	\int_{Z_s} (\Delta_G^{s\tau})^2  \sum_r Y_r(|\Delta_G^{s\tau}|^{-1} \psi) \int_{G/T_s^{s\tau}} p^* h \cdot H_r \cdot Z_r \gamma \d \bar g 	\\
	=	&
	\int_{Z_s} \Delta_G^{s\tau} \sum_{r, j}\zeta_{r, j} X_j \psi \cdot \omega_r\d t,
	\end{aligned}\end{equation}
	where
\[\omega_r(t) = \int_{G/T_s^{s\tau}}\left(p^*h \cdot H_r \right)(\bar g, t) \cdot \left(Z_r \gamma\right)(\bar g)  \d \bar g.\]
Now let
\[\chi_j = \Delta_G^{s\tau} \sum_r \zeta_{r, j} \omega_r = 
	 \Delta_G^{s\tau} \sum_r \zeta_{r, j} \cdot
	 	\int_{G/T_s^{s\tau}}p^*h \cdot H_r  \cdot Z_r \gamma  \d \bar g ,\]
	and $W_j$ be the adjoint operator of $X_j$; then \eqref{eq: restrict-Theta-to-Z-4} becomes
\begin{equation}\label{eq: restrict-Theta-to-Z-5}
	\int_{Z_s} \Delta_G^{s\tau}\Theta^{s\tau}(\pi) \psi \d t \stackrel{\eqref{eq: restrict-Theta-to-Z-4}} = 
	\int_{Z_s} \sum_j \chi_j \cdot X_j \psi \d t =
	\int_{Z_s} \psi \sum_j W_j \chi_j \d t.
	\end{equation}
	Consequently, we obtain the desired equality of distributions $\Delta_G^{s\tau} \Theta^{s\tau}(\pi) = \sum_j W_j \chi_j$.

	It is clear that each $W_j$ is also translation-invariant.
	The estimate for each $\chi_j$ follows from the compact support of $\gamma$, the estimate \eqref{estimate-on-h} for $h$, and the estimates for $H_r$ and $\zeta_{r,j}$ in Lemma \ref{lem: play-with-differential-operator}.
\end{subequations}
\end{pf}

\subsection{Analysis on homogeneous space}
	In this subsection, we deal with the decompositions and estimates of differential operators claimed in the proof of Lemma \ref{lem: asymptotic-restriction}.
	Its proof is similar to the last paragraph in \cite[proof of Lemma 6.39]{HS83}, by replacing $\xi$ with $p: G/T_s^{s\tau} \times Z_s' \to G'\tau$.
\begin{lem}\label{lem: play-with-differential-operator}\label{lem: differential-geometry}
\begin{enumerate}[fullwidth, label=$(\arabic*)$]
\item
	For any translation-invariant differential operator $Y$ on $Z_s$, there exists a finite decomposition
\[|\Delta_G^{s\tau}| Y |\Delta_G^{s\tau}|^{-1} = \sum_{j} \zeta_j X_j\]
	such that
\begin{itemize}
\item
	each $X_j$ is a translation-invariant differential operator on $Z_s$, and
\item
	each $\zeta_j$ is a smooth function on $Z_s^-$, with the estimate that for each compact subset $U \subset Z_s^-$, there is $C(U) >0$ such that for each $t \in U$, $X \in {\lag a}^{\tau, -}_0$,
\[|\zeta_j(t \exp X)| \leqslant C(U).\]
\end{itemize}

\item
	For any left translation-invariant differential operator $\Omega$ on $G\tau$, there exists a pullback $p^* \Omega$ on $G/T_s^{s\tau} \times Z_s'$ such that for each $F \in  \Smt_\cpt \left( G/T_s^{s\tau} \times Z_s' \right)$,
\[p_*\left( (p^*\Omega) F \right) = \Omega( p_* F).\]
	Moreover, there is a finite decomposition
\[p^*\Omega = \sum_r H_r Z_r \boxtimes Y_r\]
	such that
\begin{itemize}
\item
	each $Y_r$ is a translation-invariant differential operator on $Z_s$,
\item
	each $Z_r$ is a smooth differential operator on $G/T_s^{s\tau}$, and
\item
	each $H_r$ is a smooth function on $G/T_s^{s\tau} \times Z_s^-$, with the estimate that for each compact subset $S \subset G/T_s^{s\tau}$ and compact subset $U \subset Z_s^-$, there is $C(S, U) >0$ such that for each $\bar g \in S$, $t \in U$, $X \in {\lag a}^{\tau, -}_0$,
\[|H_r(\bar g, t \exp X)| \leqslant C(S, U).\]
\end{itemize}
\end{enumerate}
\end{lem}
\begin{pf}
	We mainly focus on the case of vector fields, since differential operators are tensor products of them, and the conclusion follows by induction.

	(1) Any $Y \in {\lag z}_{s,0}$ gives rise to a translation-invariant vector field on $Z_s$, which sends $\psi \in \Smt\left( Z_s \right)$ to
\[(Y \psi)(t) = \left. \frac \d {\d \varepsilon} \right|_{\varepsilon =0} \psi( t \exp \varepsilon Y).\]
	Then, as a differential operator on $Z_s^-$, we have $|\Delta_G^{s\tau}| Y |\Delta_G^{s\tau}|^{-1} = Y + \phi$, where $\phi := |\Delta_G^{s\tau}| Y(|\Delta_G^{s\tau}|^{-1})$ is a smooth function on $Z_s^-$.
	It suffices to bound the derivatives of $\phi$.

	Recall from the formula \eqref{eq: twisted-Weyl-denominator} that
\[\Delta_G^{s\tau}(t) = \prod_{\alpha_\res \in \RS_\res^+} \prod_{j=1}^{\dim {\lag g}_{\alpha_\res}}
	\left| \lambda_{\alpha_\res, j} \alpha_\res(t ) \right|^{-\frac 1 2}
	\left| 1- \lambda_{\alpha_\res, j}\alpha_\res(t ) \right|.\]
	By Leibniz's rule we calculate
\[\begin{aligned}
	\phi=	&
	 \sum_{\alpha_\res \in \RS_\res^+} \sum_{j=1}^{\dim {\lag g}_{\alpha_\res}}
	\left| \alpha_\res(t) \right|^{-\frac 1 2} Y\left( 
		|\alpha_\res(t)|^{\frac 1 2}
	\right )+
	\left| 1- \lambda_{\alpha_\res, j}\alpha_\res(t ) \right| Y\left(
		\left| 1- \lambda_{\alpha_\res, j}\alpha_\res(t ) \right|^{-1}
	\right)	\\
	=	&
	\sum_{\alpha_\res \in \RS_\res^+} \sum_{j=1}^{\dim {\lag g}_{\alpha_\res}}
	\frac 1 2 \Re \la \alpha_\res, Y \ra-
	 \Re \frac{\lambda_{\alpha_\res, j} \la \alpha_\res, Y \ra}
	{\lambda_{\alpha_\res, j} \alpha_\res(t ) - 1}.
	\end{aligned}\]
	Consequently, it suffices to bound the derivatives of $ \left(\lambda_{\alpha_\res, j} \alpha_\res(t ) - 1\right)^{-1} $ for $t\in U A^{\tau, -}$.

	We consider two cases.
	Suppose $\alpha_\res \in \RS^+({\lag g}, {\lag t}^\tau)$ does not vanish on ${\lag a}^\tau$; for each $t \in U \subset Z_s^-$, $|\lambda_{\alpha_\res, j} \alpha_\res(t )|< 1$.
	Since $U$ is compact, there exists $\varepsilon(U)>0$ such that for each $t \in U$, $| \lambda_{\alpha_\res, j}\alpha_\res(t )| \leqslant  1- \varepsilon(U)$.
	Then, based on the definition of $A^{\tau, -}$ earlier in this section, $| \lambda_{\alpha_\res, j}\alpha_\res(t)| \leqslant 1- \varepsilon(U)$ for all $t\in U A^{\tau, -}$.
	Suppose $\alpha_\res \in \RS^+({\lag g}, {\lag t}^\tau)$ vanishes on ${\lag a}^\tau$; then, for each $t\in U$, $a\in A^{\tau, -}$, $\alpha_\res(ta) = \alpha_\res(t)$.
	Since for each $t\in U \subset Z_s^-$, $ts\tau$ is regular by definition,  its eigenvalues $\lambda_{\alpha_\res, j} \alpha_\res(t)$ on ${\lag g}_{\alpha_\res}$ would never be $1$.
	Due to the compactness of $U$, there exists $\varepsilon(U) >0$ such that for each $t \in U$,  $|\lambda_{\alpha_\res, j} \alpha_\res(t) -1| \geqslant \varepsilon(U)$.
	In either case, $|\lambda_{\alpha_\res, j} \alpha_\res(t) - 1 | ^{-1}$ is bounded by some $\varepsilon(U)^{-1}$ for all $t\in U A^{\tau, -}$.
	It follows that $\phi$ is bounded on $UA^{\tau, -}$; the same argument applies to all its derivatives.

	(2) For any vector field $X$ on $G\tau$, define $p^*X$ as follows.
	Take a partition of unity $\{(\varphi_i, U_i)\}$ for $G/T_s^{s\tau} \times Z_s'$ such that each $p|_{U_i}$ is a diffeomorphism.
	Define $(p^*X): F \mapsto \sum_i \varphi_i \cdot \left( (p|_{U_i})_*^{-1} X\right) ( F|_{U_i})$.
	Indeed, if $p$ is diffeomorphic on $U \subset G/T_s^{s\tau} \times Z_s'$, then one can deduce $p^* X|_U = (p|_U)_*^{-1} X$ from $p_*((p^*X) F) = X(p_* F)$.

	Let us calculate $p^* X$ more explicitly.
	Take a small neighborhood $V$ of $1\in G/T_s^{s\tau}$ such that
\begin{enumerate}[label=(\alph*)]
\item
	the tangent bundle $\Tg \left( G/T_s^{s\tau} \right)$ is trivial over $V$,
\item
	the submersion $G \to G/T_s^{s\tau}$ has a section $\sigma$ over $V$, and
\item
	the translation $Vw$ of $V$ by $w \in W_{G, s\tau}$ is disjoint from $V$ unless $w=1$.
\end{enumerate}
	It follows from (c) that for each $g_0 \in G$, the map $p:g_0 V \times Z_s'\to p\left(g_0V \times Z_s'\right)$ is a diffeomorphism.

	Any $Y\in {\lag z}_{s, 0}$ yields a translation-invariant vector field on $Z_s$ as defined in (1), and any $Z \in {\lag q}_{s, 0}$ becomes a vector field on $g_0 V$ by sending $\gamma \in \Smt(g_0 V)$ to
\[(Z\gamma)(g_0 \bar g) = \left. \frac{\d}{\d \varepsilon} \right|_{\varepsilon = 0}\gamma ( g_0 \sigma(\bar g) \exp \varepsilon Z),\]
	where
$\gamma$ is pulled back to $G$ via $G \to G/T_s^{s\tau}$.
	In particular, if we take orthonormal bases $\{Y_i\} \subset {\lag z}_{s, 0}$ and $\{Z_j\} \subset {\lag q}_{s, 0}$, then they give rise to a local frame for $\Tg\left(G/ T_s^{s\tau} \times Z_s' \right)$ over $g_0 V \times Z_s'$.

	For a vector field $X$ on $G\tau$, we know from (c) that $p^*X |_{g_0 V \times Z_s'} = (p|_{g_0 V \times Z_s'})^{-1}_* X$.
	In particular, decompose it with respect to the local frame $\{Y_i, Z_j\}$ for $\Tg \left( G/ T_s^{s\tau} \times Z_s' \right)$ over $g_0V \times Z_s'$, say
\[\left. p^*X \right|_{g_0 V \times Z_s'} = \sum_i E_i Y_i + \sum_j F_j Z_j.\]
	Then for each $(\bar g, t) \in V \times Z_s'$,
\[\begin{aligned}
	X_{p(g_0 \bar g, t) } =	&
	\sum_i E_i(g_0 \bar g, t) p_* Y_{i, (g_0 \bar g, t)} +
	\sum_j F_j(g_0 \bar g, t) p_* Z_{j, (g_0\bar g, t)}	\\
	=	&
	\Ad(g_0 \sigma(\bar g)) \left( 
		\sum_i E_i(g_0 \bar g, t)  Y_i +
		\sum_j F_j(g_0 \bar g, t) \left( \Ad(ts \tau)^{-1} -1  \right) Z_j \right) \in {\lag g}_0.
	\end{aligned}\]
	Here we have used the calculation in the proof of Lemma \ref{lem: finite-cover}.

	If $X$ is left translation invariant, and is given by
\[(Xf)(g\tau) = \left. \frac{\d}{\d \varepsilon} \right|_{\varepsilon = 0}
	f(g \tau \exp \varepsilon X)\]
	for some $X \in {\lag g}_0$, then
\[\begin{aligned}
	E_i(g_0 \bar g, t) =	&
	\la Y_i,  \Ad(g_0 \sigma(\bar g))^{-1} X\ra,
	\\
	F_j(g_0 \bar g, t) =	&\la Z_j,  \left( 
		\left(\Ad(t s \tau)^{-1} - 1\right)^{-1} \circ 
		\pr_{{\lag q}_{s, 0}} \circ 
		\Ad(g_0 \sigma(\bar g) )^{-1} \right) X \ra.
	\end{aligned}\]
	Here one should observe that $\Ad(ts \tau)^{-1} -1$ is invertible on $\lag q_s \subset \lag g$.

	To obtain the estimate, it suffices to bound $F_j(g_0 \bar g, t)$ and its derivatives for $(\bar g, t) \in V \times U A^{\tau, -}$.
	According to the above expression for $F_j$, what matters is the eigenvalue of $\left( \Ad(t s \tau)^{-1} - 1\right)^{-1}$, which takes the form $\left( \left( \lambda_{\alpha_\res, j} \alpha_\res(t )\right) ^{-1}  - 1\right)^{-1}$ on ${\lag g}_{\alpha_\res}$, and  $\left( \lambda_{\alpha_\res, j} \alpha_\res(t ) - 1\right)^{-1}$ on ${\lag g}_{-\alpha_\res}$.
	We have shown in (1) that $|\lambda_{\alpha_\res, j} \alpha_\res(t) - 1 | ^{-1}$ is bounded by some $\varepsilon(U)^{-1}$ for all $t\in U A^{\tau, -}$.
	Since $|\lambda_{\alpha_\res, j}\alpha_\res(t)|$ is also bounded on $UA^{\tau, -}$, due to the compactness of $U$ and the definition of $A^{\tau, -}$,
\[\left| 	\left( \lambda_{\alpha_\res, j} \alpha_\res(t) \right)^{-1}
		- 1 \right| ^{-1} =
	|\lambda_{\alpha_\res, j} \alpha_\res(t) | \cdot 
	|\lambda_{\alpha_\res, j} \alpha_\res(t) - 1 | ^{-1}\]
	is bounded for $t\in UA^{\tau, -}$.
	The derivatives of $\left( \lambda_{\alpha_\res, j} \alpha_\res(t ) - 1\right)^{-1}$ and $\left( \left( \lambda_{\alpha_\res, j} \alpha_\res(t ) \right)^{-1} - 1\right)^{-1}$ can then be bounded using a similar argument.
	Now the conclusion follows.
\end{pf}

\subsection{One more condition on $\lambda$}\label{subsec: extra-condition}
	We have fixed a Langlands embedding $V \hookrightarrow \Ind_{\bar P}^G( \xi \boxtimes \nu)$, where $\bar P = \bar L \bar A \bar N$, and $\nu \in \bar {\lag a}^*$ satisfies $\Re \la \nu, \beta \ra < 0$ for all $\beta \in \RS^+({\lag g}, \bar {\lag a}) = \RS(\bar {\lag n}, \bar {\lag a})$.
	We want a stronger dominance condition on $\nu$:
\begin{equation}\label{eq: nu-dominant}
	\forall \beta \in \RS^+({\lag g}, \bar {\lag a}), \quad
	\Re\la \nu+ \rho_{\bar P}, \beta \ra <0.
	\end{equation}
	Assuming the above property, we can deduce the following analogue of \cite[Lemma 6.46]{HS83}:

\begin{lem}\label{asymptotic-character}\label{lem: exponents-of-V}
	Up to conjugation by $G$, assume the minimal parabolic subgroup $P_\min \subseteq \bar P$ is also contained in $P$, with $L_\min A_\min\subseteq M$, $A_\min \supseteq A$, $N_\min \supseteq N$.
	In particular, functionals on $\bar {\lag a}$, ${\lag a}$ can be extended to ${\lag a}_\min$ by zero to their orthogonal complements.

\begin{enumerate}[fullwidth, label=$(\arabic*)$]
\item
	For any compact subset $U \subset Z_s^-$, there exists a constant $C(U, \nu)>0$ such that for each $t\in U$ and $X\in {\lag a}^{\tau, -}_0$,
\[\left( \Delta_{G^+} \varphi_{\Re \nu} \right)(t s\tau \exp X) \leqslant 
	C(U, \nu) e^{\Re \la \nu+ \rho_{\bar P} - \rho _P, X \ra}.\]

\item
	Since $H_0({\lag n}, V)_\mu = 0$, there exists $\delta>0$ and $\beta_\res \in \RS^+({\lag g}, {\lag a}^\tau) = \RS({\lag n}, {\lag a}^\tau)$ such that for each $X \in {\lag a}^{\tau, -}_0$,
\[\Re \la \nu+ \rho_{\bar P}, X \ra \leqslant \Re \la \mu+ \rho_P+ \delta \beta_\res, X \ra.\]
\end{enumerate}
\end{lem}
\begin{pf}
(1) One calculates from \eqref{eq: descent} and $\Delta_{G^+} = |\det_{\lag n}|^{-\frac 1 2} | D_{\lag n} | \Delta_{M^+}$ that
\[\Delta_{G^+}(ts\tau \exp X) =  e^{-\la \rho_P, X \ra} 
	\left( |\det_{\lag n}|^{-\frac 1 2} \Delta_{M^+} \right) (ts\tau) \cdot
	|D_{\lag n}^{s\tau}(t\exp X)|.\]
	Since $t s \exp X \in U s A^{\tau, -} \subset M^-$, the factors $\left( |\det_{\lag n}|^{-\frac 1 2} \Delta_{M^+} \right) (ts\tau)$ and $|D_{\lag n}^{s\tau}(t\exp X)|$ are bounded.
	Then
\[\Delta_{G^+}(ts\tau \exp X)  \leqslant 
	C(U) e^{-\la \rho_P, X \ra}.\]

	As for $\varphi_{\Re \nu}(ts\tau a) = 
	\varphi_{\Re \nu}(ts a),$ with $t \in U$, $a \in A^{\tau, -}$, according to the proof of \cite[Lemma 6.46]{HS83},
\begin{enumerate}[label=(\alph*)]
\item
	from the compactness of $s^{-1}U^{-1}K$ one can deduce, as in \cite[(6.49)]{HS83}, that $\varphi_{\Re \nu}(tsa) \leqslant C(U, \nu) \varphi_{\Re \nu}(a)$, and
\item
	by manipulating with finite-dimensional representations one can deduce the inequality \cite[(6.50)]{HS83}, $e^{-\la \Re \nu + \rho_{\bar P}, H(a^{-1} k) \ra}\leqslant 
	e^{\la \Re \nu + \rho_{\bar P}, H(a) \ra}$,  under the condition \eqref{eq: nu-dominant} on $\Re \nu$.
\end{enumerate}
	It follows that
\[\varphi_{\Re \nu} (ts\tau \exp X) \stackrel{\textrm{(a)}}\leqslant
	C(U, \nu) \varphi_{\Re \nu}(\exp X)  = 
	C(U, \nu) \int_K e^{-\la \Re \nu + \rho_{\bar P}, H(\exp(-X) \cdot  k) \ra} \d k \stackrel{\textrm{(b)}} \leqslant
	C(U, \nu) e^{\la \Re \nu + \rho_{\bar P}, X \ra} \]
	for $t \in U$, $X \in {\lag a}_0^{\tau, -}$.
	Then the conclusion follows.

	(2) Extend $\nu, \rho_{\bar P} \in \bar{\lag a}^*$ and $\mu, \rho_P \in \lag a^*$ by zero to $\lag t$, respectively, on the orthogonal complements of $\bar{\lag a}$ and $\lag a$.
  Based on $H_0(\lag n, V)_\mu=0$, one can deduce 
\[\nu + \rho_{\bar P} - \mu - \rho_P = 
	\sum_{\alpha \in \RS^+({\lag g}, {\lag t})} x_\alpha \alpha,\]
	where
for each $\alpha$, if $\alpha|_{\lag a} \not = 0$, then $\Re x_\alpha \geqslant 0$, and
$\exists \alpha$ such that $\alpha|_{\lag a} \not = 0$ and $\Re x_\alpha > 0$;
see \cite[(6.55)-(6.59)]{HS83}.
	Elements in $\RS^+({\lag g}, {\lag a}) = \RS({\lag n}, {\lag a})$ come from non-zero restrictions of $\RS^+({\lag g}, {\lag t})$, and their restrictions to ${\lag a}^\tau$ do not vanish since ${\lag n}$ is stable under $\tau$.
	It follows that
\[ \left( \nu + \rho_{\bar P} - \mu - \rho_P \right) |_{{\lag a}^\tau} = 
	\sum_{\beta_\res \in \RS^+({\lag g}, {\lag a}^\tau)} x_{\beta_\res} \beta_\res\]
	with
$\Re x_{\beta_\res} \geqslant 0$ for each $\beta_\res \in \RS^+({\lag g}, {\lag a}^\tau)$, and
$\Re x_{\beta_\res} > 0$ for at least one $\beta_\res \in \RS^+({\lag g}, {\lag a}^\tau)$.
The positive $\Re x_{\beta_\res}$ gives rise to the desired $\delta$ and $\beta_\res$.
\end{pf}

	Let us explain how \eqref{eq: nu-dominant} is derived from an extra condition on $\lambda \in {\lag t}^*$.
	Recall the Borel pair $(\bar {\lag b}, \bar {\lag t})$ before Lemma \ref{lem: nu-stable-tau}.
	Extend $\beta \in \RS^+({\lag g}, \bar {\lag a})$ to $\bar {\lag t}$ by zero on the orthogonal complements of $\bar {\lag a}$, and transfer it to $\mu_\beta \in {\lag t}^*$ by $\Ad(g): \bar {\lag t} \to {\lag t}$.
	It follows from Lemma \ref{lem: nu-stable-tau} that
\[\la \nu, \beta \ra = \la (w\lambda) \circ \Ad(g) | _{\bar {\lag a}},\beta \ra =
	\la w\lambda, \mu_\beta \ra\]
	for some $w\in W({\lag g}, {\lag t})^\tau$.
	Since $w\lambda$ is stable under $\tau$, we deduce further that
\[\la \nu, \beta \ra = \la w\lambda, \mu_\beta \ra = 
	\la w\lambda, \Av(\mu_\beta) \ra, \]
	where $\Av(\mu_\beta)$ is the average of $\tau^i(\mu_\beta)$, $i = 1, \cdots, d$.

	One deduces from Fact \ref{fact: centralizer-of-A}(2) that the functional $\mu_\beta \in {\lag t}^*$ is a non-negative combination of $\Pi({\lag g}, {\lag t})$: since $(\bar {\lag b}, \bar {\lag t})$ is compatible with $\bar P = \bar M \bar N$, we have $\bar {\lag b} \supseteq \bar {\lag n}$, and it follows that any $\beta \in \RS^+({\lag g}, \bar {\lag a})$ takes the form $\alpha|_{\bar {\lag a}}$ for some $\alpha \in \RS^+({\lag g}, \bar {\lag t})$, so the extension $\pr_{\bar{\lag a}} \alpha \in \bar{\lag t}^*$ of $\beta \in \bar{\lag a}^*$ is a non-negative combination of $\Pi({\lag g}, \bar {\lag t})$, and then the conclusion for $\mu_\beta = \left( \pr_{\bar {\lag a}} \alpha \right) \circ \Ad(g)^{-1}$ follows.

	Since $\Pi({\lag g}, {\lag t})$ is stable under $\tau$, $\Av(\mu_\beta)$, as a non-negative combination of $\Pi({\lag g}, {\lag t})$, never vanishes.
	It is reasonable to require
\begin{equation}\label{eq: extra-on-lambda}
	\left| \Re \la w\lambda, \mu_\beta \ra \right| > C,
	\textrm{ for all } w \in W({\lag g}, {\lag t})^\tau, \beta \textrm{ appears as above};
	\end{equation}
	note that there are only finitely many standard parabolic subgroups $\bar P = \bar L \bar A \bar N \supseteq P_\min = L_\min A_\min N_\min$, and each $\RS^+(\lag g, \bar{\lag a}) = \RS(\bar{\lag n}, \bar{\lag a})$ is also finite.
	Consequently, the conditions in \eqref{eq: extra-on-lambda} involve only finitely many non-zero linear functionals on $(\lag t^*)^\tau$. 
  Since $\Lambda$ is a full lattice in $(\lag t^*)^\tau$, we may choose a strictly anti-dominant $\lambda_{\mathrm a} \in \Lambda$ outside their kernels.
  Then, for every $\lambda_0\in(\lag t^*)^\tau$, the weights $\lambda_0+N\lambda_{\mathrm a}$, for all sufficiently large integers $N$, are very
anti-dominant and satisfy \eqref{eq: extra-on-lambda}.
	In other words, let $\mathcal F \subset (\lag t^*)^\tau$ consist of the very anti-dominant weights satisfying \eqref{eq: extra-on-lambda}.
	We have shown that $\mathcal F + \Lambda = (\lag t^*)^\tau$.
	For such a $\lambda \in \mathcal F$, \eqref{eq: extra-on-lambda} implies that $\Re \la \nu, \beta \ra  = \Re \la w\lambda, \mu_\beta \ra <-C$, since we have required $\Re \la \nu, \beta \ra = \la \Re \nu, \beta \ra <0$.
	When $C$ is large enough, \eqref{eq: nu-dominant}, $\Re\la \nu+ \rho_{\bar P }, \beta \ra <0$, follows.

\subsection{Extraction of main part}\label{subsubsec: extract}
	In this final subsection, we finish the
\begin{pf}[ of Proposition \ref{prop: vanish}]

	Recall what we have done at the beginning of this section.
	Assume $\Theta(\pi^+)_\mu$ does not vanish at some $s\tau \in M^-\tau$; then there exist
\begin{itemize}
\item
	a neighborhood $U$ of $1 \in Z_s^-$,
\item
	a sequence $\{a_k \}\subset A^{\tau, -}$ such that for every $\beta_\res \in \RS^+({\lag g}, {\lag a}^\tau)$, ${\beta_\res}(a_k) \to 0$ as $k \to \infty$, and
\item
	a function $\psi \in \Smt_\cpt(U)$ such that
\[\Psi(\exp X):= e^{-\la \mu, X \ra}
	\int_U \left(\Delta_{G^+} \Theta(\pi^+) \right)(ts\tau \exp X) \psi(t) \d t\]
	has $\Psi(a_k) \to 1$.
\end{itemize}

	On the other hand, since $H_0({\lag n}, V)_\mu=0$, it follows from Lemmas \ref{lem: asymptotic-restriction} and \ref{asymptotic-character} that $\Delta_G^{s\tau}\Theta^{s\tau}(\pi) |_{Z_s^-} = \sum_j W_j \chi_j$, where
\begin{itemize}
\item
	each $W_j$ is a translation-invariant differential operator on $Z_s$, and
\item
	each $\chi_j$ satisfies the estimate
\[|\chi_j(t \exp X)| \leqslant C(\nu) e^{\la \Re \mu + \delta \beta_\res, X\ra}\]
	for $t \in \supp \psi$ and $X\in {\lag a}^{\tau, -}_0$.
\end{itemize}
	This would result in $\Psi(a_k) \to 0$.
	Indeed, since the adjoint $X_j$ of $W_j$ is also translation-invariant, we have
\[\Psi(\exp X) = e^{-\la \mu, X \ra} \int_U \left(\Delta_G^{s\tau} \Theta^{s\tau}(\pi) \right)(t \exp X) \psi(t) \d t =
	e^{-\la \mu, X \ra} \int_U \sum_j\chi_j (t \exp X) X_j\psi(t) \d t,\]
	and then
\[\left | \Psi(\exp X) \right| \leqslant
	e^{-\la\Re \mu, X \ra} \int_U \sum_j C(\nu)
	e^{\la \Re\mu + \delta \beta_\res, X \ra}|X_j\psi(t)| \d t \leqslant
	C e^{\delta \la \beta_\res, X \ra}.\]
	Since $\beta_\res(a_k)\to 0$, it follows that $\Psi(a_k) \to 0$, which leads to a contradiction.

	Consequently, $H_0({\lag n}, V)_\mu=0$, together with \eqref{eq: extra-on-lambda}, implies that $\Theta(\pi^+)_\mu$ vanishes on $M^-\tau$.
\end{pf}
 
\bibliography{D:/LaTeX/Research/Reference/A_packet_construction_p_adic.bib,
	D:/LaTeX/Research/Reference/A_packet_definition.bib, 
	D:/LaTeX/Research/Reference/Endoscopy_stablization.bib, 
	D:/LaTeX/Research/Reference/Personal.bib,
	D:/LaTeX/Research/Reference/Real_Lie_group.bib, 
	D:/LaTeX/Research/Reference/Rep_of_p_adic.bib, 
	D:/LaTeX/Research/Reference/Rep_of_real_Lie_group.bib}

\vspace{1em}
\begin{flushleft} \small
	Chang Huang: Institue for Advanced Study in Mathematics, Zhejiang University, Xihu District, Hangzhou 310058, China. \\
	E-mail address: \href{mailto:uhanghcang@163.com}{\texttt{uhanghcang@163.com}}
\end{flushleft}
\end{document}